\documentclass[12pt]{amsart} 

\usepackage{graphicx,amssymb,colordvi,textcomp,latexsym} 

\usepackage{tikz}
\usepackage{ulem} 

\usepackage[title]{appendix}
\usepackage[showframe=false]{geometry}
\usepackage{changepage}

\usepackage{adjustbox}
\usetikzlibrary{arrows,shapes, positioning, matrix, patterns, decorations.pathmorphing, arrows.meta,automata}

\usepackage{color}
\usepackage{colortbl}
\usepackage{soul}

\usepackage{array}
\usepackage[hidelinks]{hyperref}
\usepackage{amsfonts}
\usepackage{amsmath}
\usepackage{amsthm}

\usepackage{todonotes}
\usepackage{enumitem}   

\newcommand{\ignore}[1]{}

\newcommand{\R}{\mbox{$\mathbb{R}$}}

\newcommand{\TT}{\mathcal{T}}

\newcommand{\beqn}{\begin{eqnarray*}}
\newcommand{\eeqn}{\end{eqnarray*}}

\newtheorem{lemma}{Lemma}[section]

\newtheorem{prop}[lemma]{Proposition}

\theoremstyle{definition}
\newtheorem{Def}[lemma]{Definition}
\newtheorem{exam}[lemma]{Example}
\newtheorem{exams}[lemma]{Examples}

\theoremstyle{remark}
\newtheorem{rem}[lemma]{Remark}
\newtheorem{rems}[lemma]{Remarks}

\title[Excitatory-Inhibitory Networks]{Classification
of 3-node Restricted Excitatory-Inhibitory Networks}
\author{Manuela Aguiar}
\address{Manuela Aguiar, Centro de Matem\'atica da Universidade do Porto (CMUP), Faculdade de Ci\^encias, Universidade do Porto, Rua do Campo Alegre s/n, 4169-007 Porto, Portugal\newline 
Faculdade de Economia, Universidade do Porto, Rua Dr Roberto Frias, 4200-464 Porto, Portugal}
\email{maguiar@fep.up.pt}
\author{Ana Dias}
\address{Ana Dias, Centro de Matem\'atica da Universidade do Porto (CMUP), Departamento de Matem\'atica, Faculdade de Ci\^encias, 
Universidade do Porto, Rua do Campo Alegre s/n, 4169-007 Porto, Portugal}
\email{apdias@fc.up.pt}
\author{Ian Stewart}
\address{Ian Stewart, Mathematics Institute, University of Warwick, Coventry CV4 7AL, United Kingdom}
\email{i.n.stewart@warwick.ac.uk}


\date{\today}

\begin{document}

\allowdisplaybreaks

\begin{abstract} 

We classify connected 3-node restricted excitatory-inhibitory networks,
extending 
our previous paper (`Classification of 2-node Excitatory-Inhibitory Networks', 
{\it Mathematical Biosciences} {\bf 373} (2024) 109205). We assume that there are two node-types and two 
arrow-types, excitatory and inhibitory; all excitatory arrows are identical and all inhibitory arrows are identical; and
excitatory (resp. inhibitory) nodes can only output  excitatory (resp. inhibitory) arrows. 
The classification is performed under the following  two network perspectives: 
ODE-equivalence and minimality; and valence  $\le 2$.  
The results of this  and the previous work constitute a first step towards analysing dynamics and bifurcations
of excitatory-inhibitory networks and have potential applications to biological network models.
\end{abstract}

\maketitle

\section{Introduction}

\begin{figure}[h!]
\centerline{%
\includegraphics[width=0.9\textwidth]{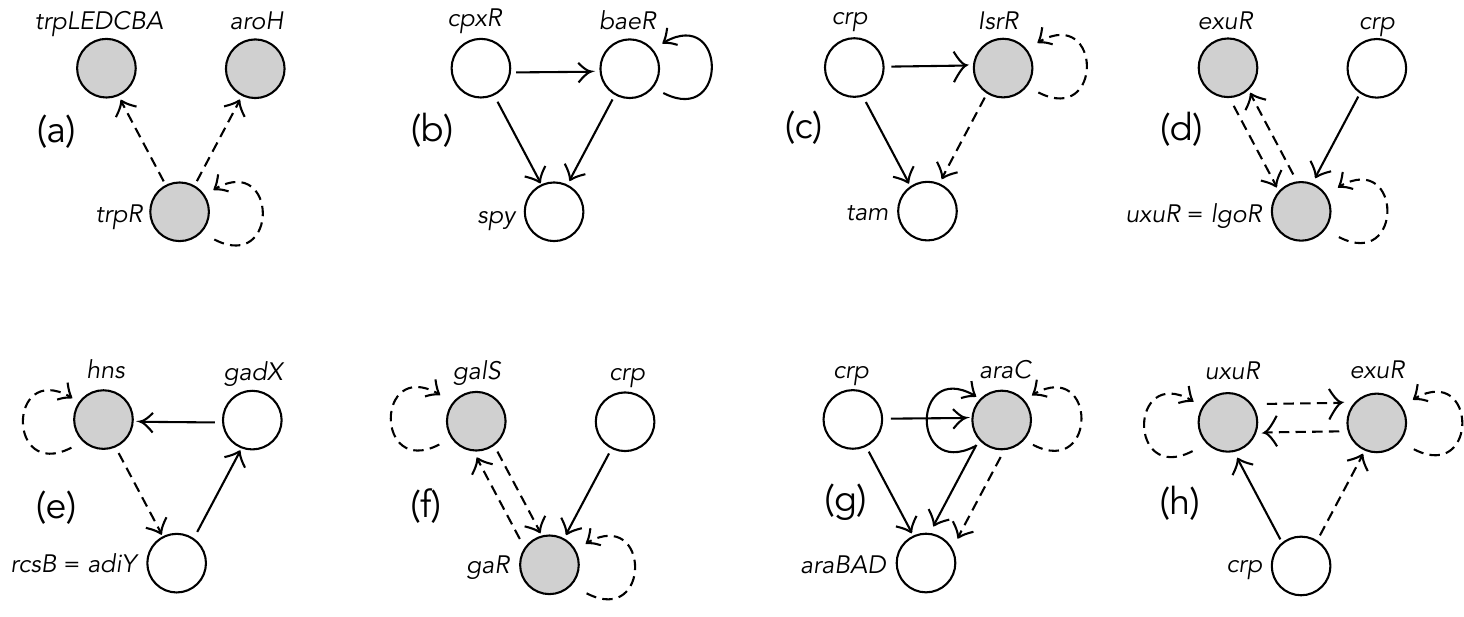}
}
\caption{Eight 3-node motifs realized in {\it E. coli}:
(a) Autoregulation loop involved in biosynthesis of tryptophan, regulated
by {\it trpR} \cite{GG87}, which represses itself, the gene {\it aroH}, and the 
{\it trpLEDCBA} operon, 
which codes for the enzymes of the 
tryptophan biosynthesis pathway. From \cite{MSDS23}.
(b) Example of a SAT-Feed-Forward-Fiber network.
From \cite{LMRASM20} Fig.1 E.
(c) Example of an UNSAT-Feed-Forward-Fiber network.
From \cite{LMRASM20} Fig.2 F.
(d) Example of a 2-FF network showing quotient by synchrony
of genes {\it uxuR} and {\it IgoR} in a 4-node network in {\it E.coli.}
From \cite{MLM20} Fig. 3B.
(e) Example of a 3-FF network showing quotient by synchrony
of genes {\it rcsB} and {\it adiY} in a 4-node network in {\it E.coli.}
From \cite{MLM20} Fig. 3B.
(f)  Example of a network where a node feeds 
forward into one node of a toggle-switch. From \cite{M20}.
(g) In the sugar utilisation transcriptional system [24], the arabinose metabolism [25] 
involves the regulation of the {\it araBAD} operon (composed of
genes {\it araB}, {\it araA}, and {\it araD}) by two transcription factors
{\it araC} and {\it crp} expressed by genes {\it araC} and {\it crp}, respectively.
From \cite{MMpre}.
(h) Example of a network where a node feeds 
forward into both nodes of a toggle-switch. From \cite{MMpre}. 
}
\label{F:EI_E_coli}
\end{figure}

Motifs, small subnetworks that carry out specific functions and occur unusually often, are important building blocks of biological networks. 
See, for example, \cite{B24, MSIKCA02,TN10}. 
Therefore, the classification of small excitatory-inhibitory networks and their dynamical analysis is a fundamental step in the understanding of the dynamics of biological networks and, consequently, in obtaining answers to important biological questions. 
Figure \ref{F:EI_E_coli} illustrates nontrivial 3-node motifs present
 in real biological networks.  
More concretely, it shows eight 3-node motifs from the
gene regulatory network of  {\it Escherichia\ coli}, an organism whose
 genetic regulatory network, compiled by RegulonDB,
  has been characterized in considerable detail
 \cite{GG16}.
For more detail and examples of biological network motifs, see \cite{ADS24}.

The importance of biological network  motifs, and their dynamics and bifurcations, leads to our interest in formalizing the structure of excitatory-inhibitory (EI) networks and to investigate small examples systematically. 
This was the motto for our work in \cite{ADS24}, where we classify connected 2-node excitatory-inhibitory networks under various conditions. 

We work in the coupled cell network formalism of 
\cite{F04, GS06, GS23,GST05, SGP03},
in which nodes (cells) and arrows (connections, directed edges) are partitioned into one or more types.
In biological networks it is common to distinguish between two types of connection:
 {\it excitatory} and {\it inhibitory}. In standard models these have
 different dynamic effects. In the  coupled cell formalism we
 represent this distinction by assuming that nodes and arrows have
two distinct types. For convenience, we call these `excitatory' and `inhibitory', but the 
classification is independent of their dynamics.

In the general theory, the
dynamics of the network can be prescribed by any system
of ordinary differential equations (ODEs) that
respects both its topology and the distinction between different types of node or arrow. Such  systems of ODEs are said to be {\it admissible} for the network. 
The dynamical interpretation of nodes or arrows as being excitatory (tending to activate the nodes to which
they connect) or inhibitory (tending to suppress such activity) is not built
into the definition of admissible ODEs, because connections can differ in other ways.
See \cite[Section 1.3]{ADS24}
 for remarks on how excitation and inhibition can
be defined within the formalism for specific ODE models.

The classification of 2-node excitatory-inhibitory networks in \cite{ADS24} considers different possibilities regarding whether the distinction between the two types of node is  maintained, or they are identified, and regarding whether a node can send only one type of output, excitatory or inhibitory, or can have both excitatory and inhibitory outputs. This leads to four different types of excitatory-inhibitory networks: restricted, partially restricted, unrestricted and completely unrestricted.
For each type we give in \cite{ADS24} two different classifications. Using 
results on ODE-equivalence and minimality, we classify the ODE-classes and present a minimal representative 
for each ODE-class. We also classify all the networks with valence $\le 2$.

In this work, as a continuation of \cite{ADS24}, we extend the classification to 3-node excitatory-inhibitory networks. 
However, here  we assume the type of connection is determined by its tail node, as happens for general neuronal networks. 
In other words, excitatory nodes output excitatory signals and inhibitory nodes output inhibitory signals.  This is what we call restricted excitatory-inhibitory (REI)  networks in Definition \ref{def:EIN} below. 
In Figure \ref{F:EI_E_coli}, networks (a)-(b) have arrows (and nodes) of a single type.
Networks (c)-(f) are REI networks. 
Networks (g)-(h) are not REI networks: some node outputs arrows of both types.

\subsection*{An Example}

 \begin{figure}
 \begin{tabular}{cc}
\begin{tikzpicture}
 [scale=.15,auto=left, node distance=2cm
 ]
 \node[fill=black!30,style={circle,draw}] (n1) at (4,0) {\small{1}};
  \node[fill=white,style={circle,draw}] (n3) at (20,0) {\small{3}};
 \node[fill=black!30,style={circle,draw}] (n2) at (12,-8) {\small{2}};
 \path (n1) [line width=1pt,->,dashed] edge[bend left=20, below] node {} (n2);
 \path (n2) [line width=1pt,->,dashed] edge[bend left=20] node {} (n1);
 \path (n3) [->,thick] edge[] node {} (n2);
\path (n1)  [line width=1pt,->]  edge[loop left=90,dashed] node {} (n1);
\path (n2)  [line width=1pt,->]  edge[loop right=-90,dashed] node {} (n2);
 \end{tikzpicture}  
 \quad & \quad 
 \begin{tikzpicture}
 [scale=.15,auto=left, node distance=2cm
 ]
 \node[fill=black!30,style={circle,draw}] (n1) at (4,0) {\small{1}};
  \node[fill=white,style={circle,draw}] (n3) at (20,0) {\small{3}};
 \node[fill=black!30,style={circle,draw}] (n2) at (12,-8) {\small{2}};
 \path (n1) [line width=1pt,->,dashed] edge[bend left=20, below] node {} (n2);
 \path (n2) [line width=1pt,->,dashed] edge[bend left=20] node {} (n1);
 \path (n3) [->,thick] edge[] node {} (n2);
 \end{tikzpicture} 
 \end{tabular}
\caption{Two $3$-node REI networks which are ODE-equivalent  to the 3-gene GRN  motif in Figure~\ref{F:EI_E_coli} (f) where a node feeds 
forward into one node of a toggle-switch. The network on the right is minimal.} \label{fig:ex_intro}
\end{figure}
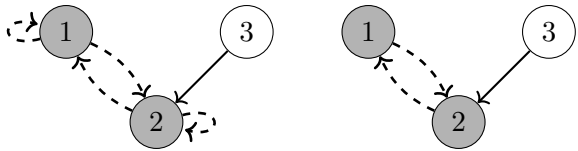

The 3-node network motif (f) of Figure~ \ref{F:EI_E_coli} is an example of an REI network. The black shaded nodes are  of one type (say, inhibitory) and the third node is  of different type (say, excitatory). Both inhibitory nodes send two  inhibitory outputs, which in this network, are directed to the  two inhibitory nodes; 
the excitatory node sends an excitatory signal to one of the inhibitory nodes. In the coupled cell network formalism the main features we retain from this particular network are that it has 3 nodes, two of them are of one type and the third one is of different type. The equal type nodes output arrows of the same type. Different node types output different arrow types. See Figure~\ref{fig:ex_intro} left. A general admissible system of ODEs consistent with this network has the form 
\begin{equation}\label{eq:intro_example}
\begin{array}{l}
\dot{x}_1 = f(x_1; \overline{x_1, x_2}),\\
\dot{x}_2 = g(x_2; \overline{x_2, x_1},x_3),\\
\dot{x}_3 =h(x_3),
\end{array}
\end{equation}
where $f,g,h$ are smooth functions. Each such function captures how the evolution of each node depends on the 
other nodes. The overbar notation over two variables in the functions $f$ and $g$ denotes their invariance under permutation of the two variables, which
occurs because the corresponding input arrows have the same type.
Assuming nodes $1,2$ have internal phase space $\R^k$ and node $3$ has internal phase space $\R^l$, then
$f:\, (\R^k)^3 \to \R^k$, $g:\, (\R^k)^3  \times \R^l \to \R^k$ and $h:\, \R^l \to \R^l$.

Interpreting the network as 
a 3-gene  {\it Escherichia coli} GRN, we may assume that the variable $x_i =( x_i^R, x_i^P) \in \R^2$ is associated with gene $i$, for $i=1,2,3$. Here, 
$x_i^R$ is the concentration of mRNA in gene $i$ and $x_i^P$ is the concentration of protein in gene $i$. We also assume that the time evolution of the cellular concentration of proteins 
and mRNA molecules is determined by an ODE.(We use this term for a single ODE and for a system.)
Moreover, there must be the constraint that a concentration cannot be negative. In this 
modeling approach, two components of the ODE
 are associated with each  gene $i$. The equation for $x_i^R$  determines  the rate of change
of  the concentration of the transcribed mRNA; the equation for  $x_i^P$  describes the rate of change of the concentration of its
corresponding translated protein. As in \cite{MRS24}, a simple example
of an admissible system of the form (\ref{eq:intro_example}), where all 3 genes have 2-dimensional node spaces (that is, $k=l=2$), arises by choosing the following functions $f,g,h$: 
\begin{equation}\label{eq:intro_GRN_example}
\begin{array}{rcl}
f(x_1; \overline{x_1; x_2}) & = &  
\left[ 
\begin{array}{c}
-\delta_1 x_1^R \\
\beta_1 x_1^R - \alpha_1 x_1^P
\end{array}
\right]
+ 
\left[ 
\begin{array}{c}
H_1^{-} (x_1^P) +  H_1^{-} (x_2^P)\\
0
\end{array}
\right]; \\
 \\
g(x_2; \overline{x_2, x_1},x_3) & = & 
\left[ 
\begin{array}{c}
 -\delta_2  x_2^R\\
 \beta_2 x_2^R - \alpha_2 x_2^P
 \end{array}
 \right] +
 \left[ 
\begin{array}{c}
H_2^{-} (x_2^P) +  H_2^{-} (x_1^P) +  H_2^{+} (x_3^P)\\
0
\end{array}
 \right]; \\
  \\
h(x_3) & = &   
\left[ 
\begin{array}{c}
-\delta_3 x_3^R \\
\beta_3 x_3^R - \alpha_3 x_3 ^P
\end{array}
\right]\, .
\end{array} 
\end{equation}

Here, as genes $1,2$ are of the same type, we take $\beta_1 = \beta_2$, $\alpha_1 = \alpha_2$ and $\delta_1 = \delta_2$. 
 Also, $\delta_i, \alpha_i$ represent, respectively, degradation of mRNA and protein for gene $i$,
 and are assumed to be independent of
the concentrations of the other molecules in the cell. 
  The function  $H_i^{-} (x_j^P)$ (resp.  $H_i^{+} (x_j^P)$)  
 in the equation for gene $i$  describes how protein $j$ inhibits (resp. activates) mRNA $i$. In this model equation we assume that the effects of the proteins are additive; 
an alternative typical modeling assumption is that they are multiplicative.  See for example \cite{PHB09}. 
These functions $H_i^-$ and $H_i^+$ are generally nonlinear. Typical choices for $H_i^{-}$ are the Hill functions: 
 $$
 H_i^{-} (z) = \frac{1}{1 + z^{n_i}}
 $$
 where $n_i$ is a positive integer. Assuming $z \geq 0$, since it represents a concentration, $H_i^{-} (z)$ converges to $0$ as $z$ converges to $+\infty$ 
 and $H_i^{-} (0) = 1$. This property encodes inhibition
into the equations. Assuming the inhibitory edges to be of the  same type corresponds to taking 
 $H_1^{-} = H_2^{-} = H_3^{-}$.  A choice for excitation
is the function 
 $$
 H_i^{+} (z) = 1 -  H_i^{-} (z) = \frac{z^{n_i^{*}}}{1 + z^{n^{*}_{i}}}, 
 $$ 
where $n^{*}_i$ is not necessarily equal to $n_i$.
With the functions $f,g,h$ as in (\ref{eq:intro_GRN_example}), and taking into account 
 the structure of network on the left of Figure~\ref{fig:ex_intro} (or  the 3-gene GRN  motif in Figure~\ref{F:EI_E_coli} (f)), equations (\ref{eq:intro_example}) 
 take the form:
 \begin{equation}\label{eq:intro_GRN_net_and_example}
\begin{array}{rcl}
\dot{x}_1^R & = & -\delta_1 x_1^R + H_1^{-} (x_1^P) +  H_1^{-} (x_2^P),\\
\dot{x}_1^P & = & \beta_1 x_1^R - \alpha_1 x_1^P,\\
\\
\dot{x}_2^R & =  & -\delta_1  x_2^R + H_1^{-} (x_2^P) +  H_1^{-} (x_1^P) +  H_2^{+} (x_3^P)\,\\
\dot{x}_2^P & =  & \beta_1 x_2^R - \alpha_1 x_2^P,\\
\\
\dot{x}_3^R & = & -\delta_3 x_3^R,\\
\dot{x}_3^P & = & \beta_3 x_3^R - \alpha_3 x_3 ^P\, .\
\end{array}
\end{equation}

Thus, for $i=1,2$, the rate of change of  the concentration of the transcribed mRNA $i$, given by $x_i^R$, is the difference between 
 the `synthesis term' ($H_1^{-} (x_1^P) +  H_1^{-} (x_2^P)$ for $i=1$ and $H_1^{-} (x_1^P) +  H_1^{-} (x_2^P) +  H_2^{+} (x_3^P)$ for $i=2$), and the `degradation term' $\delta_1 x_i^R$.  In fact, we can think that the evolution of gene $i$ given by $x_i = ( x_i^R, x_i^P)$ is a sum of two parts: one determines the internal dynamics of the gene $i$ and the other determines the coupling effect. For $i=1$,
we can consider the internal dynamics to be determined by
$$
\left[ 
\begin{array}{c}
-\delta_1 x_1^R  \\
\beta_1 x_1^R - \alpha_1 x_1^P
\end{array}
\right], 
$$
and the coupling part by
$$
\left[ 
\begin{array}{c}
  H_1^{-} (x_1^P) + H_1^{-} (x_2^P)\\
0
\end{array}
\right]\, .
$$
Alternatively, we can consider the internal dynamics to be determined by
$$
\left[ 
\begin{array}{c}
-\delta_1 x_1^R +  H_1^{-} (x_1^P) \\
\beta_1 x_1^R - \alpha_1 x_1^P
\end{array}
\right], 
$$
and the coupling part by
$$
\left[ 
\begin{array}{c}
  H_1^{-} (x_2^P)\\
0
\end{array}
\right]\, .
$$

This can be interpreted as considering different  gene internal dynamics of gene $1$. Similarly, we have two analogous options for the internal dynamics of gene $2$. 
Taking the second option for the internal dynamics of genes $1$ and $2$,
we may rewrite (\ref{eq:intro_GRN_example}) as 
\begin{equation}\label{eq:intro_GRN_nsc_example}
\begin{array}{rcl}
f(x_1; \overline{x_1; x_2}) & = &  
\left[ 
\begin{array}{c}
-\delta_1 x_1^R + H_1^{-} (x_1^P)  \\
\beta_1 x_1^R - \alpha_1 x_1^P
\end{array}
\right]
+ 
\left[ 
\begin{array}{c}
 H_1^{-} (x_2^P)\\
0
\end{array}
\right]; \\
 \\
g(x_2; \overline{x_2, x_1},x_3) & = & 
\left[ 
\begin{array}{c}
 -\delta_1  x_2^R + H_1^{-} (x_2^P) \\
 \beta_1 x_2^R - \alpha_1 x_2^P
 \end{array}
 \right] +
 \left[ 
\begin{array}{c}
H_1^{-} (x_1^P) +  H_2^{+} (x_3^P)\\
0
\end{array}
 \right]; \\
  \\
h(x_3) & = &   
\left[ 
\begin{array}{c}
-\delta_3 x_3^R \\
\beta_3 x_3^R - \alpha_3 x_3 ^P
\end{array}
\right]\, .
\end{array} 
\end{equation}

In the coupled cell network formalism, the vector field (\ref{eq:intro_GRN_nsc_example}) determines an admissible coupled cell system for the network on the right of Figure~\ref{fig:ex_intro},  
which has the general form 
\begin{equation}\label{eq:intro_nsc_example}
\begin{array}{l}
\dot{x}_1 = F(x_1; x_2),\\
\dot{x}_2 = G(x_2; x_1,x_3),\\
\dot{x}_3 =h(x_3), 
\end{array}
\end{equation}
where
$$
F(x_1; x_2) = f(x_1; \overline{x_1, x_2}), \quad G(x_2; x_1,x_3) = g(x_2; \overline{x_2, x_1},x_3)\, .
$$
In the coupled cell network formalism, we say that the two networks of Figure~\ref{fig:ex_intro} are ODE-equivalent, precisely because every 
admissible ODE for the network on the right of Figure~\ref{fig:ex_intro} can be seen as 
an admissible ODE for the network on the left of Figure~\ref{fig:ex_intro}, and conversely,
assuming the node phase spaces of the two networks are the same. 
Moreover, the network on the right of Figure~\ref{fig:ex_intro} is the minimal network in terms of number of edges among all 3-node networks 
that are ODE-equivalent to the networks in Figure~\ref{fig:ex_intro}. See Subsection \ref{S:AODE} for formal definitions and main results on network 
admissible ODEs, ODE-equivalence and minimality. In this paper, we use results on network ODE-equivalence and minimality to classify the set of 
3-node REI networks into ODE-classes and present minimal representatives for each ODE-class.

Our classification of 3-node REI networks is made under a variety of extra conditions,
summarized in Table \ref{T:classifications}. 
This classification, together with that for connected 2-node EI networks in \cite{ADS24},
 are a preparatory step towards a
systematic analysis of dynamics and bifurcations in EI networks.

\subsection*{Summary of Paper and Main Results}

We characterize and classify connected 3-node REI networks.
We give a 
classification under the relation of ODE-equivalence, where
two networks are ODE-equivalent if they have the same space
of admissible ODEs.
Sometimes we consider a restriction on the valence of the nodes.
To organize and summarize these results,
Table \ref{T:classifications} lists the main classifications obtained in this paper,
with columns for type of network, bounds on the valence, number of networks in
the classification, plus
references to associated Figures, Tables and Theorems.

\begin{table}[!htb] 
\begin{center}
{\tiny 
\begin{tabular}{|c|c|c|c|}
\hline
network& number of &figure & theorem  \\
 type & networks & & \\
\hline
\hline
  REI & $\infty$ & Figure \ref{fig:generalREI_3node} & Proposition \ref{prop:3_node_connected}  \\
 & & Table \ref{table:ODE_3node_connected} &  \\
   \hline 
  REI   (ODE) & $\infty$ & Figure \ref{fig:generalREI_3node} & Proposition \ref{prop:3_node_ODE}   \\
\hline 
\hline
  REI (ODE) val $\leq 2$ & 92 & Figure \ref{fig:generalREI_3node}               & Proposition \ref{P:4.3}   \\
   no auto 2 arrow-types &      & Table \ref{table:ODE_3node_valence_up_2_a}   & \\
  \hline  
 REI (ODE) val $\leq 2$ & 38 &   Figure \ref{fig:generalREI_3node}         & Proposition \ref{P:4.3}  \\
 no auto 1 arrow-type & & Table \ref{table:ODE_3node_valence_up_2_b}         &   \\
\hline
\hline
 REI (ODE) val $\leq 2$ & 62  &    Figure \ref{fig:generalREI_3node}   & Proposition \ref{P:4.3}  \\
  auto 2 arrow-types  &  & Table \ref{table:ODE_3node_valence_up_2_c} &   \\
\hline
 REI (ODE) val $\leq 2$ & 35 &    Figure \ref{fig:generalREI_3node}       & Proposition \ref{P:4.3} \\
 auto 1 arrow-type & &   Table \ref{table:ODE_3node_valence_up_2_d}                &\\
\hline
\hline
REI val $\leq 2$ & $> 227$ & Figure \ref{fig:3NCNREIV2} & Proposition \ref{prop:REI3}   \\
\hline
\hline
REI val $= 2$, different conditions & --- & Figures  {\rm \ref{fig:hom_3NCNREIV2}}, {\rm \ref{fig:12_3NCNREIV2}}, {\rm \ref{fig:13_3NCNREIV2}}, {\rm \ref{fig:123_3NCNREIV2}}& Propositions \ref{prop:hom_REI3min}, \ref{prop:hom_REI3min2}, \ref{prop:hom_REI3}, \ref{prop:non_hom_REI123}\\
\hline
\end{tabular}
}
\vspace{5mm}
\caption{List of classifications of connected 3-node REI networks and their locations. (ODE): ODE-equivalence classes. val: valence. auto: with autoregulation. no auto: without autoregulation.
In the penultimate line of the table, the exact number of  3-node connected REI networks of valence $\leq 2$ can be obtained by taking all combinations of the multiplicities in Figure~\ref{fig:3NCNREIV2}.}
\label{T:classifications}
\end{center}
\end{table}

Section \ref{S:CEIN} discusses REI networks from
the point of view of the general network formalism of \cite{GS23,GST05,SGP03}.
Subsection \ref{S:FD} gives a formal definition of `restricted excitatory-inhibitory' (REI)  networks.
Subsection \ref{S:AODE} defines the class of admissible ODEs associated with an
REI network. Adjacency matrices
are also discussed. 

Section \ref{S:CC3EIN} characterizes connected 3-node REI networks and classifies
them up to ODE-equivalence.
Corresponding admissible ODEs are not listed, for reasons of space, but
can be deduced algorithmically from the network diagrams.
Subsection \ref{S:C3REIVLE2} classifies the connected $3$-node REI networks with valence $\leq 2$ and also classifies their ODE-classes. 
Subsection \ref{S:C3REIV2} classifies connected $3$-node REI networks with valence $2$ under four different conditions: (i) every node receives one arrow of each type; 
(ii) only the two excitatory nodes receive one arrow of each type; (iii) only the inhibitory node and one excitatory node receive one arrow of each type; 
(iv) given any two nodes there is no arrow-type preserving bijection between their input sets.

\section{
Restricted Excitatory-Inhibitory Networks}
\label{S:CEIN}

In this section we define the class of {\it restricted} EI-networks (REI).  We assume the networks have two distinct node-types $N^E,N^I$ and two different arrow-types $A^E,A^I$, which we may think of as excitatory/inhibitory  nodes and excitatory/inhibitory arrows. Moreover, we make the standard simplified modeling assumption 
that all excitatory arrows are identical and all inhibitory arrows are identical.
Without this last assumption, the lists of networks becomes much larger, already for the class of 3-node networks.

In some areas of biology, notably
neuroscience, a given node cannot output both an excitatory arrow and an
inhibitory one. We make that assumption here. 
Also, as in \cite{ADS24}, we work in the modified network formalism presented in \cite{GS23}, 
which allows arrows of the same type to have heads of different types. 
This differs from
 the formalism of \cite{GST05,SGP03}, in which arrows of the same type have heads (and tails) of the same type. We remove that condition so
that an excitatory (resp. inhibitory) node can send excitatory (resp. inhibitory) arrows to excitatory and/or inhibitory nodes.  See \cite[Section 9.3]{GS23} for technical details where it is pointed out that  the main network theorems and their proofs remain valid in the more general
formalism. See also \cite[Remarks 2.1]{ADS24} for a discussion of this approach.

\subsection{Formal Definitions} 
\label{S:FD}

We define restricted excitatory-inhibitory (REI) networks, state our conventions for representing
them in diagrams, and give examples.

\begin{Def} \label{def:EIN} 

A network ${\mathcal G}$ is a {\it restricted excitatory-inhibitory network} ({\it REI network}) if it satisfies the following four conditions:

(a) There are two distinct node-types, $N^E$ and $N^I$.

(b) There are two distinct arrow-types, $A^E$ and $A^I$.

(c) If $e \in A^E$ then $\TT (e) \in N^E$.

(d) If $e \in A^I$ then $\TT (e) \in  N^I$, 

\noindent
where $\TT(e)$ indicates the tail node of arrow $e$. 
\hfill $\Diamond$
\end{Def}

\vspace{5pt}
\paragraph{\bf Conventions}

The following conventions are used throughout the paper without further mention,
except as an occasional reminder for clarity. 

(a)	We represent type $N^E$ nodes by white circles and type $N^I$ nodes by grey circles. Type $A^E$ arrows are solid and type
$A^I$ arrows are dashed. (Various other conventions for excitatory/inhibitory arrows
are found in the literature; this one is chosen for convenience.) 
 
 (b) All classifications are stated up to {\it renumbering} of nodes and {\it duality}; that is,
 interchange of `excitatory' and `inhibitory' on nodes and arrows: 
 $N^E \leftrightarrow N^I$ and $A^E \leftrightarrow A^I$.
 \hfill $\Diamond$
 
 \begin{exam}
 The networks (c)-(f)  in Figure \ref{F:EI_E_coli} are REI networks. However, 
networks (g)-(h) are not REI networks as  some node (the {\it araC} gene in network (g)  and  the {\it crp} gene in network (h)) outputs arrows of both types. 
\hfill $\Diamond$
 \end{exam}
 
\begin{Def} \label{Def:input_equiv}
(a) 
In an REI network, every node $i$ can receive excitatory and inhibitory arrows: here, the sets of excitatory and inhibitory arrows directed to $i$  are denoted by  $I^E(i)$ and  $I^I(i)$, and called  the {\it excitatory} and {\it inhibitory input sets} of $i$, respectively. The union $I(i) = I^E(i) \cup I^I(i)$ is the {\it input set} of $i$ and the cardinality  $\#I(i)$ of $I(i)$ is the {\it valence (degree, in-degree)} of $i$. 

(b) 
Two nodes $i$ and $j$ with the same node-type and valence are said to be  {\it input equivalent} when $\#I^E(i) = \#I^E(j)$ and $\#I^I(i) = \#I^I(j)$.
We write $i \sim_I j$.
Trivially,
 the relation $\sim_I$ is an equivalence relation, which partitions the set of nodes into
 disjoint {\it input classes}. \\
(c) 
A network where the nodes are not all input equivalent is {\it inhomogeneous}. Otherwise, it is {\it homogeneous}. 
 \hfill $\Diamond$
\end{Def}

\begin{rems}
\label{R:2.3}
(a)
Every REI network is inhomogeneous as by definition it has two distinct node-types,  $N^E $ and $N^I$. \\
(b) The definition of (robust) synchrony in \cite{GS23,GST05,SGP03} implies that
synchronous nodes must be input equivalent. Thus for EI networks, 
nodes of type $N^E$ cannot synchronize with nodes of type $N^I$. See also Subsection~\ref{sub:synchrony}. 
\hfill $\Diamond$
\end{rems}

In this paper we consider {\it connected} networks in the sense there is an undirected path between every pair of nodes. We distinguish connected 
networks according to the existence of a closed directed arrow-path containing every node, or not. In the first case, the network is {\it transitive}. Otherwise, it is {\it feedforward}. 

\begin{exam}
Consider the two networks (e)-(f) in Figure \ref{F:EI_E_coli}. Network (e)  is transitive and network (f) is feedforward. 
\hfill $\Diamond$
\end{exam}

\subsection{Admissible ODEs}
\label{S:AODE}
We adopt the general form of admissible ODEs for a network as defined in  \cite{GS23,GST05,SGP03} 
with the assumption in this paper that all nodes have the same state space, say $P = \R^m$ for some $m > 0$.
Given an EI network with a finite set of nodes, node $i$ is represented in the ODE system by the variable 
$x_i$ which is governed by a system of  ordinary differential equations. The word `admissible' is used in the sense that 
the ODE system encodes information about the node and arrow types. Specifically, 
when two input equivalent nodes have the same numbers,  say $n_e$, of excitatory arrows and  $n_i$ of inhibitory 
arrows, targeting the two nodes, we specify their dynamics by the same smooth function, 
say $f : P^{k+1} \to P$, evaluated at the node and at the corresponding tail nodes of the 
arrows targeting the node. We follow \cite[Definition~2.8]{ADS24}: 

\begin{Def}
A system of ODEs is {\it admissible} for an EI network if it has the form
\[
\dot{x}^s_i = f_i(x^s_i; \overline{x^+_{i_1}, \ldots, x^+_{i_{n_e}}}; \overline{x^-_{i_{n_e+1}}, \ldots, x^-_{i_{n_e+n_i}}}) 
\]
where $x^s_i \in \{x^+_i,x^-_i \}$ and 
the overlines indicate that the function $f_i$ is symmetric in the 
overlined variables.
The node variables are indexed by $i$. 
The {\it multiset} of all tail nodes of
input arrows is the union of two subsets:  
the multiset $\{i_1, \ldots, i_{n_e}\}$ of all tail nodes of the
excitatory input set of node $i$, and the multiset $\{i_{n_e+1}, \ldots, i_{n_e+n_i}\}$
of all tail nodes of  the inhibitory input set of node $i$. 
The functional notation converts
these multisets into tuples of the corresponding variables. 
We use the superscripts $+$ and $-$, as a notation convention, to make the distinction between the input variables corresponding to tail nodes in the excitatory and in the inhibitory input sets, respectively. Analogously, when there are two distinct node-types $N^E$ and $N^I$, we use the superscripts $+$ and $-$ to make the distinction between the state variable of excitatory and inhibitory nodes. 

Moreover, if nodes $i,j$ of the same node-type 
are  in the same input class, 
that is,  there is an arrow-type preserving bijection between the corresponding input sets, 
then $f_i=f_j$. The evolution of nodes in different input classes is governed by 
 different functions $f_i$, one for each input class. 
\hfill $\Diamond$
\end{Def}

\begin{rem}
Observe that {\it multiple arrows} are permitted as there can be distinct excitatory (resp. inhibitory) arrows with the same tail node directed to the same node. 
Moreover, {\it self-loops} are also permitted as a node can input an arrow to itself.  In biology, the term {\it autoregulation} is used when a node influences
its own state. 
\hfill $\Diamond$
\end{rem}

\begin{exam} \label{Ex:smolen}
The UNSAT-Feed-Forward-Fiber network in Figure~\ref{F:EI_E_coli} (c), which is one of the 3-node motifs from the
gene regulatory network of  {\it Escherichia\ coli}, is an REI (inhomogeneous) network. 
Nodes `crp' and `tam' are type $N^E$ and node `IsrR' is type $N^I$. We number them as nodes 1, 2 and 3, respectively.
There are two  type $A^E$ arrows; one from $1$ to $2$ and the other from $1$ to $3$. 
There are two type $A^I$ arrows; one from $3$ to itself and the other from $3$ to $2$. 

Node $1$ has empty input set. 
Nodes $2$ and $3$ have excitatory and inhibitory input sets with cardinality $1$.
Node $3$ is autoregulatory. 
Thus, although nodes $1$ and $2$ are of same type, they are not input equivalent, since they have different valences.
On the other hand,  although nodes $2$ and $3$ have same excitatory and inhibitory input valences, they are not input equivalent, since they are of different types.

Admissible ODEs are:
\begin{equation}
\label{E:smolen}
\begin{array}{l}
\dot{x}_1^+ 
= f(x^+_1)\\
\dot{x}_2^+
= g(x^+_2; x^+_1; x^-_3) \\
\dot{x}_3^- 
= h(x^-_3; x^+_1; x^-_3)
\end{array}\, .
\end{equation}
Here, $x^+_1,x^+_2, x^-_3 \in P$, where $P$ is the node state space, and $f:\, P  \to P$ and $g, h:\, P^3  \to P$ are smooth functions. 
\hfill $\Diamond$
\end{exam}

An $n$-node network can be represented by its {\it adjacency matrix}, which is the $n \times n$
matrix $A = (a_{ij})$ such that 
$a_{ij}$ is the number of arrows from node $j$ to node $i$. (In the graph-theoretic literature the opposite convention is often used, which gives the
transpose of the adjacency matrix defined here.) 
For an REI network, conditions (c)-(d) of Definition~\ref{def:EIN} allow us to 
deduce the arrow-types from its adjacency matrix, provided we know the node-types of nodes $i$ and $j$. In fact, we consider two {\it node-type} 
$n \times n$ matrices, which are both diagonal: given one node-type matrix,   the diagonal entry $ii$ is 1 if node $i$ is of that 
type and zero otherwise. 
When we need to distinguish the different arrow-types, as is the case in this paper when classifying networks using 
ODE-equivalence, see Subsection~\ref{sub:ODE_equiv}, we 
consider {\it arrow-type} adjacency matrices,  one for each arrow type. 
For example, for REI-networks, we will consider two  arrow-type adjacency matrices, one for  excitatory arrows and the other for inhibitory arrows.

\begin{exam} 
\label{ex:smolen_adj}
The adjacency matrix of the UNSAT-Feed-Forward-Fiber network in Figure~\ref{F:EI_E_coli} (c) is
$$
\left[
\begin{array}{ccc}
0 & 0& 0 \\
1 & 0 & 1 \\
1 & 0 & 1 
\end{array}
\right]\, .
$$

We may also distinguish node- and arrow-types and equip each with its own adjacency matrix. Here there are four:
$$
\begin{array}{ll} 
\mbox{Node-type $N^E$: } 
\left[
\begin{array}{ccc}
1 & 0 & 0\\
0 & 1 & 0 \\
0 & 0& 0
\end{array}
\right]; \quad 
\mbox{Node-type $N^I$: } 
\left[
\begin{array}{ccc}
0 & 0 & 0\\
0 & 0 & 0\\
0 & 0 & 1 
\end{array}
\right]; \\ 
\ \\
\mbox{Arrow-type $A^E$: } 
\left[
\begin{array}{ccc}
0 & 0& 0 \\
1 & 0 & 0 \\
1 & 0 & 0 
\end{array}
\right]; \quad 
\mbox{Arrow-type $A^I$: } 
\left[
\begin{array}{ccc}
0 & 0& 0 \\
0 & 0 & 1 \\
0 & 0 & 1 
\end{array}
\right]\, .
\end{array}
$$
\hfill $\Diamond$
\end{exam}

\subsection{ODE-equivalent Networks} \label{sub:ODE_equiv}
As mentioned and exemplified in the Introduction, different networks 
with the same number of nodes  are 
said to be {\it ODE-equivalent} if they have the same set of admissible ODEs,
 for any choice of node state spaces, when their nodes are identified
 by a suitable bijection that preserves node state spaces. See~\cite{DS05,GS23,GST05}. 
 
 \begin{rems} \label{rem:ODEequiv}
 (a)  
A necessary and sufficient condition for two 
networks to be ODE-equivalent, using the  associated node and arrow adjacency matrices,
is proved in \cite[Theorem 7.1, Corollary 7.9]{DS05}.
Specifically, two networks with the same number of nodes are ODE-equivalent if and only if, for a suitable identification of nodes, they have
 the same vector spaces of {\it linear} admissible maps when node state spaces are $\R$.
 Equivalently, the adjacency matrices of all node- and arrow-types span the same space. \\
  (b) For REI networks, as mentioned above, the node-types determine the arrow-types,
  which implies that the adjacency matrices naturally decompose into four blocks.
The linear condition in (a) preserves this decomposition, so two REI networks are
ODE-equivalent if and only if these components are separately ODE-equivalent. \\
  (c) In fact,  using the results in \cite{AD07,AD07a} on network minimality, it follows that given an ODE-class of 
REI networks, we can distinguish a subclass containing the 
REI networks in the ODE-class that have a minimal number of arrows. This is a {\it minimal subclass} which  
in general need not be a singleton. 
\hfill $\Diamond$
\end{rems}

\begin{figure}
\begin{tikzpicture}
 [scale=.15,auto=left, node distance=2cm
 ]
 \node[fill=white,style={circle,draw}] (n1) at (4,0) {\small{1}};
  \node[fill=black!30,style={circle,draw}] (n3) at (20,0) {\small{3}};
 \node[fill=white,style={circle,draw}] (n2) at (12,-8) {\small{2}};
 \path (n1) [->,thick] edge[] node {} (n2);
 \path (n1) [->,thick] edge[] node {} (n3);
 \path (n3) [->,dashed] edge[] node {} (n2);
 \end{tikzpicture}  
\caption{The minimal $3$-node network ODE-equivalent to the UNSAT-Feed-Forward-Fiber network in Figure~\ref{F:EI_E_coli} (c). } \label{fig:ODE_equiv_Smolen}
\end{figure}
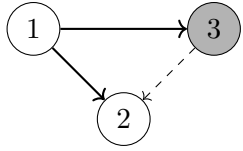

 \begin{exams}
 (a) The  REI network in Figure~\ref{fig:ODE_equiv_Smolen} is ODE-equivalent to the  REI 
UNSAT-Feed-Forward-Fiber network in Figure~\ref{F:EI_E_coli} (c) and it is minimal. 
 Moreover, the admissible ODE \eqref{E:smolen} determines an arbitrary dynamical
system in $(x^+_1,x^+_2, x_3^-)$. \\
(b) The REI network on the right of Figure~\ref{fig:ex_intro} is ODE-equivalent to the REI  3-gene GRN  motif in Figure~\ref{F:EI_E_coli} (f) 
and it is minimal. 
 \hfill $\Diamond$
 \end{exams}
 
\subsection{Robust Network Synchrony Subspaces}\label{sub:synchrony}

 \begin{figure}
 \begin{tabular}{ll}
  \begin{tikzpicture}
 [scale=.15,auto=left, node distance=2cm
 ]
 \node[fill=white,style={circle,draw}] (n1) at (4,0) {\small{1}};
  \node[fill=white,style={circle,draw}] (n3) at (20,0) {\small{2}};
 \node[fill=black!30,style={circle,draw}] (n2) at (12,-8) {\small{3}};
 \path (n1) [line width=1pt,->] edge[bend left=20, below] node {} (n2);
 \path (n2) [line width=1pt,->,dashed] edge[bend left=20] node {} (n1);
  \path (n2) [line width=1pt,->,dashed] edge[bend left=20] node {} (n3);
 \end{tikzpicture} \qquad & \qquad 
 \begin{tikzpicture}
 [scale=.15,auto=left, node distance=2cm
 ]
 \node[fill=white,style={circle,draw}] (n1) at (4,0) {\small{1}};
 \node[fill=black!30,style={circle,draw}] (n2) at (12,-8) {\small{3}};
 \path (n1) [line width=1pt,->] edge[bend left=20, below] node {} (n2);
 \path (n2) [line width=1pt,->,dashed] edge[bend left=20] node {} (n1);
 \end{tikzpicture} 
 \end{tabular}
\caption{
(Left) A $3$-node minimal REI network where nodes $1,2$ can synchronize robustly. (Right) A $2$-node REI network which is the quotient of the 3-node network on the left by taking the equivalence relation on the 3-node network set with classes $\{1,2\}$ and $\{3\}$.} \label{fig:ex_sync}
\end{figure}
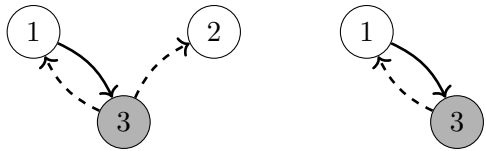

Consider the 3-node REI network on the left of Figure~\ref{fig:ex_sync}. The two excitatory nodes $1,2$ are input equivalent as both receive only one inhibitory arrow. A general admissible ODE-system associated with this network has the form
\begin{equation}\label{eq:sync_example}
\begin{array}{l}
\dot{x}_1^+ = f_1(x_1^+; x_3^-),\\
\dot{x}_2^+ = f_1(x_2^+; x_3^-),\\
\dot{x}_3 ^-=f_3(x_3^-; x_1^+),
\end{array}
\end{equation}
where $f_1,f_3:\, \R^l \times \R^l \to \R^l$ are smooth functions. We see that any solution $( x_1^+(t), x_2^+(t), x_3^-(t))$ of (\ref{eq:sync_example}) with initial 
condition satisfying say $x_1^+(0) = x_2^+(0)$ has nodes $1,2$ synchronized for all time, that is, 
$$
x_1^+(0) = x_2^+(0) \Rightarrow x_1^+(t) =x_2^+(t),\quad \forall t\, .
$$ 
This property does not depend on the choices of the functions 
$f_1,f_3$ neither the internal node phase spaces $\R^l$. It is determined {\it only} by the structure of the network on the left of  
Figure~\ref{fig:ex_sync}; concretely, the two nodes $1,2$ are of the same node type and 
each receives one inhibitory arrow from the  inhibitory node $3$, which in this example is the unique inhibitory 
node. Equivalently, we see that the vector field $F(x_1; x_2; x_3) = ( f_1(x_1; x_3), f_1(x_2; x_3); f_3(x_3; x_1))$ 
leaves invariant the space $\Delta = \{ (x_1, x_1, x_3)\}$, that is, 
$$
 F(\Delta) \subseteq \Delta\, .
 $$
 In this case, we say that $\Delta$ is a {\it robust network synchrony space}. 
Restricting (\ref{eq:sync_example}) to $\Delta$, we obtain the system 
\begin{equation}\label{eq:sync_quo_example}
\begin{array}{l}
\dot{x}_1^+ = f_1(x_1^+; x_3^-),\\
\dot{x}_3 ^-=f_3(x_3^-; x_1^+),
\end{array}
\end{equation}
which is admissible for the 2-node network on the right of Figure~\ref{fig:ex_sync} which is also an REI network. 
In the terminology of \cite{SGP03}, the 2-node network on the right of Figure~\ref{fig:ex_sync} is the {\it quotient} of the network on the left of Figure~\ref{fig:ex_sync} by the equivalence relation on the node set of the 3-node network with equivalence 
classes $\{1,2\}$ and $\{3\}$. This relation is said to be {\it balanced}, which is equivalent to   
the invariance of  $\Delta$ under the node and arrow
adjacency matrices. 

These ideas generalize to $n$-node networks and it is proved in \cite{SGP03} 
that the admissible vector fields for a network leave invariant a linear subspace defined in terms of equalities of certain node coordinates if and only if the equivalence relation on the network node set with classes given by the clusters of nodes whose coordinates are identified is balanced. See  
\cite[Definition~6.4]{SGP03} for the definition of network balanced relation, 
\cite[Proposition 10.20 ]{GS23} or \cite[Section 5]{GST05}, and  \cite[Theorem 6.5]{SGP03} for 
the definition of  quotient network by a balanced equivalence relation. 

For REI networks, it is trivial to show that the restriction of any admissible ODE for an REI network to a robust 
synchrony subspace is admissible for a smaller network, which is also 
an REI network. That is, the quotient of an REI network by a  balanced equivalence relation on the network node set 
is also an REI network.

\section{Classification of Connected 3-node REI Networks}
\label{S:CC3EIN}

We now classify REI networks with three nodes, which we assume are connected.
Moreover, up to duality and numbering of the nodes, we can assume that the
networks have nodes $1$ and $2$  of type $N^E$ and node $3$ of type $N^I$.

\subsection{Connected 3-node REI Networks}
\label{S:C3REIN}
In this section we characterize the connected $3$-node REI networks, without imposing any restrictions, and classify them up to ODE-equivalence.

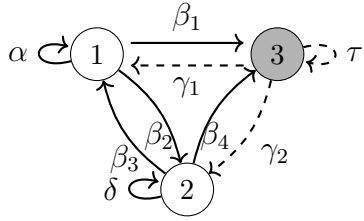
\begin{figure}
\begin{tikzpicture}
 [scale=.15,auto=left, node distance=2cm
 ]
 \node[fill=white,style={circle,draw}] (n1) at (4,0) {\small{1}};
  \node[fill=black!30,style={circle,draw}] (n2) at (20,0) {\small{3}};
  \node[fill=white,style={circle,draw}] (n3) at (12,-12) {\small{2}};
  \path
   (7,1)  [->] edge[thick] node {$\beta_1$}  (17,1)
 (n1)  [line width=1pt,->]  edge[loop left=90,thick] node {$\alpha$} (n1);
 \path 
 (17,-1)  [dashed,->] edge[thick] node {$\gamma_1$}  (7,-1);
 \path 
 (n2)  [dashed,->]  edge[loop right=90,thick] node {$\tau$} (n2);
 \path 
 (n3)  [line width=1pt,->]  edge[loop left=90,thick] node {$\delta$} (n3);
  \path (n2) [->,dashed] edge[bend left=20, thick] node {$\gamma_2$} (n3);
 \path (n3) [->,thick] edge[bend left=20, below] node {$\beta_4$} (n2);
  \path (n1) [->,thick] edge[bend left=20, below] node {$\beta_2$} (n3);
 \path (n3) [->,thick] edge[bend left=20, below] node {$\beta_3$} (n1);
 \end{tikzpicture} 
 \caption{$3$-node REI network: nodes $1$ and $2$ are excitatory and node $3$ is inhibitory. The nonnegative integer arrow multiplicities are $\alpha$, $\delta$,  $\tau$, $\beta_i$, $\gamma_j$, $i=1,2,3,4$, $j=1,2$.} \label{fig:generalREI_3node}
\end{figure}

Up to duality any $3$-node REI network is as shown in Figure~\ref{fig:generalREI_3node}, for a suitable choice of nonnegative integer arrow multiplicities $\alpha$, $\delta$,  $\tau$, $\beta_i$, $\gamma_j$, where $i=1,2,3,4$, $j=1,2$.

The adjacency matrices are  
$$
\begin{array}{ll} 
\mbox{Node-type $N^E$: } 
A_1 = 
\left[
\begin{array}{ccc}
1 & 0 & 0 \\
0 & 1  & 0 \\
0 & 0 & 0
\end{array}
\right]; \quad 
\mbox{Node-type $N^I$: } 
A_2 = 
\left[
\begin{array}{ccc}
0& 0 & 0\\
0 & 0 & 0 \\
0& 0 & 1
\end{array}
\right]; \\ 
\ \\
\mbox{Arrow-type $A^E$: } 
A_3 = 
\left[
\begin{array}{ccc}
\alpha & \beta_3  & 0\\
\beta_2 & \delta & 0 \\
\beta_1 & \beta_4 & 0 
\end{array}
\right]; \quad 
\mbox{Arrow-type $A^I$: } 
A_4= 
\left[
\begin{array}{ccc}
0& 0 & \gamma_1 \\
0& 0 & \gamma_2 \\
0 & 0 & \tau 
\end{array}
\right].
\end{array}
$$

\begin{prop} \label{prop:3_node_connected}
Any $3$-node REI network is as shown in {\rm Figure}~{\rm \ref{fig:generalREI_3node}}, for a suitable choice of nonnegative integer arrow multiplicities $\alpha$, $\delta$, $\tau$, $\beta_i$,  $\gamma_j$, where $i=1,2,3,4$, $j=1,2$.
A $3$-node REI network is connected if and only if its nonzero arrow multiplicities, excluding 
autoregulation arrows, are listed in {\rm Table~\ref{table:ODE_3node_connected}}.
\end{prop}

\begin{proof}
A $3$-node REI network is connected if and only if the union of the input and output sets of each node, excluding self-coupling arrows, is nonempty. That  is, if and only if
at least one multiplicity is nonzero in each of the sets
$$
 \{\beta_1, \beta_2, \beta_3, \gamma_1 \}, \quad  \{\beta_2, \beta_3, \beta_4, \gamma_2 \},  \quad \mbox{and} \quad  \{\beta_1, \beta_4, \gamma_1, \gamma_2 \}.
$$
The possible combinations are listed in Table~\ref{table:ODE_3node_connected}.
\end{proof}

\begin{table}[!ht]
\begin{tabular}{|l|l|l|l|l|}
\hline 
$\beta_1, \beta_2$ & $\beta_1, \beta_2, \beta_3$ & $\beta_1, \beta_2, \beta_3, \beta_4$ & $\beta_1, \beta_2, \beta_3, \gamma_1$ &  $\beta_1, \beta_2, \beta_3,  \gamma_1,\gamma_2$  \\
 \hline 
$\beta_1, \beta_2, \beta_3,  \gamma_2$  & $\beta_1, \beta_2, \beta_3, \beta_4, \gamma_1$ & $\beta_1, \beta_2, \beta_3, \beta_4, \gamma_1, \gamma_2$ & $\beta_1, \beta_2, \beta_3, \beta_4,  \gamma_2$ & $\beta_1, \beta_2,  \beta_4$ \\
 \hline 
  $\beta_1, \beta_2,  \beta_4 , \gamma_1$ & $\beta_1, \beta_2,  \beta_4 ,  \gamma_1,\gamma_2$ & $\beta_1, \beta_2,  \beta_4 ,\gamma_2$ & $\beta_1, \beta_2, \gamma_1$ & $\beta_1, \beta_2, \gamma_1, \gamma_2$ \\
 \hline 
$\beta_1, \beta_2, \gamma_2$ & $\beta_1, \beta_3$ & $\beta_1, \beta_3, \beta_4$ & $\beta_1, \beta_3,  \beta_4,  \gamma_1$ & $\beta_1, \beta_3, \beta_4,  \gamma_1, \gamma_2$  \\
 \hline 
 $\beta_1, \beta_3, \beta_4, \gamma_2$ & $\beta_1, \beta_3, \gamma_1$ & $\beta_1, \beta_3, \gamma_1, \gamma_2$ & $\beta_1, \beta_3, \gamma_2$  & $\beta_1, \beta_4$ \\
 \hline 
$\beta_1, \beta_4, \gamma_1$ & $\beta_1, \beta_4, \gamma_1,  \gamma_2$ & $\beta_1, \beta_4, \gamma_2$ & $\beta_1, \gamma_1,  \gamma_2$ & $\beta_1,  \gamma_2$ \\
 \hline 
$\beta_2, \beta_3, \beta_4$ & $\beta_2, \beta_3, \beta_4, \gamma_1$ & $\beta_2, \beta_3, \beta_4, \gamma_1,  \gamma_2$ & $\beta_2, \beta_3, \beta_4, \gamma_2$ & $\beta_2,\beta_4$ \\
 \hline 
$\beta_2,\beta_4, \gamma_1$  & $\beta_2,\beta_4, \gamma_1,  \gamma_2$ & $\beta_2,\beta_4, \gamma_2$ & $\beta_2, \gamma_1$ & $\beta_2, \gamma_1, \gamma_2$ \\
 \hline 
$\beta_2, \gamma_2$ & $\beta_3, \beta_4$ & $\beta_3, \beta_4, \gamma_1$ & $\beta_3, \beta_4, \gamma_1 , \gamma_2$  & $\beta_3, \beta_4,  \gamma_2$ \\
 \hline 
 $\beta_3, \gamma_1$ &  $\beta_3, \gamma_1,  \gamma_2$ &  $\beta_3,  \gamma_2$& $\beta_4, \gamma_1$ &  $\beta_4, \gamma_1,  \gamma_2$\\
 \hline 
$\gamma_1,  \gamma_2$ & & & & \\
 \hline 
\end{tabular}
\vspace{.2cm}
\caption{Possible nonzero multiplicities of the arrows of a connected REI network as shown in Figure~\ref{fig:generalREI_3node}.}
\label{table:ODE_3node_connected}
\end{table}

\begin{prop} \label{prop:3_node_ODE}
The $3$-node REI networks are those in {\rm Figure}~{\rm \ref{fig:generalREI_3node}}, for a suitable choice of nonnegative integer arrow multiplicities $\alpha$, $\delta$, $\tau$, $\beta_i$,  $\gamma_j$, where $i=1,2,3,4$, $j=1,2$.

Up to ODE-equivalence and minimality, we can assume that $\tau$ is zero and at least one of $\alpha$ or $\delta$ is zero.
\begin{itemize}
\item[]  Moreover, either
\begin{itemize}
\item[] $\gamma_1$ and $\gamma_2$ are coprime, if both nonzero, or
\item[] $\gamma_1 = 1$ and $\gamma_2=0$, or $\gamma_1 = 0$ and $\gamma_2=1$, or
\item[] $\gamma_1 = \gamma_2=0$.
\end{itemize}
\item[] If $\alpha = \delta=0$, then
\begin{itemize}
\item[] the nonzero $\beta_i$, $i=1,\ldots,4$, are coprime, or
\item[] $\beta_i=1$, if $\beta_j =0$, $j \ne i$, $i,j=1,\ldots,4$, or
\item[] $\beta_i=0$, $i=1,\ldots,4$.
\end{itemize}
\item[] If $\alpha \ne 0$ and  $\delta=0$, then
\begin{itemize}
\item[] $\alpha$ and the nonzero $\beta_i$, $i=1,\ldots,4$, are coprime, or
\item[] $\beta_i=0$, $i=1,\ldots,4$.
\end{itemize}
\end{itemize}
\end{prop}
\begin{proof}
Up to ODE-equivalence, we can assume that $\tau=0$, so
$$
\langle  A_1, A_2, A_3, A_4\rangle = \left\langle  A_1, A_2, A_3, 
\left[
\begin{array}{ccc}
0& 0 & \gamma_1 \\
0& 0 & \gamma_2 \\
0 & 0 & 0
\end{array}
\right] \right\rangle\, .
$$

Thus, if $\gamma_2 =0$, we can set $\gamma_1=1$, and {\it vice versa}. If both $\gamma_1$ and $\gamma_2$ are nonzero, we can assume they are coprime.

Moreover, up to ODE-equivalence, we can assume that, at least one, of $\alpha$ or $\delta$ is zero. If $\alpha \ne 0$ and $\delta =0$ we can assume that $\alpha$ and the nonzero $\beta_i$, for $i=1,2,3,4$, are coprime.

If both $\alpha$ and $\delta$ are zero then
$$
\langle  A_1, A_2, A_3, A_4 \rangle = \left\langle  A_1, A_2, 
\left[
\begin{array}{ccc}
0 & \beta_3  & 0\\
\beta_2 & 0 & 0 \\
\beta_1 & \beta_4 & 0 
\end{array}
\right], 
\left[
\begin{array}{ccc}
0& 0 & \gamma_1 \\
0& 0 & \gamma_2 \\
0 & 0 & 0
\end{array}
\right] 
\right\rangle \, ,
$$
then we can assume that the nonzero $\beta_i$, for $i=1,\ldots,4$, are coprime.
\end{proof}

\subsection{Connected 3-node REI Networks with Valence $\le 2$}
\label{S:C3REIVLE2}

In this section we classify the connected $3$-node REI networks with valence $\leq 2$. 
We start by classifying them up to ODE-equivalence.

\begin{table}[htb]
{\tiny 
\begin{tabular}{|l|l|}
\hline 
{nonzero multiplicities}
 & 
{$\#$ ODE-classes}
\\
 \hline 
$1 \le \beta_1, \beta_2 \le 2, \quad \beta_3 = \gamma_1 = 1$ & 4  \\
 \hline 
$1 \le \beta_1 \le 2, \quad \beta_2=\beta_3 = \gamma_1  = \gamma_2= 1$ & 2  \\
\hline  
$1 \le \beta_1, \beta_3 \le 2, \quad \beta_2 = \gamma_2 = 1$ & 4  \\
\hline
$1 \le \beta_2 \le 2, \quad \beta_1=\beta_3 = \beta_4  = \gamma_1= 1$ & 2  \\
\hline
$1 \le \beta_3 \le 2, \quad \beta_1=\beta_2 = \beta_4  = \gamma_2= 1$ & 2  \\
\hline
$\beta_1=\beta_2 = \beta_3 = \beta_4 = \gamma_1  = \gamma_2= 1$ & 1  \\
\hline
$1 \le \beta_2 \le 2, \quad \beta_1 = \beta_4 = \gamma_1 = 1$ & 2  \\
\hline
$\beta_1=\beta_2  = \beta_4  = \gamma_2= 1$ & 1  \\
\hline
$1 \le \gamma_1 \le 2, \quad \beta_1 = \beta_2 = \beta_4 = \gamma_2 = 1$ & 2  \\
\hline
$1 \le \beta_1, \beta_2,  \le 2, \mbox{ excluding } \beta_1 = \beta_2  = 2, \quad \gamma_1 = 1 $ & 3  \\
\hline
$1 \le \beta_1, \gamma_1 \le 2, \quad \beta_2 = \gamma_2 = 1$ & 4  \\
\hline
$1 \le \beta_1 \le 2, \quad \beta_2 = \gamma_2 = 1$ & 2  \\
\hline
$\beta_1=\beta_3  = \beta_4  = \gamma_1= 1$ & 1  \\
\hline
$1 \le \gamma_2 \le 2, \quad \beta_1=\beta_3  = \beta_4  = \gamma_1= 1$ & 2  \\
\hline
$1 \le \beta_3 \le 2, \quad \beta_1 = \beta_4  = \gamma_2= 1$ & 2  \\
\hline
$1 \le \beta_1 \le 2, \quad \beta_3  = \gamma_1= 1$ & 2  \\
\hline
$1 \le \beta_1, \gamma_2 \le 2, \quad \beta_3 = \gamma_1 = 1$ & 4  \\
\hline
$1 \le \beta_1, \beta_3 \le 2, \mbox{ excluding } \beta_1 = \beta_3  = 2, \quad \gamma_2 = 1 $ & 3  \\
\hline
$ \beta_1 = \beta_4 = \gamma_1 = 1$ & 1  \\
\hline
$1 \le  \gamma_1,  \gamma_2 \le 2, \mbox{ excluding } \gamma_1=  \gamma_2=2, \quad \beta_1 = \beta_4 = 1  $ & 3  \\
\hline
$ \beta_1 = \beta_4 = \gamma_2 = 1$ & 1  \\
\hline
$\beta_1 = \gamma_2 = 1$ & 1  \\
\hline
$1 \le  \gamma_1,  \gamma_2 \le 2, \mbox{ excluding } \gamma_1=  \gamma_2=2, \quad  \beta_1 = 1 $ & 3  \\
\hline
$1 \le \beta_2, \beta_4 \le 2, \quad \beta_3 = \gamma_1 = 1$ & 4  \\
\hline
$1 \le  \beta_4 \le 2, \quad \beta_2=\beta_3 = \gamma_1 = \gamma_2 =1$ & 2  \\
\hline
$1 \le \beta_3, \beta_4 \le 2, \quad \beta_2 = \gamma_2 = 1$ & 4  \\
\hline
$1 \le \beta_2, \beta_4 \le 2, \mbox{ excluding } \beta_2 = \beta_4  = 2, \quad  \gamma_1 = 1 $ & 3  \\
\hline
$1 \le \beta_4, \gamma_1 \le 2, \quad  \beta_2 = \gamma_2 = 1$ & 4  \\
\hline
$1 \le \beta_4 \le 2, \quad  \beta_2 = \gamma_2 = 1$ & 2  \\
\hline
$   \beta_2= \gamma_1 = 1$ & 1  \\
\hline
$1 \le \gamma_1 \le 2, \quad   \beta_2 =\gamma_2 = 1$ & 2  \\
\hline
$ \beta_2 =\gamma_2 = 1$ & 1  \\
\hline
$1 \le \beta_4 \le 2, \quad  \beta_3 = \gamma_1 = 1$ & 2  \\
\hline
$1 \le \beta_4, \gamma_2 \le 2, \quad  \beta_3 = \gamma_1 = 1$ & 4  \\
\hline
$1 \le  \beta_3,  \beta_4 \le 2, \mbox{ excluding } \beta_3 = \beta_4  = 2, \quad  \gamma_2 = 1 $ & 3  \\
\hline
$ \beta_3 =\gamma_1 = 1$ & 1  \\
\hline
$1 \le \gamma_2 \le 2, \quad   \beta_3 =\gamma_1 = 1$ & 2  \\
\hline
$ \beta_3 = \gamma_2 = 1$ & 1  \\
\hline
$ \beta_4 = \gamma_1 = 1$ & 1  \\
\hline
$1 \le  \gamma_1, \gamma_2 \le 2, \mbox{ excluding }  \gamma_1 = \gamma_2 = 2,  \quad \beta_4 =1 $ & 3  \\
\hline
\end{tabular}
}
\vspace{.2cm}
\caption{The $92$ ODE-classes of connected $3$-node REI networks with valence $\le 2$ without autoregulation having both excitatory and inhibitory arrows.  
See Figure~\ref{fig:generalREI_3node}.
}
\label{table:ODE_3node_valence_up_2_a}
\end{table}

\begin{prop} 
\label{P:4.3}
Any connected $3$-node REI network with valence $\leq 2$ is ODE-equivalent to
the network in {\rm Figure}~{\rm \ref{fig:generalREI_3node}}, where, under minimality, $\delta = \tau=0$ and

{\rm (a)} If there is no autoregulation, the nonzero arrow multiplicities $\beta_i$ $(i=1,2,3,4)$ and $\gamma_j$ $(j=1,2)$ appear in {\rm Tables}~{\rm \ref{table:ODE_3node_valence_up_2_a}} and {\rm \ref{table:ODE_3node_valence_up_2_b}}.

{\rm (b)} If there is autoregulation, the nonzero arrow multiplicities $\alpha$, $\beta_i$ $(i=1,2,3,4)$ and $\gamma_j$ $(j=1,2)$, appear in {\rm Tables}~{\rm \ref{table:ODE_3node_valence_up_2_c}} and {\rm \ref{table:ODE_3node_valence_up_2_d}}.
\end{prop}

\begin{proof}
The result follows from Propositions~\ref{prop:3_node_connected} and \ref{prop:3_node_ODE}, since a $3$-node REI network with valence $\leq 2$ must satisfy
$$
0 \le \alpha + \beta_3 +\gamma_1 \le 2, \quad 0 \le \tau + \beta_1 + \beta_4 \le 2 \quad \mbox{and} \quad 0 \le \delta + \beta_2 +\gamma_2 \le 2.
$$
\end{proof}

\begin{table}[htb]
{\tiny 
\begin{tabular}{|l|l|}
\hline 
{ nonzero multiplicities}
 & 
{$\#$ ODE-classes}
\\
 \hline
$1 \le  \beta_1,  \beta_2 \le 2, \mbox{ excluding } \beta_1=  \beta_2=2 $ & 3  \\
 \hline 
$1 \le  \beta_1,  \beta_2,  \beta_3 \le 2, \mbox{ excluding } \beta_1=  \beta_2=  \beta_3=2 $ & 7  \\
 \hline 
$1 \le  \beta_2 \le 2, \quad \beta_1=  \beta_4=1 $ & 2  \\
 \hline 
$1 \le  \beta_2,  \beta_3 \le 2, \quad \beta_1=  \beta_4=1 $ & 4  \\
 \hline 
$1 \le  \beta_2,  \beta_3, \beta_4 \le 2, \mbox{ excluding } \beta_2=  \beta_3=  \beta_4=2 $ & 7  \\
 \hline 
$1 \le  \beta_2,  \beta_4 \le 2, \mbox{ excluding } \beta_2=  \beta_4=2 $ & 3  \\
 \hline
$1 \le  \beta_1,  \beta_3 \le 2, \mbox{ excluding } \beta_1=  \beta_3=2 $ & 3  \\
 \hline 
$1 \le   \beta_3 \le 2, \quad \beta_1=  \beta_4=1 $ & 2  \\
 \hline 
$1 \le  \beta_3,  \beta_4 \le 2, \mbox{ excluding } \beta_3=  \beta_4=2 $ & 3  \\
 \hline
$ \beta_1=  \beta_4=1 $ & 1  \\
 \hline 
$1 \le  \gamma_1,  \gamma_2 \le 2, \mbox{ excluding } \gamma_1=  \gamma_2=2 $ & 3  \\
 \hline
\end{tabular}
}
\vspace{.2cm}
\caption{The $38$ ODE-classes of connected REI networks with valence $\le 2$ without autoregulation having only excitatory or inhibitory arrows. 
See Figure~\ref{fig:generalREI_3node}
}
\label{table:ODE_3node_valence_up_2_b}
\end{table}

\begin{lemma}
If $\mathcal{G}$ is a  connected $3$-node REI network with input valence $\leq 2$, where nodes $1,2$ are of type $N^E$ and node $3$ is of type $N^I$, then 
the subnetwork of $\mathcal{G}$ containing nodes $2,3$ 
and all arrows between these two nodes
is a $2$-node REI  network with input valence $\leq 2$. 
\label{lemma:strategy_2_to_3}
\end{lemma}
\begin{proof}
The subnetwork $S$ of $\mathcal{G}$ containing nodes $2,3$ is a $2$-node network 
where node $2$ is of type $N^E$ and node $3$ is of type $N^I$. Since 
$\mathcal{G}$ is REI then node $2$ outputs only excitatory arrows and 
node $3$ outputs only inhibitory arrows. Therefore $S$ is also an REI network. 
\end{proof}

\begin{table}[!ht]
{\tiny 
\begin{tabular}{|l|l|}
\hline 
{nonzero multiplicities}
 & 
{$\#$ ODE-classes}
\\
\hline 
\hline  
$1 \le \beta_1 \le 2, \quad \beta_3 = \alpha = 1, \beta_2 = \gamma_2 = 1$ & 2  \\
\hline
$\beta_3 = \alpha=1,  \beta_1=\beta_2 = \beta_4  = \gamma_2= 1$ & 1  \\
\hline
$1 \le \beta_2 \le 2, \quad \gamma_1  = \alpha=1, \beta_1 = \beta_4 = 1$ & 2  \\
\hline
$1 \le \alpha \le 2, \quad\beta_1=\beta_2  = \beta_4  = \gamma_2= 1$ & 2  \\
\hline
$\gamma_1 = \alpha=1, \beta_1 = \beta_2 = \beta_4 = \gamma_2 = 1$ & 2  \\
\hline
$1 \le \beta_1, \beta_2,  \le 2,  \quad \gamma_1 =\alpha = 1 $ & 4  \\
\hline
$1 \le \beta_1 \le 2, \quad \gamma_1 = \alpha=1, \beta_2 = \gamma_2 = 1$ & 2  \\
\hline
$1 \le \alpha, \beta_1 \le 2, \quad \beta_2 = \gamma_2 = 1$ & 4  \\
\hline
$\beta_3 =\alpha =1, \beta_1 = \beta_4  = \gamma_2= 1$ & 1  \\
\hline
$1 \le \beta_1 \le 2, \quad \beta_3 = \alpha=1,  \gamma_2 = 1 $ & 2  \\
\hline
$ \gamma_1 = \alpha = 1, \beta_1 = \beta_4 = 1$ & 1  \\
\hline
$1 \le  \gamma_2 \le 2, \quad  \gamma_1 = \alpha =1, \beta_1 = \beta_4 = 1  $ & 2  \\
\hline
$1 \le  \alpha \le 2, \quad \beta_1 = \beta_4 = \gamma_2 = 1$ & 2  \\
\hline
$1 \le  \alpha, \beta_1 \le 2,  \mbox{ excluding } \alpha=  \beta_1=2 \quad \gamma_2 = 1$ & 3 \\
\hline
$1 \le  \beta_1,  \gamma_2 \le 2, \quad  \gamma_1 = \alpha= 1 $ & 4  \\
\hline
$1 \le \beta_4 \le 2, \quad \beta_3 = \alpha=1, \beta_2 = \gamma_2 = 1$ & 2  \\
\hline
$1 \le \beta_2, \beta_4 \le 2, \quad  \gamma_1 = \alpha = 1 $ & 4  \\
\hline
$1 \le \beta_4 \le 2, \quad  \gamma_1 = \alpha=1, \beta_2 = \gamma_2 = 1$ & 2  \\
\hline
$1 \le \alpha, \beta_4 \le 2, \quad  \beta_2 = \gamma_2 = 1$ & 4  \\
\hline
$  1 \le \beta_2 \le 2, \quad  \gamma_1 = \alpha = 1$ & 2  \\
\hline
$\gamma_1 =\alpha = 1,  \beta_2 =\gamma_2 = 1$ & 1  \\
\hline
$ 1 \le \alpha \le 2, \quad\beta_2 =\gamma_2 = 1$ & 2  \\
\hline
$1 \le  \beta_4 \le 2, \quad  \beta_3 = \alpha = 1, \gamma_2 = 1 $ & 2  \\
\hline
$ \beta_3 = \alpha =1,  \gamma_2 = 1$ & 1  \\
\hline
$ 1 \le  \beta_4 \le 2, \quad  \gamma_1 = \alpha = 1$ & 2  \\
\hline
$1 \le  \beta_4, \gamma_2 \le 2,  \quad \gamma_1 = \alpha  =1 $ & 4  \\
\hline
$1 \le   \gamma_2 \le 2, \quad  \gamma_1 = \alpha = 1 $ & 2  \\
 \hline
\end{tabular}
}
\vspace{.2cm}
\caption{The $62$ ODE-classes of connected REI networks with valence $\le 2$ with autoregulation having both excitatory and inhibitory arrows. 
See Figure~\ref{fig:generalREI_3node}.
}
\label{table:ODE_3node_valence_up_2_c}
\end{table}

\begin{table}[!ht]
{\tiny 
\begin{tabular}{|l|l|}
\hline 
{\footnotesize nonzero multiplicities}
 & 
{$\#$  ODE-classes}
\\
 \hline
$1 \le  \alpha, \beta_1,  \beta_2 \le 2, \mbox{ excluding } \alpha=\beta_1=  \beta_2 = 2 $ & 7  \\
 \hline 
$1 \le  \beta_1,  \beta_2 \le 2, \quad \beta_3 = \alpha=1  $ & 4  \\
 \hline 
$1 \le \alpha, \beta_2 \le 2, \quad \beta_1=  \beta_4=1 $ & 4  \\
 \hline 
$1 \le  \beta_2 \le 2, \quad \beta_3 = \alpha =1, \beta_1=  \beta_4=1 $ & 2  \\
 \hline 
$1 \le  \beta_2,   \beta_4 \le 2, \quad \beta_3 = \alpha = 1 $ & 4  \\
 \hline 
$1 \le  \alpha, \beta_2,  \beta_4 \le 2, \mbox{ excluding } \alpha = \beta_2=  \beta_4=2 $ & 7  \\
 \hline
$1 \le  \beta_1 \le 2, \quad   \beta_3 = \alpha = 1 $ & 2  \\
 \hline 
$  \beta_3 = \alpha = 1,  \beta_1=  \beta_4=1 $ & 1  \\
 \hline 
$1 \le   \beta_4 \le 2, \quad \beta_3 = \alpha =1 $ & 2  \\
 \hline
$ 1 \le   \alpha \le 2, \quad \beta_1=  \beta_4=1 $ & 2  \\
 \hline 
\end{tabular}
}
\vspace{.2cm}
\caption{The $35$ ODE-classes of connected REI networks with valence $\le 2$ with autoregulation having only excitatory or inhibitory arrows.
See Figure~\ref{fig:generalREI_3node}.
}
\label{table:ODE_3node_valence_up_2_d}
\end{table}

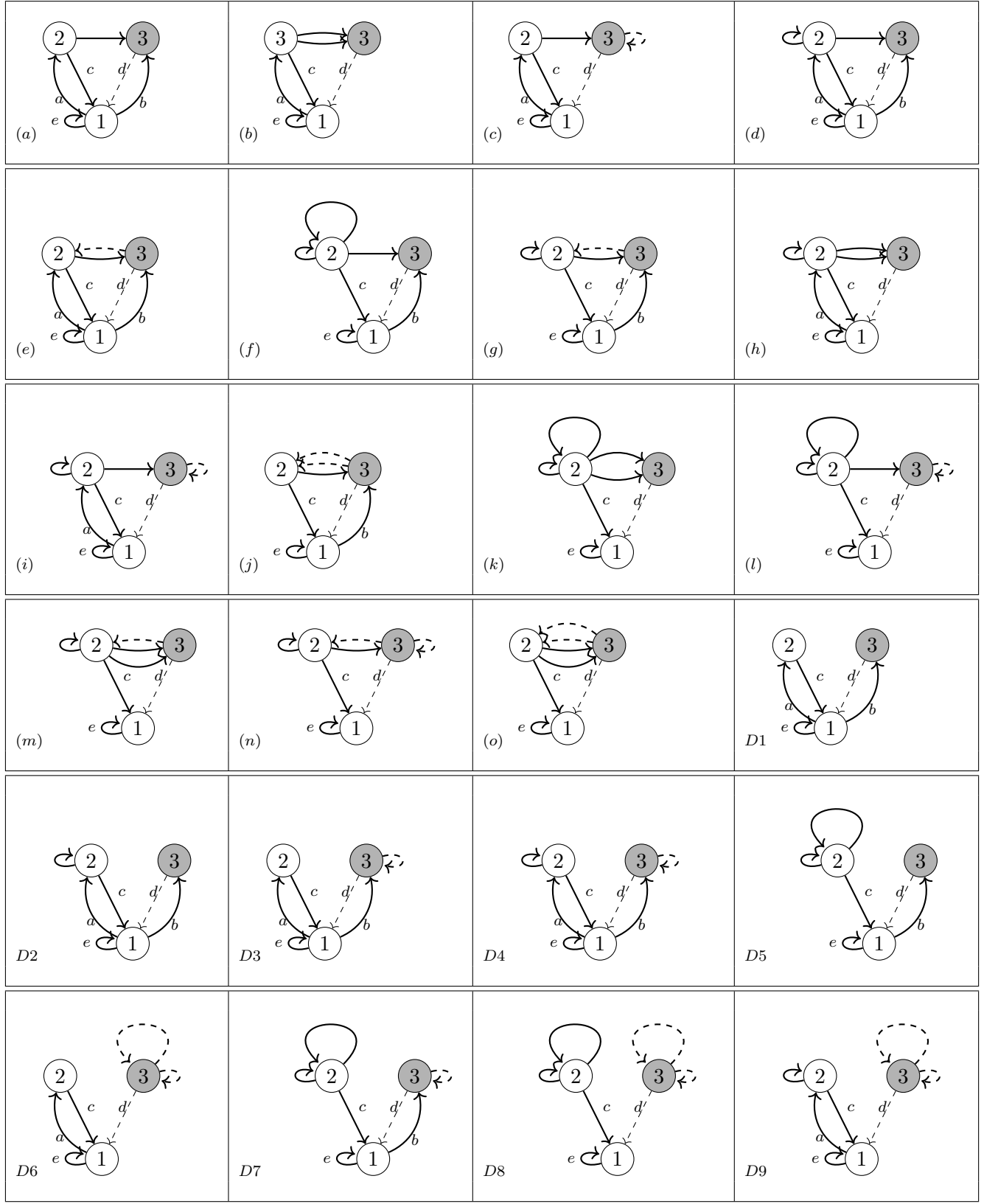
\begin{figure}[]
\begin{center}
{\tiny 
\begin{tabular}{|l|l|l|l|} 
 \hline 
  & & &  \\
$(a)$ \begin{tikzpicture}
 [scale=.15,auto=left, node distance=1.5cm, 
 ]
 \node[fill=white,style={circle,draw}] (n0) at (9,-10) {\small{1}};
 \node[fill=white,style={circle,draw}] (n1) at (4,0) {\small{2}};
  \node[fill=black!30,style={circle,draw}] (n2) at (14,0) {\small{3}};
 \draw[->, thick] (n1) edge[thick] node {}  (n2); 
  \path (n1) [->,thick] edge[] node {{\tiny $c$}} (n0);
 \path (n2) [->,dashed] edge[above] node {{\tiny $d$}} (n0);
  \path (n0) [->,thick] edge[loop left=90] node {{\tiny $e$}} (n0);   
  \path  (n0) [->,thick] edge[bend left=40,below] node {{\tiny $a$}} (n1);   
   \path       (n0) [->,thick] edge[bend right=40,below] node {{\tiny $b$}} (n2);  
 \end{tikzpicture}  
 &  
$(b)$ \begin{tikzpicture}
 [scale=.15,auto=left, node distance=1.5cm, 
 ]
 \node[fill=white,style={circle,draw}] (n0) at (9,-10) {\small{1}};
 \node[fill=white,style={circle,draw}] (n1) at (4,0) {\small{3}};
  \node[fill=black!30,style={circle,draw}] (n2) at (14,0) {\small{3}};
 \draw[->, thick] (n1) edge[bend right=10, thick] node {}  (n2); 
 \draw[->, thick] (n1) edge[bend right=-10,thick] node {}  (n2); 
  \path (n1) [->,thick] edge[] node {{\tiny $c$}} (n0);
 \path (n2) [->,dashed] edge[above] node {{\tiny $d$}} (n0);
  \path (n0) [->,thick] edge[loop left=90] node {{\tiny $e$}} (n0);   
  \path  (n0) [->,thick] edge[bend left=40,below] node {{\tiny $a$}} (n1);   
 \end{tikzpicture}  
  &  
$(c)$ \begin{tikzpicture}
 [scale=.15,auto=left, node distance=1.5cm,
 ]
 \node[fill=white,style={circle,draw}] (n0) at (9,-10) {\small{1}};
 \node[fill=white,style={circle,draw}] (n1) at (4,0) {\small{2}};
  \node[fill=black!30,style={circle,draw}] (n2) at (14,0) {\small{3}};
\path
        (n1) [->,solid, thick]  edge[thick] node { } (n2);
\path        
        (n2) [->,dashed, thick]  edge[loop right=90,thick] node { } (n2);
 \path (n1) [->,thick] edge[] node {{\tiny $c$}} (n0);
 \path (n2) [->,dashed] edge[above] node {{\tiny $d$}} (n0);
  \path (n0) [->,thick] edge[loop left=90] node {{\tiny $e$}} (n0);   
  \path  (n0) [->,thick] edge[bend left=40,below] node {{\tiny $a$}} (n1);   
 \end{tikzpicture} 
  &  
$(d)$ \begin{tikzpicture}
 [scale=.15,auto=left, node distance=1.5cm, 
 ]
 \node[fill=white,style={circle,draw}] (n0) at (9,-10) {\small{1}};
 \node[fill=white,style={circle,draw}] (n1) at (4,0) {\small{2}};
  \node[fill=black!30,style={circle,draw}] (n2) at (14,0) {\small{3}};
\path
        (n1) [->,thick]  edge[thick] node { } (n2)
        (n1) [->,thick]  edge[loop left=90,thick] node { } (n1);
         \path (n1) [->,thick] edge[] node {{\tiny $c$}} (n0);
 \path (n2) [->,dashed] edge[above] node {{\tiny $d$}} (n0);
  \path (n0) [->,thick] edge[loop left=90] node {{\tiny $e$}} (n0);   
  \path  (n0) [->,thick] edge[bend left=40,below] node {{\tiny $a$}} (n1);   
   \path       (n0) [->,thick] edge[bend right=40,below] node {{\tiny $b$}} (n2);  
 \end{tikzpicture} \\
 & & &  \\
 \hline 
 \hline 
  & & &  \\
$(e)$ \begin{tikzpicture}
 [scale=.15,auto=left, node distance=1.5cm,
 ]
 \node[fill=white,style={circle,draw}] (n0) at (9,-10) {\small{1}};
 \node[fill=white,style={circle,draw}] (n1) at (4,0) {\small{2}};
  \node[fill=black!30,style={circle,draw}] (n2) at (14,0) {\small{3}};
\path
        (n1) [->,thick]  edge[bend right=10] node { } (n2);
 \path 
        (n2) [->,dashed, thick]  edge[bend left=-10] node { } (n1);
         \path (n1) [->,thick] edge[] node {{\tiny $c$}} (n0);
 \path (n2) [->,dashed] edge[above] node {{\tiny $d$}} (n0);
  \path (n0) [->,thick] edge[loop left=90] node {{\tiny $e$}} (n0);   
  \path  (n0) [->,thick] edge[bend left=40,below] node {{\tiny $a$}} (n1);   
   \path       (n0) [->,thick] edge[bend right=40,below] node {{\tiny $b$}} (n2);  
 \end{tikzpicture} 
&   
 $(f)$ \begin{tikzpicture}
 [scale=.15,auto=left, node distance=1.5cm, 
 ]
 \node[fill=white,style={circle,draw}] (n0) at (9,-10) {\small{1}};
 \node[fill=white,style={circle,draw}] (n1) at (4,0) {\small{2}};
  \node[fill=black!30,style={circle,draw}] (n2) at (14,0) {\small{3}};
\path
        (n1)  [->]  edge[loop,thick] node {} (n1)
        (n1)  [->]  edge[loop left=90,thick] node {} (n1)
        (n1) [->,thick]  edge[thick] node { } (n2);
  \path (n1) [->,thick] edge[] node {{\tiny $c$}} (n0);
 \path (n2) [->,dashed] edge[above] node {{\tiny $d$}} (n0);
  \path (n0) [->,thick] edge[loop left=90] node {{\tiny $e$}} (n0);   
  \path       (n0) [->,thick] edge[bend right=40,below] node {{\tiny $b$}} (n2);     
 \end{tikzpicture}
  &  
 $(g)$ \begin{tikzpicture}
 [scale=.15,auto=left, node distance=1.5cm, 
 ]
 \node[fill=white,style={circle,draw}] (n0) at (9,-10) {\small{1}};
 \node[fill=white,style={circle,draw}] (n1) at (4,0) {\small{2}};
  \node[fill=black!30,style={circle,draw}] (n2) at (14,0) {\small{3}};
\path
        (n1)  [->]  edge[loop left=90,thick] node {} (n1)
        (n1) [->,thick]  edge[bend right=10, thick] node { } (n2);
\path 
        (n2) [->,dashed, thick]  edge[bend left=-10] node { } (n1);
        \path (n1) [->,thick] edge[] node {{\tiny $c$}} (n0);
 \path (n2) [->,dashed] edge[above] node {{\tiny $d$}} (n0);
  \path (n0) [->,thick] edge[loop left=90] node {{\tiny $e$}} (n0);   
  \path       (n0) [->,thick] edge[bend right=40,below] node {{\tiny $b$}} (n2);   
 \end{tikzpicture}
&   
 $(h)$  \begin{tikzpicture}
 [scale=.15,auto=left, node distance=1.5cm, 
 ]
 \node[fill=white,style={circle,draw}] (n0) at (9,-10) {\small{1}};
 \node[fill=white,style={circle,draw}] (n1) at (4,0) {\small{2}};
  \node[fill=black!30,style={circle,draw}] (n2) at (14,0) {\small{3}};
\path
        (n1)  [->]  edge[loop left=90,thick] node {} (n1)
        (n1) [->,thick]  edge[bend right=10, thick] node { } (n2)
        (n1) [->,thick]  edge[bend right=-10, thick] node { } (n2);
  \path (n1) [->,thick] edge[] node {{\tiny $c$}} (n0);
 \path (n2) [->,dashed] edge[above] node {{\tiny $d$}} (n0);
  \path (n0) [->,thick] edge[loop left=90] node {{\tiny $e$}} (n0);   
  \path  (n0) [->,thick] edge[bend left=40,below] node {{\tiny $a$}} (n1);   
 \end{tikzpicture}\\
 & & &  \\
 \hline 
 \hline 
  & & &  \\
$(i)$ \begin{tikzpicture}
 [scale=.15,auto=left, node distance=1.5cm, 
 ]
 \node[fill=white,style={circle,draw}] (n0) at (9,-10) {\small{1}};
 \node[fill=white,style={circle,draw}] (n1) at (4,0) {\small{2}};
  \node[fill=black!30,style={circle,draw}] (n2) at (14,0) {\small{3}};
\path
        (n1) [->,thick]  edge[thick] node { } (n2)
         (n1) [->,thick]  edge[loop left=90, thick] node { } (n1);
\path 
        (n2) [->,dashed, thick]  edge[loop right=90, thick] node { } (n2);
 \path (n1) [->,thick] edge[] node {{\tiny $c$}} (n0);
 \path (n2) [->,dashed] edge[above] node {{\tiny $d$}} (n0);
  \path (n0) [->,thick] edge[loop left=90] node {{\tiny $e$}} (n0);   
  \path  (n0) [->,thick] edge[bend left=40,below] node {{\tiny $a$}} (n1);   
 \end{tikzpicture} 
  &   
 $(j)$ \begin{tikzpicture}
 [scale=.15,auto=left, node distance=1.5cm, 
 ]
 \node[fill=white,style={circle,draw}] (n0) at (9,-10) {\small{1}};
 \node[fill=white,style={circle,draw}] (n1) at (4,0) {\small{2}};
  \node[fill=black!30,style={circle,draw}] (n2) at (14,0) {\small{3}};
\path
         (n1) [->,thick]  edge[bend right=10, thick] node { } (n2);
\path 
        (n2) [->,thick, dashed]  edge[bend left =-10] node { } (n1)
        (n2) [->,thick, dashed]  edge[bend left=-30] node { } (n1);
 \path (n1) [->,thick] edge[] node {{\tiny $c$}} (n0);
 \path (n2) [->,dashed] edge[above] node {{\tiny $d$}} (n0);
  \path (n0) [->,thick] edge[loop left=90] node {{\tiny $e$}} (n0);   
  \path       (n0) [->,thick] edge[bend right=40,below] node {{\tiny $b$}} (n2);     
 \end{tikzpicture}
 &   
 $(k)$ 
  \begin{tikzpicture}
 [scale=.15,auto=left, node distance=1.5cm,
 ]
 \node[fill=white,style={circle,draw}] (n0) at (9,-10) {\small{1}};
 \node[fill=white,style={circle,draw}] (n1) at (4,0) {\small{2}};
  \node[fill=black!30,style={circle,draw}] (n2) at (14,0) {\small{3}};
\path
        (n1)  [->]  edge[loop,thick] node {} (n1)
        (n1)  [->]  edge[loop left=90,thick] node {} (n1)
        (n1) [->,thick]  edge[bend right=20, thick] node { } (n2)
        (n1) [->,thick]  edge[bend right=-30, thick] node { } (n2);
        \path (n1) [->,thick] edge[] node {{\tiny $c$}} (n0);
   \path (n2) [->,dashed] edge[above] node {{\tiny $d$}} (n0);
    \path (n0) [->,thick] edge[loop left=90] node {{\tiny $e$}} (n0);  
 \end{tikzpicture}
&  
 $(l)$  \begin{tikzpicture}
 [scale=.15,auto=left, node distance=1.5cm, 
 ]
 \node[fill=white,style={circle,draw}] (n0) at (9,-10) {\small{1}};
 \node[fill=white,style={circle,draw}] (n1) at (4,0) {\small{2}};
  \node[fill=black!30,style={circle,draw}] (n2) at (14,0) {\small{3}};
\path
        (n1)  [->]  edge[loop,thick] node {} (n1)
        (n1)  [->]  edge[loop left=90,thick] node {} (n1)
        (n1) [->,thick]  edge[thick] node { } (n2);
\path 
         (n2)  [->, dashed]  edge[loop right=90,thick] node {} (n2); 
 \path (n1) [->,thick] edge[] node {{\tiny $c$}} (n0);
   \path (n2) [->,dashed] edge[above] node {{\tiny $d$}} (n0);
    \path (n0) [->,thick] edge[loop left=90] node {{\tiny $e$}} (n0);                     
 \end{tikzpicture}\\
  & & &  \\
 \hline 
 \hline 
  & & &  \\
 $(m)$   \begin{tikzpicture}
 [scale=.15,auto=left, node distance=1.5cm, 
 ]
 \node[fill=white,style={circle,draw}] (n0) at (9,-10) {\small{1}};
 \node[fill=white,style={circle,draw}] (n1) at (4,0) {\small{2}};
  \node[fill=black!30,style={circle,draw}] (n2) at (14,0) {\small{3}};
\path
        (n1)  [->]  edge[loop left=90,thick] node {} (n1)
        (n1) [->,thick]  edge[bend right=10, thick] node { } (n2)
         (n1) [->,thick]  edge[bend right=40, thick] node { } (n2);
\path 
        (n2) [->,thick, dashed]  edge[bend left=-10, thick] node { } (n1);   
 \path (n1) [->,thick] edge[] node {{\tiny $c$}} (n0);
   \path (n2) [->,dashed] edge[above] node {{\tiny $d$}} (n0);
    \path (n0) [->,thick] edge[loop left=90] node {{\tiny $e$}} (n0);             
 \end{tikzpicture}
 &  
 $(n)$  \begin{tikzpicture}
 [scale=.15,auto=left, node distance=1.5cm,
 ]
 \node[fill=white,style={circle,draw}] (n0) at (9,-10) {\small{1}};
 \node[fill=white,style={circle,draw}] (n1) at (4,0) {\small{2}};
  \node[fill=black!30,style={circle,draw}] (n2) at (14,0) {\small{3}};
\path
        (n1)  [->]  edge[loop left=90,thick] node {} (n1)
         (n1) [->,thick]  edge[bend right=10] node { } (n2);
\path          
         (n2)  [->, dashed]  edge[loop right=90,thick] node {} (n2)       
        (n2) [->,thick, dashed]  edge[bend left=-10] node { } (n1);
  \path (n1) [->,thick] edge[] node {{\tiny $c$}} (n0);
   \path (n2) [->,dashed] edge[above] node {{\tiny $d$}} (n0);
    \path (n0) [->,thick] edge[loop left=90] node {{\tiny $e$}} (n0);              
 \end{tikzpicture} 
 & 
 $(o)$ \begin{tikzpicture}
 [scale=.15,auto=left, node distance=1.5cm,
 ]
 \node[fill=white,style={circle,draw}] (n0) at (9,-10) {\small{1}};
 \node[fill=white,style={circle,draw}] (n1) at (4,0) {\small{2}};
  \node[fill=black!30,style={circle,draw}] (n2) at (14,0) {\small{3}};
\path
        (n1)  [->]  edge[bend right=10,thick] node {} (n2)
         (n1) [->,thick]  edge[bend right=40, thick] node { } (n2);
 \path 
         (n2)  [->,dashed]  edge[bend left=-10,thick] node {} (n1)      
        (n2) [->,thick,dashed]  edge[bend left=-40, thick] node { } (n1);
        \path (n1) [->,thick] edge[] node {{\tiny $c$}} (n0);
   \path (n2) [->,dashed] edge[above] node {{\tiny $d$}} (n0);
    \path (n0) [->,thick] edge[loop left=90] node {{\tiny $e$}} (n0);        
 \end{tikzpicture} & 
$D1$ \begin{tikzpicture}
 [scale=.15,auto=left, node distance=1.5cm, 
 ]
 \node[fill=white,style={circle,draw}] (n0) at (9,-10) {\small{1}};
 \node[fill=white,style={circle,draw}] (n1) at (4,0) {\small{2}};
  \node[fill=black!30,style={circle,draw}] (n2) at (14,0) {\small{3}};
   \path (n1) [->,thick] edge[] node {{\tiny $c$}} (n0);
 \path (n2) [->,dashed] edge[above] node {{\tiny $d$}} (n0);
  \path (n0) [->,thick] edge[loop left=90] node {{\tiny $e$}} (n0);   
  \path  (n0) [->,thick] edge[bend left=40,below] node {{\tiny $a$}} (n1);   
   \path       (n0) [->,thick] edge[bend right=40,below] node {{\tiny $b$}} (n2);  
 \end{tikzpicture}  \\
  & & &  \\
 \hline 
 \hline 
  & & &  \\
$D2$ \begin{tikzpicture}
 [scale=.15,auto=left, node distance=1.5cm, 
 ]
 \node[fill=white,style={circle,draw}] (n0) at (9,-10) {\small{1}};
 \node[fill=white,style={circle,draw}] (n1) at (4,0) {\small{2}};
  \node[fill=black!30,style={circle,draw}] (n2) at (14,0) {\small{3}};
\path
 (n1) [->,thick] edge[loop left=90,thick] node {} (n1);
  \path (n1) [->,thick] edge[] node {{\tiny $c$}} (n0);
 \path (n2) [->,dashed] edge[above] node {{\tiny $d$}} (n0);
  \path (n0) [->,thick] edge[loop left=90] node {{\tiny $e$}} (n0);   
  \path  (n0) [->,thick] edge[bend left=40,below] node {{\tiny $a$}} (n1);   
   \path       (n0) [->,thick] edge[bend right=40,below] node {{\tiny $b$}} (n2);  
 \end{tikzpicture}   
 & 
$D3$ \begin{tikzpicture}
 [scale=.15,auto=left, node distance=1.5cm, 
 ]
 \node[fill=white,style={circle,draw}] (n0) at (9,-10) {\small{1}};
 \node[fill=white,style={circle,draw}] (n1) at (4,0) {\small{2}};
  \node[fill=black!30,style={circle,draw}] (n2) at (14,0) {\small{3}};
  \path 
  (n2) [->,dashed] edge[loop right=90,thick] node {} (n2);
   \path (n1) [->,thick] edge[] node {{\tiny $c$}} (n0);
 \path (n2) [->,dashed] edge[above] node {{\tiny $d$}} (n0);
  \path (n0) [->,thick] edge[loop left=90] node {{\tiny $e$}} (n0);   
  \path  (n0) [->,thick] edge[bend left=40,below] node {{\tiny $a$}} (n1);   
   \path       (n0) [->,thick] edge[bend right=40,below] node {{\tiny $b$}} (n2);  
 \end{tikzpicture} &
$D4$ \begin{tikzpicture}
 [scale=.15,auto=left, node distance=1.5cm, 
 ]
 \node[fill=white,style={circle,draw}] (n0) at (9,-10) {\small{1}};
 \node[fill=white,style={circle,draw}] (n1) at (4,0) {\small{2}};
  \node[fill=black!30,style={circle,draw}] (n2) at (14,0) {\small{3}};
\path
 (n1) [->,thick] edge[loop left=90,thick] node {} (n1);
 \path 
(n2) [->,dashed] edge[loop right=90,thick] node {} (n2);
 \path (n1) [->,thick] edge[] node {{\tiny $c$}} (n0);
 \path (n2) [->,dashed] edge[above] node {{\tiny $d$}} (n0);
  \path (n0) [->,thick] edge[loop left=90] node {{\tiny $e$}} (n0);   
  \path  (n0) [->,thick] edge[bend left=40,below] node {{\tiny $a$}} (n1);   
   \path       (n0) [->,thick] edge[bend right=40,below] node {{\tiny $b$}} (n2);  
 \end{tikzpicture} 
  &  
  $D5$ \begin{tikzpicture}
 [scale=.15,auto=left, node distance=1.5cm, 
 ]
 \node[fill=white,style={circle,draw}] (n0) at (9,-10) {\small{1}};
 \node[fill=white,style={circle,draw}] (n1) at (4,0) {\small{2}};
  \node[fill=black!30,style={circle,draw}] (n2) at (14,0) {\small{3}};
\path
        (n1)  [->]  edge[loop,thick] node {} (n1)
        (n1)  [->]  edge[loop left=90,thick] node {} (n1);
         \path (n1) [->,thick] edge[] node {{\tiny $c$}} (n0);
 \path (n2) [->,dashed] edge[above] node {{\tiny $d$}} (n0);
  \path (n0) [->,thick] edge[loop left=90] node {{\tiny $e$}} (n0);   
   \path       (n0) [->,thick] edge[bend right=40,below] node {{\tiny $b$}} (n2);  
 \end{tikzpicture}   \\
 & & &  \\
 \hline 
 \hline 
  & & &  \\
 $D6$ \begin{tikzpicture}
 [scale=.15,auto=left, node distance=1.5cm, 
 ]
 \node[fill=white,style={circle,draw}] (n0) at (9,-10) {\small{1}};
 \node[fill=white,style={circle,draw}] (n1) at (4,0) {\small{2}};
  \node[fill=black!30,style={circle,draw}] (n2) at (14,0) {\small{3}};
\path
        (n2)  [->,dashed]  edge[loop,thick] node {} (n2)
        (n2)  [->,dashed]  edge[loop right=90,thick] node {} (n2);
        \path (n1) [->,thick] edge[] node {{\tiny $c$}} (n0);
 \path (n2) [->,dashed] edge[above] node {{\tiny $d$}} (n0);
  \path (n0) [->,thick] edge[loop left=90] node {{\tiny $e$}} (n0);   
  \path  (n0) [->,thick] edge[bend left=40,below] node {{\tiny $a$}} (n1);   
 \end{tikzpicture}   &
 $D7$ \begin{tikzpicture}
 [scale=.15,auto=left, node distance=1.5cm, 
 ]
 \node[fill=white,style={circle,draw}] (n0) at (9,-10) {\small{1}};
 \node[fill=white,style={circle,draw}] (n1) at (4,0) {\small{2}};
  \node[fill=black!30,style={circle,draw}] (n2) at (14,0) {\small{3}};
\path
  (n1)  [->]  edge[loop,thick] node {} (n1)
        (n1)  [->]  edge[loop left=90,thick] node {} (n1);
 \path 
        (n2)  [->,dashed]  edge[loop right=90,thick] node {} (n2);
        \path (n1) [->,thick] edge[] node {{\tiny $c$}} (n0);
 \path (n2) [->,dashed] edge[above] node {{\tiny $d$}} (n0);
  \path (n0) [->,thick] edge[loop left=90] node {{\tiny $e$}} (n0);   
   \path       (n0) [->,thick] edge[bend right=40,below] node {{\tiny $b$}} (n2);  
 \end{tikzpicture} 
   & 
  $D8$ \begin{tikzpicture}
 [scale=.15,auto=left, node distance=1.5cm, 
 ]
 \node[fill=white,style={circle,draw}] (n0) at (9,-10) {\small{1}};
 \node[fill=white,style={circle,draw}] (n1) at (4,0) {\small{2}};
  \node[fill=black!30,style={circle,draw}] (n2) at (14,0) {\small{3}};
\path 
        (n1)  [->]  edge[loop,thick] node {} (n1)
        (n1)  [->]  edge[loop left=90,thick] node {} (n1);
\path
        (n2)  [->,dashed]  edge[loop,thick] node {} (n2)
        (n2)  [->,dashed]  edge[loop right=90,thick] node {} (n2);
         \path (n1) [->,thick] edge[] node {{\tiny $c$}} (n0);
 \path (n2) [->,dashed] edge[above] node {{\tiny $d$}} (n0);
  \path (n0) [->,thick] edge[loop left=90] node {{\tiny $e$}} (n0);   
 \end{tikzpicture} 
  & 
  $D9$ \begin{tikzpicture}
 [scale=.15,auto=left, node distance=1.5cm, 
 ]
 \node[fill=white,style={circle,draw}] (n0) at (9,-10) {\small{1}};
 \node[fill=white,style={circle,draw}] (n1) at (4,0) {\small{2}};
  \node[fill=black!30,style={circle,draw}] (n2) at (14,0) {\small{3}};
\path
 (n1) [->,thick] edge[loop left=90,thick] node {} (n1);
 \path
        (n2)  [->,dashed]  edge[loop,thick] node {} (n2)
        (n2)  [->,dashed]  edge[loop right=90,thick] node {} (n2);
         \path (n1) [->,thick] edge[] node {{\tiny $c$}} (n0);
 \path (n2) [->,dashed] edge[above] node {{\tiny $d$}} (n0);
  \path (n0) [->,thick] edge[loop left=90] node {{\tiny $e$}} (n0);   
  \path  (n0) [->,thick] edge[bend left=40,below] node {{\tiny $a$}} (n1);   
 \end{tikzpicture} \\
 & & &  \\
 \hline 
\end{tabular}}
\end{center}
\caption{The connected $3$-node REI networks with valence $\leq 2$. 
Here $c,d,e$ are nonnegative integers such that $c + d +e\in \{0,1,2\}$.  Also, $a \in \{0,1,2\}$ (respectively  $b \in \{0,1,2\}$) is such that the sum of $a$ (respectively $b$) and the valence of node $2$ (respectively node $3$) is $\leq 2$.
}
\label{fig:3NCNREIV2}
\end{figure}

\begin{prop}
The set of connected $3$-node REI networks with valence $\leq 2$ 
comprises the networks in {\rm Figure}~{\rm \ref{fig:3NCNREIV2}}.
\label{prop:REI3}
\end{prop}

\begin{proof} 
We enumerate the set of connected $3$-node REI networks $\mathcal{G}$ with valence $\leq 2$ using Lemma~\ref{lemma:strategy_2_to_3}. We can assume that nodes $1$ and $2$ have type $N^E$ and node $3$ has type $N^I$.  

Consider the subnetwork $S$ of $\mathcal{G}$ containing node $2$ (of type $N^E$) and node $3$ (of type $N^I$) and all arrows between these nodes. This is a $2$-node REI network with valence $\leq 2$. 
If $S$ is connected then it is one of the $15$ networks in Figure~\ref{fig:2NCNREIV2} (Figure 7 in \cite{ADS24}), where node $2$ is of type $N^E$ and node $3$ is of type $N^I$. 
If $S$ is not connected then $S$ is one of the $9$ networks in Figure~\ref{fig:2N_disc_REIV2}.
The options for  arrows from $S$ to node $1$ and from node $1$ to $S$
 are shown in Figure \ref{fig:options}. 
 
 Since node $1$ has valence $\leq 2$, multiplicities $c,d,e$ satisfy $c + d +e\in \{0,1,2\}$.  Also, $a \in \{0,1,2\}$ (respectively  $b \in \{0,1,2\}$) is such that the sum of $a$ 
 (respectively  $b$) and the valence of node $2$ (respectively  node $3$) in $S$ is up to two. Combining this information with the networks in Figures~\ref{fig:2NCNREIV2} 
 and \ref{fig:2N_disc_REIV2} we obtain Figure~\ref{fig:3NCNREIV2}.
\end{proof}

\begin{figure}[!h]
{\tiny 
\begin{center}
\begin{tabular}{|l|l|l|l|} 
 \hline 
  & & &  \\
$(a)$ \begin{tikzpicture}
 [scale=.15,auto=left, node distance=1.5cm, 
 ]
 \node[fill=white,style={circle,draw}] (n1) at (4,0) {\small{1}};
  \node[fill=black!30,style={circle,draw}] (n2) at (14,0) {\small{2}};
 \draw[->, thick] (n1) edge[thick] node {}  (n2); 
 \end{tikzpicture}  
 &  
$(b)$ \begin{tikzpicture}
 [scale=.15,auto=left, node distance=1.5cm, 
 ]
 \node[fill=white,style={circle,draw}] (n1) at (4,0) {\small{1}};
  \node[fill=black!30,style={circle,draw}] (n2) at (14,0) {\small{2}};
 \draw[->, thick] (n1) edge[bend right=10, thick] node {}  (n2); 
 \draw[->, thick] (n1) edge[bend right=-10,thick] node {}  (n2); 
 \end{tikzpicture}  
  &  
$(c)$ \begin{tikzpicture}
 [scale=.15,auto=left, node distance=1.5cm, 
 ]
 \node[fill=white,style={circle,draw}] (n1) at (4,0) {\small{1}};
  \node[fill=black!30,style={circle,draw}] (n2) at (14,0) {\small{2}};
\path
        (n1) [->,solid, thick]  edge[thick] node { } (n2);
\path        
        (n2) [->,dashed, thick]  edge[loop right=90,thick] node { } (n2);
 \end{tikzpicture} 
  &  
$(d)$ \begin{tikzpicture}
 [scale=.15,auto=left, node distance=1.5cm, 
 ]
 \node[fill=white,style={circle,draw}] (n1) at (4,0) {\small{1}};
  \node[fill=black!30,style={circle,draw}] (n2) at (14,0) {\small{2}};
\path
        (n1) [->,thick]  edge[thick] node { } (n2)
        (n1) [->,thick]  edge[loop left=90,thick] node { } (n1);
 \end{tikzpicture} \\
 & & &  \\
 \hline 
 \hline 
  & & &  \\
 $(e)$ \begin{tikzpicture}
 [scale=.15,auto=left, node distance=1.5cm, 
 ]
 \node[fill=white,style={circle,draw}] (n1) at (4,0) {\small{1}};
  \node[fill=black!30,style={circle,draw}] (n2) at (14,0) {\small{2}};
\path
        (n1)  [->]  edge[loop,thick] node {} (n1)
        (n1)  [->]  edge[loop left=90,thick] node {} (n1)
        (n1) [->,thick]  edge[thick] node { } (n2);
 \end{tikzpicture}
  &  
 $(f)$  \begin{tikzpicture}
 [scale=.15,auto=left, node distance=1.5cm, 
 ]
 \node[fill=white,style={circle,draw}] (n1) at (4,0) {\small{1}};
  \node[fill=black!30,style={circle,draw}] (n2) at (14,0) {\small{2}};
\path
        (n1)  [->]  edge[loop left=90,thick] node {} (n1)
        (n1) [->,thick]  edge[bend right=10, thick] node { } (n2)
        (n1) [->,thick]  edge[bend right=-10, thick] node { } (n2);
 \end{tikzpicture}
&
$(g)$ \begin{tikzpicture}
 [scale=.15,auto=left, node distance=1.5cm, 
 ]
 \node[fill=white,style={circle,draw}] (n1) at (4,0) {\small{1}};
  \node[fill=black!30,style={circle,draw}] (n2) at (14,0) {\small{2}};
\path
        (n1) [->,thick]  edge[thick] node { } (n2)
         (n1) [->,thick]  edge[loop left=90, thick] node { } (n1);
\path 
        (n2) [->,dashed, thick]  edge[loop right=90, thick] node { } (n2);
 \end{tikzpicture} 
& 
 $(h)$ 
  \begin{tikzpicture}
 [scale=.15,auto=left, node distance=1.5cm, 
 ]
 \node[fill=white,style={circle,draw}] (n1) at (4,0) {\small{1}};
  \node[fill=black!30,style={circle,draw}] (n2) at (14,0) {\small{2}};
\path
        (n1)  [->]  edge[loop,thick] node {} (n1)
        (n1)  [->]  edge[loop left=90,thick] node {} (n1)
        (n1) [->,thick]  edge[bend right=20, thick] node { } (n2)
        (n1) [->,thick]  edge[bend right=-30, thick] node { } (n2);
 \end{tikzpicture}\\
 & & &  \\
 \hline 
 \hline 
  & & &  \\
 $(i)$  \begin{tikzpicture}
 [scale=.15,auto=left, node distance=1.5cm, 
 ]
 \node[fill=white,style={circle,draw}] (n1) at (4,0) {\small{1}};
  \node[fill=black!30,style={circle,draw}] (n2) at (14,0) {\small{2}};
\path
        (n1)  [->]  edge[loop,thick] node {} (n1)
        (n1)  [->]  edge[loop left=90,thick] node {} (n1)
        (n1) [->,thick]  edge[thick] node { } (n2);
\path 
         (n2)  [->, dashed]  edge[loop right=90,thick] node {} (n2); 
 \end{tikzpicture}
&
$(j)$ \begin{tikzpicture}
 [scale=.15,auto=left, node distance=1.5cm, 
 ]
 \node[fill=white,style={circle,draw}] (n1) at (4,0) {\small{1}};
  \node[fill=black!30,style={circle,draw}] (n2) at (14,0) {\small{2}};
\path
        (n1) [->,thick]  edge[bend right=10] node { } (n2);
 \path 
        (n2) [->,dashed, thick]  edge[bend left=-10] node { } (n1);
 \end{tikzpicture}
&
 $(k)$ \begin{tikzpicture}
 [scale=.15,auto=left, node distance=1.5cm, 
 ]
 \node[fill=white,style={circle,draw}] (n1) at (4,0) {\small{1}};
  \node[fill=black!30,style={circle,draw}] (n2) at (14,0) {\small{2}};
\path
        (n1)  [->]  edge[loop left=90,thick] node {} (n1)
        (n1) [->,thick]  edge[bend right=10, thick] node { } (n2);
\path 
        (n2) [->,dashed, thick]  edge[bend left=-10] node { } (n1);
 \end{tikzpicture}
  &   
 $(l)$ \begin{tikzpicture}
 [scale=.15,auto=left, node distance=1.5cm, 
 ]
 \node[fill=white,style={circle,draw}] (n1) at (4,0) {\small{1}};
  \node[fill=black!30,style={circle,draw}] (n2) at (14,0) {\small{2}};
\path
         (n1) [->,thick]  edge[bend right=10, thick] node { } (n2);
\path 
        (n2) [->,thick, dashed]  edge[bend left =-10] node { } (n1)
        (n2) [->,thick, dashed]  edge[bend left=-30] node { } (n1);
 \end{tikzpicture}\\
  & & &  \\
 \hline 
 \hline 
  & & &  \\
 $(m)$   \begin{tikzpicture}
 [scale=.15,auto=left, node distance=1.5cm, 
 ]
 \node[fill=white,style={circle,draw}] (n1) at (4,0) {\small{1}};
  \node[fill=black!30,style={circle,draw}] (n2) at (14,0) {\small{2}};
\path
        (n1)  [->]  edge[loop left=90,thick] node {} (n1)
        (n1) [->,thick]  edge[bend right=10, thick] node { } (n2)
         (n1) [->,thick]  edge[bend right=40, thick] node { } (n2);
\path 
        (n2) [->,thick, dashed]  edge[bend left=-10, thick] node { } (n1);      
 \end{tikzpicture}
 &  
 $(n)$  \begin{tikzpicture}
 [scale=.15,auto=left, node distance=1.5cm, 
 ]
 \node[fill=white,style={circle,draw}] (n1) at (4,0) {\small{1}};
  \node[fill=black!30,style={circle,draw}] (n2) at (14,0) {\small{2}};
\path
        (n1)  [->]  edge[loop left=90,thick] node {} (n1)
         (n1) [->,thick]  edge[bend right=10] node { } (n2);
\path          
         (n2)  [->, dashed]  edge[loop right=90,thick] node {} (n2)       
        (n2) [->,thick, dashed]  edge[bend left=-10] node { } (n1);
 \end{tikzpicture} 
 & 
 $(o)$ \begin{tikzpicture}
 [scale=.15,auto=left, node distance=1.5cm, 
 ]
 \node[fill=white,style={circle,draw}] (n1) at (4,0) {\small{1}};
  \node[fill=black!30,style={circle,draw}] (n2) at (14,0) {\small{2}};
\path
        (n1)  [->]  edge[bend right=10,thick] node {} (n2)
         (n1) [->,thick]  edge[bend right=40, thick] node { } (n2);
 \path 
         (n2)  [->,dashed]  edge[bend left=-10,thick] node {} (n1)      
        (n2) [->,thick,dashed]  edge[bend left=-40, thick] node { } (n1);
 \end{tikzpicture} 
&	\\
 & & &  \\
 \hline 
\end{tabular}
\end{center}
}
\caption{Connected 2-node REI networks with input valence $\leq 2$. This corresponds to \cite[Figure 7]{ADS24}.}
\label{fig:2NCNREIV2}
\end{figure}
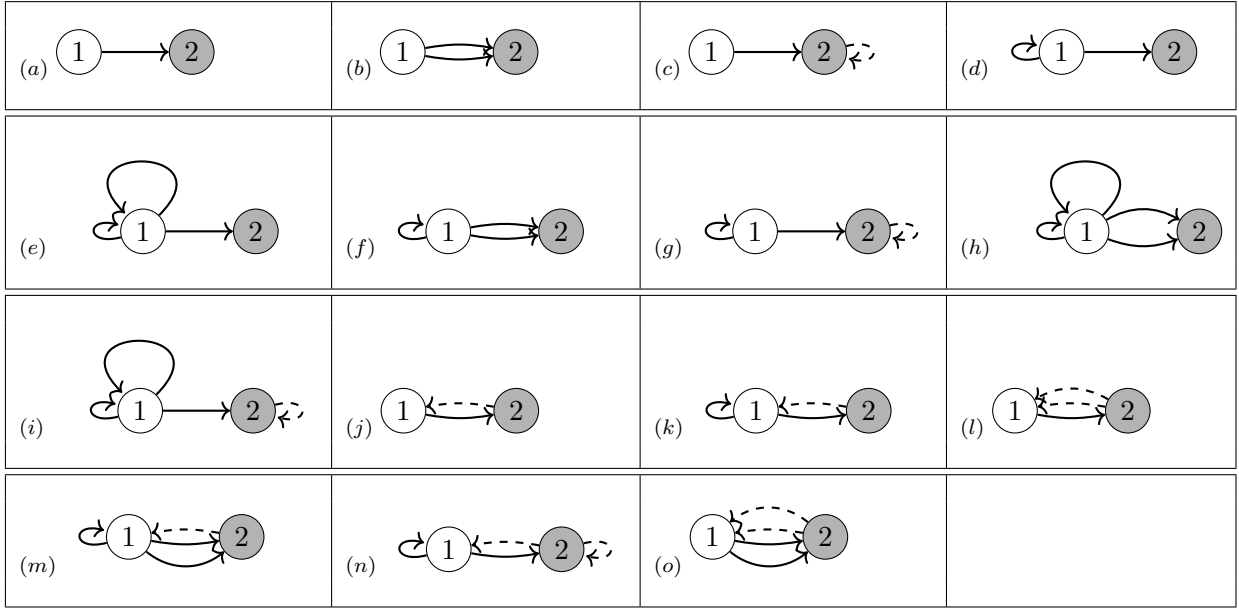

\begin{figure}
\begin{center}
{\tiny 
\begin{tabular}{|l|l|l|} 
 \hline 
  & &   \\
$D1$ \begin{tikzpicture}
 [scale=.15,auto=left, node distance=1.5cm, 
 ]
 \node[fill=white,style={circle,draw}] (n1) at (4,0) {\small{2}};
  \node[fill=black!30,style={circle,draw}] (n2) at (14,0) {\small{3}};
 \end{tikzpicture}  
\quad & \quad 
$D2$ \begin{tikzpicture}
 [scale=.15,auto=left, node distance=1.5cm, 
 ]
 \node[fill=white,style={circle,draw}] (n1) at (4,0) {\small{2}};
  \node[fill=black!30,style={circle,draw}] (n2) at (14,0) {\small{3}};
\path
 (n1) [->,thick] edge[loop left=90,thick] node {} (n1);
 \end{tikzpicture}   
 \quad & \quad 
$D3$ \begin{tikzpicture}
 [scale=.15,auto=left, node distance=1.5cm, 
 ]
 \node[fill=white,style={circle,draw}] (n1) at (4,0) {\small{2}};
  \node[fill=black!30,style={circle,draw}] (n2) at (14,0) {\small{3}};
  \path 
  (n2) [->,dashed] edge[loop right=90,thick] node {} (n2);
 \end{tikzpicture}  \\
 & &   \\
 \hline 
 \hline 
  & &   \\
$D4$ \begin{tikzpicture}
 [scale=.15,auto=left, node distance=1.5cm,
 ]
 \node[fill=white,style={circle,draw}] (n1) at (4,0) {\small{2}};
  \node[fill=black!30,style={circle,draw}] (n2) at (14,0) {\small{3}};
\path
 (n1) [->,thick] edge[loop left=90,thick] node {} (n1);
 \path 
(n2) [->,dashed] edge[loop right=90,thick] node {} (n2);
 \end{tikzpicture} 
  \quad & \quad 
  $D5$ \begin{tikzpicture}
 [scale=.15,auto=left, node distance=1.5cm, 
 ]
 \node[fill=white,style={circle,draw}] (n1) at (4,0) {\small{2}};
  \node[fill=black!30,style={circle,draw}] (n2) at (14,0) {\small{3}};
\path
        (n1)  [->]  edge[loop,thick] node {} (n1)
        (n1)  [->]  edge[loop left=90,thick] node {} (n1);
 \end{tikzpicture}   
 \quad & \quad 
 $D6$ \begin{tikzpicture}
 [scale=.15,auto=left, node distance=1.5cm, 
 ]
 \node[fill=white,style={circle,draw}] (n1) at (4,0) {\small{2}};
  \node[fill=black!30,style={circle,draw}] (n2) at (14,0) {\small{3}};
\path
        (n2)  [->,dashed]  edge[loop,thick] node {} (n2)
        (n2)  [->,dashed]  edge[loop right=90,thick] node {} (n2);
 \end{tikzpicture}   \\
  & &   \\
 \hline 
 \hline 
  & &   \\
 $D7$ \begin{tikzpicture}
 [scale=.15,auto=left, node distance=1.5cm, 
 ]
 \node[fill=white,style={circle,draw}] (n1) at (4,0) {\small{2}};
  \node[fill=black!30,style={circle,draw}] (n2) at (14,0) {\small{3}};
\path
  (n1)  [->]  edge[loop,thick] node {} (n1)
        (n1)  [->]  edge[loop left=90,thick] node {} (n1);
 \path 
        (n2)  [->,dashed]  edge[loop right=90,thick] node {} (n2);
 \end{tikzpicture} 
  \quad & \quad 
  $D8$ \begin{tikzpicture}
 [scale=.15,auto=left, node distance=1.5cm, 
 ]
 \node[fill=white,style={circle,draw}] (n1) at (4,0) {\small{2}};
  \node[fill=black!30,style={circle,draw}] (n2) at (14,0) {\small{3}};
\path 
        (n1)  [->]  edge[loop,thick] node {} (n1)
        (n1)  [->]  edge[loop left=90,thick] node {} (n1);
\path
        (n2)  [->,dashed]  edge[loop,thick] node {} (n2)
        (n2)  [->,dashed]  edge[loop right=90,thick] node {} (n2);
 \end{tikzpicture} 
  \quad & \quad 
  $D9$ \begin{tikzpicture}
 [scale=.15,auto=left, node distance=1.5cm, 
 ]
 \node[fill=white,style={circle,draw}] (n1) at (4,0) {\small{2}};
  \node[fill=black!30,style={circle,draw}] (n2) at (14,0) {\small{3}};
\path
 (n1) [->,thick] edge[loop left=90,thick] node {} (n1);
 \path
        (n2)  [->,dashed]  edge[loop,thick] node {} (n2)
        (n2)  [->,dashed]  edge[loop right=90,thick] node {} (n2);
 \end{tikzpicture} \\
 & &   \\
 \hline 
\end{tabular}}
\end{center}
\caption{The $2$-node disconnected REI networks with valence $\leq 2$.
}
\label{fig:2N_disc_REIV2}
\end{figure}
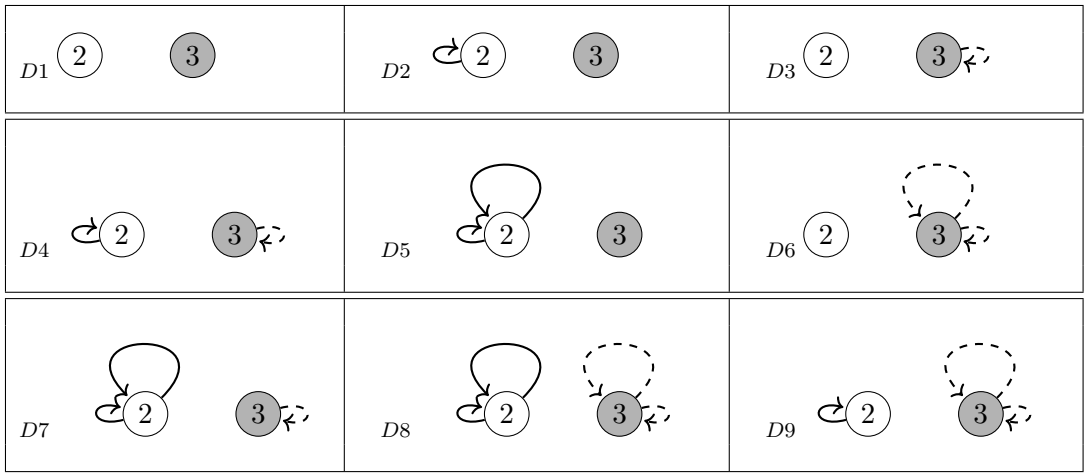

 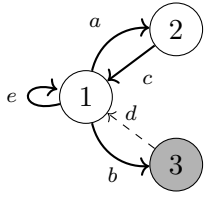
\begin{figure}
\begin{center}
 \begin{tikzpicture}
 [scale=.15,auto=left, node distance=1.5cm, 
 ]
 \node[fill=white,style={circle,draw}] (n1) at (4,6) {\small{2}};
 \node[fill=black!30,style={circle,draw}] (n3) at (4,-6) {\small{3}};
  \node[fill=white,style={circle,draw}] (n2) at (-4,0) {\small{1}};
  \path (n1) [->,thick] edge[] node {{\tiny $c$}} (n2);
   \path (n3) [->,dashed] edge[above] node {{\tiny $d$}} (n2);
    \path (n2) [->,thick] edge[loop left=90] node {{\tiny $e$}} (n2);
    \path (n2) [->,thick] edge[bend left=40] node {{\tiny $a$}} (n1);
    \path (n2) [->,thick] edge[bend right=40,below] node {{\tiny $b$}} (n3);
 \end{tikzpicture} 
 \end{center}
 \caption{Options for arrows from $S$ to node $1$ and from node $1$ to $S$.}
 \label{fig:options}
 \end{figure}
 
\subsection{Connected 3-node REI Networks with Valence 2}
\label{S:C3REIV2}

In this section we classify
connected $3$-node REI networks with valence $2$. 
We consider four different cases:
\begin{itemize}
\item[(i)] Every node receives one arrow of each type; 
\item[(ii)] Only the two excitatory nodes receive one arrow of each type; 
\item[(iii)] Only the inhibitory node and one excitatory node receive one arrow of each type; 
\item[(iv)] Given any two nodes there is no arrow-type preserving bijection between their input sets.
 \end{itemize}

We start 
by classifying the  $3$-node REI networks of valence $2$ that are  {\it almost homogeneous}; 
that is, where every node receives exactly one excitatory and one inhibitory arrow.
(The obstacle to exact homogeneity is that the nodes have different types.)

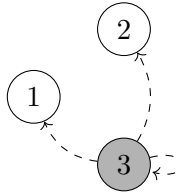
\begin{figure}[!h]
\begin{center}
 \begin{tikzpicture}
 [scale=.15,auto=left, node distance=1.5cm, 
 ]
 \node[fill=white,style={circle,draw}] (n1) at (4,6) {\small{2}};
 \node[fill=black!30,style={circle,draw}] (n3) at (4,-6) {\small{3}};
  \node[fill=white,style={circle,draw}] (n2) at (-4,0) {\small{1}};
   \path (n3) [->,dashed] edge[loop right=90] node {} (n3); 
     \path (n3) [->,dashed] edge[bend right=30] node {} (n1);
      \path (n3) [->,dashed] edge[bend left=30] node {} (n2);
 \end{tikzpicture} 
 \end{center}
 \caption{A $3$-node network where node $3$ outputs an inhibitory arrow to every node.}
 \label{f:arrow_I}
 \end{figure}

\begin{lemma}\label{lemma:strategy_sub_inhibition}
If $\mathcal{G}$ is 
an almost  
homogeneous connected $3$-node REI network of valence $2$, with $N^E = \{1,2\}$ and $N^I = \{ 3\}$ and arrow-types $A^E$ and $A^I$, then the subnetwork of $\mathcal{G}$ containing only arrows of type $A^I$ is the network in {\rm Figure}~{\rm \ref{f:arrow_I}}.
\end{lemma}

\begin{proof}
Since $\mathcal{G}$ is an almost homogeneous REI and node $3$ is the only one of type $N^I$, every node receives one arrow of type $A^I$ from node $3$.
\end{proof}

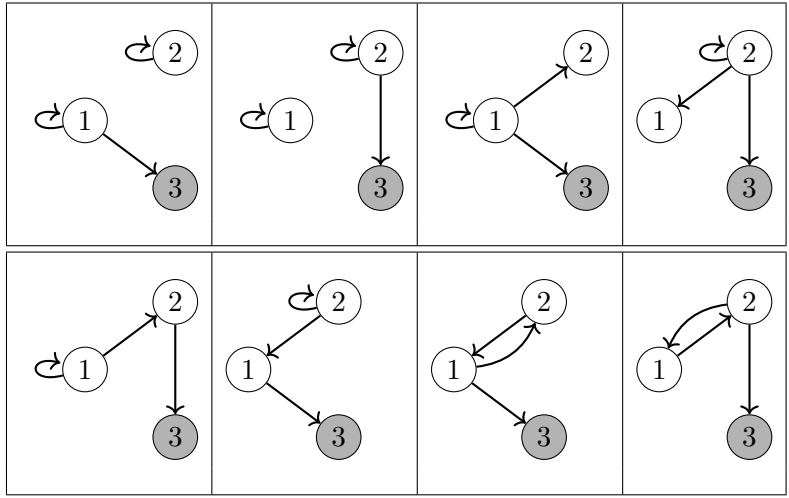
\begin{figure}
\begin{center}
{\tiny 
\begin{tabular}{|l|l|l|l|}
\hline
 & & & \\
 \begin{tikzpicture}
 [scale=.15,auto=left, node distance=1.5cm, 
 ]
 \node[fill=white,style={circle,draw}] (n1) at (4,6) {\small{2}};
 \node[fill=black!30,style={circle,draw}] (n3) at (4,-6) {\small{3}};
  \node[fill=white,style={circle,draw}] (n2) at (-4,0) {\small{1}};
   \path (n2) [->,thick] edge[] node {} (n3);
    \path (n2) [->,thick] edge[loop left=90] node {} (n2);
    \path (n1) [->,thick] edge[loop left=90] node {} (n1);
 \end{tikzpicture} &
 \begin{tikzpicture}
 [scale=.15,auto=left, node distance=1.5cm, 
 ]
 \node[fill=white,style={circle,draw}] (n1) at (4,6) {\small{2}};
 \node[fill=black!30,style={circle,draw}] (n3) at (4,-6) {\small{3}};
  \node[fill=white,style={circle,draw}] (n2) at (-4,0) {\small{1}};
   \path (n1) [->,thick] edge[] node {} (n3);
    \path (n2) [->,thick] edge[loop left=90] node {} (n2);
    \path (n1) [->,thick] edge[loop left=90] node {} (n1);
 \end{tikzpicture} &
 \begin{tikzpicture}
 [scale=.15,auto=left, node distance=1.5cm, 
 ]
 \node[fill=white,style={circle,draw}] (n1) at (4,6) {\small{2}};
 \node[fill=black!30,style={circle,draw}] (n3) at (4,-6) {\small{3}};
  \node[fill=white,style={circle,draw}] (n2) at (-4,0) {\small{1}};
   \path (n2) [->,thick] edge[] node {} (n3);
    \path (n2) [->,thick] edge[loop left=90] node {} (n2);
     \path (n2) [->,thick] edge[] node {} (n1);
 \end{tikzpicture} &
 \begin{tikzpicture}
 [scale=.15,auto=left, node distance=1.5cm, 
 ]
 \node[fill=white,style={circle,draw}] (n1) at (4,6) {\small{2}};
 \node[fill=black!30,style={circle,draw}] (n3) at (4,-6) {\small{3}};
  \node[fill=white,style={circle,draw}] (n2) at (-4,0) {\small{1}};
    \path (n1) [->,thick] edge[] node {} (n2);
    \path (n1) [->,thick] edge[loop left=90] node {} (n1);
     \path (n1) [->,thick] edge[] node {} (n3);
 \end{tikzpicture}   \\
 & & &  \\
 \hline 
 \hline 
  & & &  \\
  \begin{tikzpicture}
 [scale=.15,auto=left, node distance=1.5cm, 
 ]
 \node[fill=white,style={circle,draw}] (n1) at (4,6) {\small{2}};
 \node[fill=black!30,style={circle,draw}] (n3) at (4,-6) {\small{3}};
  \node[fill=white,style={circle,draw}] (n2) at (-4,0) {\small{1}};
   \path (n2) [->,thick] edge[] node {} (n1);
    \path (n2) [->,thick] edge[loop left=90] node {} (n2);
     \path (n1) [->,thick] edge[] node {} (n3);
 \end{tikzpicture} &
 \begin{tikzpicture}
 [scale=.15,auto=left, node distance=1.5cm, 
 ]
 \node[fill=white,style={circle,draw}] (n1) at (4,6) {\small{2}};
 \node[fill=black!30,style={circle,draw}] (n3) at (4,-6) {\small{3}};
  \node[fill=white,style={circle,draw}] (n2) at (-4,0) {\small{1}};
    \path (n1) [->,thick] edge[] node {} (n2);
    \path (n1) [->,thick] edge[loop left=90] node {} (n1);
     \path (n2) [->,thick] edge[] node {} (n3);
 \end{tikzpicture} &  
 \begin{tikzpicture}
 [scale=.15,auto=left, node distance=1.5cm, 
 ]
 \node[fill=white,style={circle,draw}] (n1) at (4,6) {\small{2}};
 \node[fill=black!30,style={circle,draw}] (n3) at (4,-6) {\small{3}};
  \node[fill=white,style={circle,draw}] (n2) at (-4,0) {\small{1}};
   \path (n2) [->,thick] edge[] node {} (n3);
     \path (n1) [->,thick] edge[] node {} (n2);
     \path (n2) [->,thick] edge[bend right=30] node {} (n1);
 \end{tikzpicture} &
 \begin{tikzpicture}
 [scale=.15,auto=left, node distance=1.5cm, 
 ]
 \node[fill=white,style={circle,draw}] (n1) at (4,6) {\small{2}};
 \node[fill=black!30,style={circle,draw}] (n3) at (4,-6) {\small{3}};
  \node[fill=white,style={circle,draw}] (n2) at (-4,0) {\small{1}};
    \path (n2) [->,thick] edge[] node {} (n1);
    \path (n1) [->,thick] edge[] node {} (n3);
       \path (n1) [->,thick] edge[bend right=30] node {} (n2);
 \end{tikzpicture} \\
  & & & \\
 \hline
 \end{tabular}}
 \end{center}
 \caption{The $3$-node networks in which every node receives an excitatory  arrow,
  which can be from node $1$ or node $2$.}
 \label{f:arrow_E}
 \end{figure}

\begin{lemma} \label{lemma:strategy_sub_activation}
If $\mathcal{G}$ is 
an almost  
homogeneous connected $3$-node REI network of valence $2$, with $N^E = \{1,2\}$ and $N^I = \{ 3\}$, and arrow-types $A^E$ and $A^I$, then 
the subnetwork of $\mathcal{G}$ containing only arrows of type $A^E$ is one of the networks in {\rm Figure~{\rm \ref{f:arrow_E}}}.
\end{lemma}

\begin{proof}
Since $\mathcal{G}$ is an almost 
homogeneous REI and node $3$ is the only of type $N^I$, every node receives
one arrow of type $A^E$ from nodes $1$ or $2$.
\end{proof}

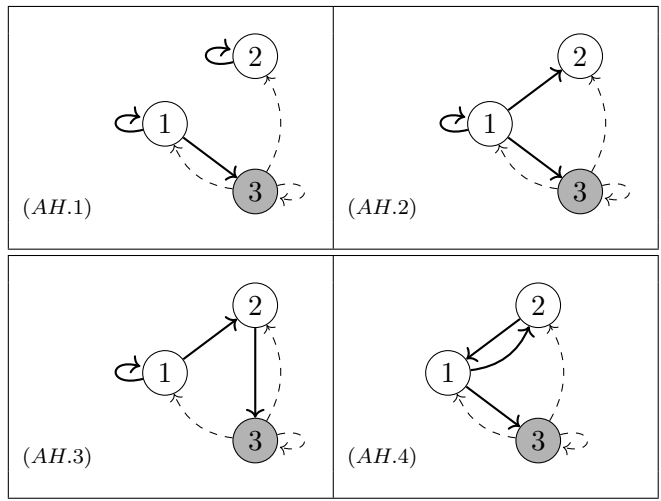
\begin{figure}
\begin{center}
{\tiny 
\begin{tabular}{|l|l|}
\hline
 & \\
 $(AH.1)$ \begin{tikzpicture}
 [scale=.15,auto=left, node distance=1.5cm, 
 ]
 \node[fill=white,style={circle,draw}] (n1) at (4,6) {\small{2}};
 \node[fill=black!30,style={circle,draw}] (n3) at (4,-6) {\small{3}};
  \node[fill=white,style={circle,draw}] (n2) at (-4,0) {\small{1}};
   \path (n2) [->,thick] edge[] node {} (n3);
    \path (n2) [->,thick] edge[loop left=90] node {} (n2);
    \path (n1) [->,thick] edge[loop left=90] node {} (n1);
     \path (n3) [->,dashed] edge[loop right=90] node {} (n3); 
     \path (n3) [->,dashed] edge[bend right=30] node {} (n1);
      \path (n3) [->,dashed] edge[bend left=30] node {} (n2);
 \end{tikzpicture}  &
$(AH.2)$ \begin{tikzpicture}
 [scale=.15,auto=left, node distance=1.5cm, 
 ]
 \node[fill=white,style={circle,draw}] (n1) at (4,6) {\small{2}};
 \node[fill=black!30,style={circle,draw}] (n3) at (4,-6) {\small{3}};
  \node[fill=white,style={circle,draw}] (n2) at (-4,0) {\small{1}};
   \path (n2) [->,thick] edge[] node {} (n3);
    \path (n2) [->,thick] edge[loop left=90] node {} (n2);
     \path (n2) [->,thick] edge[] node {} (n1);
     \path (n3) [->,dashed] edge[loop right=90] node {} (n3); 
     \path (n3) [->,dashed] edge[bend right=30] node {} (n1);
      \path (n3) [->,dashed] edge[bend left=30] node {} (n2);
 \end{tikzpicture}    \\
 &   \\
 \hline 
 \hline 
  &   \\
 $(AH.3)$ \begin{tikzpicture}
 [scale=.15,auto=left, node distance=1.5cm, 
 ]
 \node[fill=white,style={circle,draw}] (n1) at (4,6) {\small{2}};
 \node[fill=black!30,style={circle,draw}] (n3) at (4,-6) {\small{3}};
  \node[fill=white,style={circle,draw}] (n2) at (-4,0) {\small{1}};
   \path (n2) [->,thick] edge[] node {} (n1);
    \path (n2) [->,thick] edge[loop left=90] node {} (n2);
     \path (n1) [->,thick] edge[] node {} (n3);
     \path (n3) [->,dashed] edge[loop right=90] node {} (n3); 
     \path (n3) [->,dashed] edge[bend right=30] node {} (n1);
      \path (n3) [->,dashed] edge[bend left=30] node {} (n2);
 \end{tikzpicture}  &  
 $(AH.4)$ \begin{tikzpicture}
 [scale=.15,auto=left, node distance=1.5cm, 
 ]
 \node[fill=white,style={circle,draw}] (n1) at (4,6) {\small{2}};
 \node[fill=black!30,style={circle,draw}] (n3) at (4,-6) {\small{3}};
  \node[fill=white,style={circle,draw}] (n2) at (-4,0) {\small{1}};
   \path (n2) [->,thick] edge[] node {} (n3);
     \path (n1) [->,thick] edge[] node {} (n2);
     \path (n2) [->,thick] edge[bend right=30] node {} (n1);
     \path (n3) [->,dashed] edge[loop right=90] node {} (n3); 
     \path (n3) [->,dashed] edge[bend right=30] node {} (n1);
      \path (n3) [->,dashed] edge[bend left=30] node {} (n2);
 \end{tikzpicture}  \\
  &  \\
 \hline
 \end{tabular}}
 \end{center}
 \caption{
The almost homogeneous connected $3$-node REI networks with valence $2$, where nodes $1,2$ are of type $N^E$, node $3$ is of type $N^I$, and there are two arrow-types $A^E$ and $A^I$.
All networks have a unique $2$-dimensional robust synchrony space where only nodes $1,2$ are synchronized. }
 \label{fig:hom_3NCNREIV2}
 \end{figure}

\begin{prop} \label{prop:hom_REI3min}
Any almost 
homogeneous connected $3$-node REI network 
 of valence $2$ 
with $N^E=\{1,2\}$, $N^I = \{ 3\}$, and two arrow-types $A^E$ and $A^I$,
is one of the $4$ networks in {\rm Figure}~{\rm \ref{fig:hom_3NCNREIV2}}. These are not ODE-equivalent. 
Each of these networks has a unique $2$-dimensional 
robust synchrony subspace where only nodes $1,2$ are synchronized; see {\rm Remark {\rm \ref{R:2.3}(b)}} and {\rm Subsection~\ref{sub:synchrony}}.

\end{prop}

\begin{proof} We can assume that REI networks have nodes $1$ and $2$  of type $N^E$ and node $3$ of type $N^I$. If $\mathcal{G}$ is an almost  homogeneous connected $3$-node REI network with valence $2$ then the subnetwork containing only the arrow-type $A^I$ is 
the network in Figure~\ref{f:arrow_I}, see 
Lemma~\ref{lemma:strategy_sub_inhibition}, and the subnetwork of $\mathcal{G}$ containing only
arrow-type $A^E$ is one of the networks listed in 
Figure~\ref{f:arrow_E}, see
Lemma~\ref{lemma:strategy_sub_activation}. The subnetwork containing only arrow-type $A^I$ is symmetric under transposition 
of nodes $1$ and $2$. 
We obtain the networks in Figure~\ref{fig:hom_3NCNREIV2}.
\end{proof}

We consider now 3-node REI networks $\mathcal{G}$ of valence 2 which are inhomogeneous,
where nodes $1,2$ are input equivalent, 
each receives one arrow of each type, but node $3$ does not receive 
one arrow of each type.

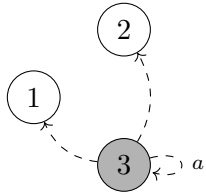
\begin{figure}[!h]
\begin{center}
 \begin{tikzpicture}
 [scale=.15,auto=left, node distance=1.5cm, 
 ]
 \node[fill=white,style={circle,draw}] (n1) at (4,6) {\small{2}};
 \node[fill=black!30,style={circle,draw}] (n3) at (4,-6) {\small{3}};
  \node[fill=white,style={circle,draw}] (n2) at (-4,0) {\small{1}};
   \path (n3) [->,dashed] edge[bend left=30] node {} (n2);
    \path (n3) [->,dashed] edge[loop right=90] node {{\tiny $a$}} (n3);
    \path (n3) [->,dashed] edge[bend right=30] node {} (n1);
 \end{tikzpicture} 
 \end{center}
 \caption{A $3$-node network in which node $3$ sends an 
 inhibitory arrow to nodes $1,2$. Here $a \in \{0,1,2\}$ is the number of inhibitory self-inputs of node $3$.}
 \label{f:12_arrow_I}
 \end{figure}

\begin{lemma}\label{lemma:12_strategy_sub_inhibition}
Let $\mathcal{G}$ be a connected $3$-node REI network of valence $2$, with $N^E = \{1,2\}$,  $N^I = \{ 3\}$, and arrow-types $A^E$ and $A^I$. Assume that 
nodes $1$ and $2$ are input equivalent, receiving one arrow of each type. 
Then the subnetwork of $\mathcal{G}$ containing only the arrow-type $A^I$ is the network in {\rm Figure}~{\rm \ref{f:12_arrow_I}}.
\end{lemma}
\begin{proof}
Since nodes $1,2$ of $\mathcal{G}$ are input equivalent and $\mathcal{G}$ is REI, node $3$ is the only one of type $N^I$, and nodes $1,2$ receive one arrow of type $A^I$ from node $3$.
\end{proof}

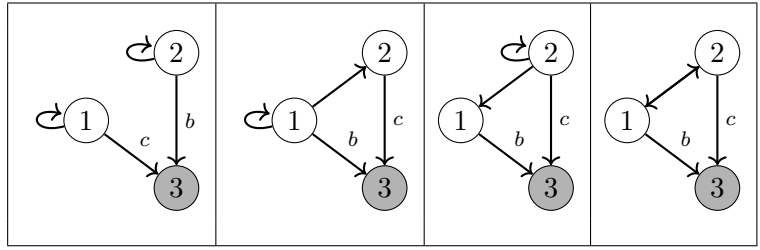
\begin{figure}
\begin{center}
{\tiny 
\begin{tabular}{|l|l|l|l|l|}
\hline
 & & & \\
 \begin{tikzpicture}
 [scale=.15,auto=left, node distance=1.5cm, 
 ]
 \node[fill=white,style={circle,draw}] (n1) at (4,6) {\small{2}};
 \node[fill=black!30,style={circle,draw}] (n3) at (4,-6) {\small{3}};
  \node[fill=white,style={circle,draw}] (n2) at (-4,0) {\small{1}};
   \path (n1) [->,thick] edge[] node {{\tiny $b$}} (n3);
    \path (n2) [->,thick] edge[] node {{\tiny $c$}} (n3);
    \path (n2) [->,thick] edge[loop left=90] node {} (n2);
    \path (n1) [->,thick] edge[loop left=90] node {} (n1);
 \end{tikzpicture} &
\begin{tikzpicture}
 [scale=.15,auto=left, node distance=1.5cm, 
 ]
 \node[fill=white,style={circle,draw}] (n1) at (4,6) {\small{2}};
 \node[fill=black!30,style={circle,draw}] (n3) at (4,-6) {\small{3}};
  \node[fill=white,style={circle,draw}] (n2) at (-4,0) {\small{1}};
   \path (n2) [->,thick] edge[] node {{\tiny $b$}} (n3);
    \path (n1) [->,thick] edge[] node {{\tiny $c$}} (n3);
    \path (n2) [->,thick] edge[loop left=90] node {} (n2);
      \path (n2) [->,thick] edge[] node {} (n1);
 \end{tikzpicture} &
 \begin{tikzpicture}
 [scale=.15,auto=left, node distance=1.5cm, 
 ]
 \node[fill=white,style={circle,draw}] (n1) at (4,6) {\small{2}};
 \node[fill=black!30,style={circle,draw}] (n3) at (4,-6) {\small{3}};
  \node[fill=white,style={circle,draw}] (n2) at (-4,0) {\small{1}};
   \path (n1) [->,thick] edge[] node {{\tiny $c$}} (n3);
    \path (n2) [->,thick] edge[] node {{\tiny $b$}} (n3);
    \path (n1) [->,thick] edge[] node {} (n2);
    \path (n1) [->,thick] edge[loop left=90] node {} (n1);
 \end{tikzpicture} &
 \begin{tikzpicture}
 [scale=.15,auto=left, node distance=1.5cm, 
 ]
 \node[fill=white,style={circle,draw}] (n1) at (4,6) {\small{2}};
 \node[fill=black!30,style={circle,draw}] (n3) at (4,-6) {\small{3}};
  \node[fill=white,style={circle,draw}] (n2) at (-4,0) {\small{1}};
   \path (n1) [->,thick] edge[] node {{\tiny $c$}} (n3);
    \path (n2) [->,thick] edge[] node {{\tiny $b$}} (n3);
    \path (n1) [->,thick] edge[] node {} (n2);
 \path (n2) [->,thick] edge[] node {} (n1);
 \end{tikzpicture}  \\
  & & &  \\
 \hline
 \end{tabular}}
 \end{center}
 \caption{The $3$-node networks  
where nodes $1, 2$ are excitatory and node $3$ is inhibitory, and
where nodes $1,2$ receive an excitatory arrow which can be from node $1$ or node $2$.  Here $b,c$ are nonnegative integers such that $b+c \in \{0,1,2\}$, representing the total number of 
excitatory inputs that node $3$ receives (from nodes $1,2$).}
 \label{f:12_arrow_E}
 \end{figure}

\begin{lemma} \label{lemma:12_strategy_sub_activation}
Let $\mathcal{G}$ be a connected $3$-node REI network of valence $2$ with $N^E = \{1,2\}$, $N^I = \{ 3\}$, and arrow-types $A^E$ and $A^I$. Assume that 
nodes $1$ and $2$ are input equivalent, each receiving one arrow of each type. 
Then the subnetwork of $\mathcal{G}$ containing only the arrow-type $A^E$ is one of the networks in {\rm Figure}~{\rm \ref{f:12_arrow_E}}. 
\end{lemma}
\begin{proof}
Since $\mathcal{G}$ is REI and nodes  $1$ and $2$ are of type $N^E$, each of nodes $1,2$ receives one arrow of type $A^E$ from nodes $1$ or $2$.
\end{proof}

\begin{figure}
\begin{center}
{\tiny 
\begin{tabular}{|l|l|l|}
\hline
 & &  \\
(NH.1) \begin{tikzpicture}
 [scale=.15,auto=left, node distance=1.5cm, 
 ]
 \node[fill=white,style={circle,draw}] (n1) at (4,6) {\small{2}};
 \node[fill=black!30,style={circle,draw}] (n3) at (4,-6) {\small{3}};
  \node[fill=white,style={circle,draw}] (n2) at (-4,0) {\small{1}};
   \path (n1) [->,thick] edge[] node [above left] {{\tiny $b$}} (n3);
    \path (n2) [->,thick] edge[] node {{\tiny $c$}} (n3);
    \path (n2) [->,thick] edge[loop left=90] node {} (n2);
    \path (n1) [->,thick] edge[loop left=90] node {} (n1);
    \path (n3) [->,dashed] edge[bend left=30] node {} (n2);
    \path (n3) [->,dashed] edge[loop right=90] node {{\tiny $a$}} (n3);
    \path (n3) [->,dashed] edge[bend right=30] node {} (n1);
 \end{tikzpicture} &
(NH.2) \begin{tikzpicture}
 [scale=.15,auto=left, node distance=1.5cm, 
 ]
 \node[fill=white,style={circle,draw}] (n1) at (4,6) {\small{2}};
 \node[fill=black!30,style={circle,draw}] (n3) at (4,-6) {\small{3}};
  \node[fill=white,style={circle,draw}] (n2) at (-4,0) {\small{1}};
   \path (n2) [->,thick] edge[] node {{\tiny $b$}} (n3);
    \path (n1) [->,thick] edge[] node [above left] {{\tiny $c$}} (n3);
    \path (n2) [->,thick] edge[loop left=90] node {} (n2);
      \path (n2) [->,thick] edge[] node {} (n1);
      \path (n3) [->,dashed] edge[bend left=30] node {} (n2);
    \path (n3) [->,dashed] edge[loop right=90] node {{\tiny $a$}} (n3);
    \path (n3) [->,dashed] edge[bend right=30] node {} (n1);
 \end{tikzpicture}  &
 (NH.3) \begin{tikzpicture}
 [scale=.15,auto=left, node distance=1.5cm, 
 ]
 \node[fill=white,style={circle,draw}] (n1) at (4,6) {\small{2}};
 \node[fill=black!30,style={circle,draw}] (n3) at (4,-6) {\small{3}};
  \node[fill=white,style={circle,draw}] (n2) at (-4,0) {\small{1}};
   \path (n1) [->,thick] edge[] node  [above left] {{\tiny $c$}} (n3);
    \path (n2) [->,thick] edge[] node {{\tiny $b$}} (n3);
    \path (n1) [->,thick] edge[] node {} (n2);
 \path (n2) [->,thick] edge[] node {} (n1);
 \path (n3) [->,dashed] edge[bend left=30] node {} (n2);
    \path (n3) [->,dashed] edge[loop right=90] node {{\tiny $a$}} (n3);
    \path (n3) [->,dashed] edge[bend right=30] node {} (n1);
 \end{tikzpicture}  \\
  & &   \\
 \hline
 \end{tabular}}
 \end{center}
 \caption{The connected $3$-node REI networks with valence $2$, 
two arrow-types $A^E$ and $A^I$, and where nodes $1,2$ are input equivalent receiving one input of each arrow-type. Here $a,b,c$ are nonnegative integers such that $a+b+c =2$ and $a \not=1$. That is, $a=2,\, b=c=0$ or $a=0,\, b+c=2$.}
 \label{fig:12_3NCNREIV2}
 \end{figure}
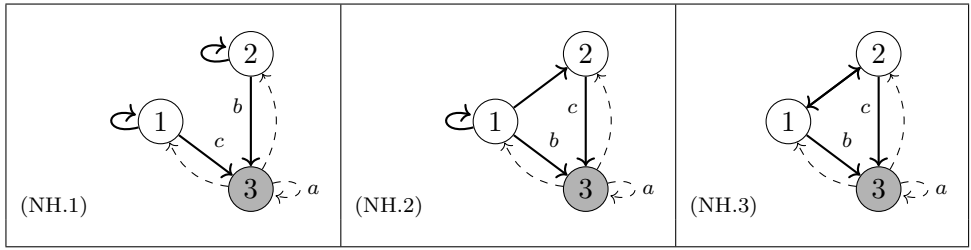

\begin{prop} \label{prop:hom_REI3min2}
Any connected $3$-node REI network of valence $2$, 
where the two excitatory nodes are input equivalent receiving one arrow of each type and the inhibitory node does not receive one arrow of each type,
is one of the networks in {\rm Figure}~{\rm \ref{fig:12_3NCNREIV2}}. All these networks have exactly
one $2$-dimensional robust
synchrony subspace where nodes $1$ and $2$ are synchronized.
\end{prop}

\begin{proof} We can assume that the REI network has nodes $1$ and $2$ of type $N^E$ and node $3$ of type $N^I$. Suppose that $\mathcal{G}$ is a minimal connected $3$-node REI network with valence $2$, input equivalence relation $\sim_I = \left\{ \{ 1,2\}, \{ 3\}\right\}$, and where nodes $1$ and $2$ receive one arrow of each type. Then the subnetwork containing only 
arrow-type $A^I$ is the network in Figure~\ref{f:12_arrow_I}, see Lemma~\ref{lemma:12_strategy_sub_inhibition}, and the subnetwork of $\mathcal{G}$ containing only arrow-type $A^E$ is one of the networks in Figure~\ref{f:12_arrow_E}, see Lemma~\ref{lemma:12_strategy_sub_activation}. The subnetwork containing only arrow-type $A^I$ is symmetric under transposition of nodes $1$ and $2$. We obtain the networks in Figure~\ref{fig:12_3NCNREIV2}. Clearly the only possible 
robust synchrony subspace must have  nodes 1 and 2 synchronized.
\end{proof}

We consider 3-node REI networks of valence 2, where we assume now that all three nodes are not input equivalent, but nodes 
$1,3$ receive one arrow of each type. 

\begin{figure}[!h]
\begin{center}
 \begin{tikzpicture}
 [scale=.15,auto=left, node distance=1.5cm, 
 ]
 \node[fill=white,style={circle,draw}] (n1) at (4,6) {\small{2}};
 \node[fill=black!30,style={circle,draw}] (n3) at (4,-6) {\small{3}};
  \node[fill=white,style={circle,draw}] (n2) at (-4,0) {\small{1}};
    \path (n3) [->,dashed] edge[bend left=30] node {} (n2);
    \path (n3) [->,dashed] edge[loop right=90] node {} (n3);
    \path (n3) [->,dashed] edge[bend right=70] node {{\tiny $a$}} (n1);
 \end{tikzpicture} 
 \end{center}
 \caption{$3$-node network in which 
nodes $1, 2$ are excitatory and node $3$ is inhibitory, and
 node $3$ sends an inhibitory arrow to nodes $1,3$. 
 Here, $a \in \{0,1,2\}$ is the number of 
  inhibitory inputs  to node $2$ (from node $3$).}
 \label{f:13_arrow_I}
 \end{figure}
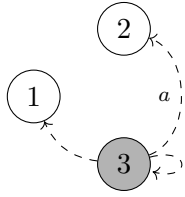

\begin{lemma}\label{lemma:13_strategy_sub_inhibition}
Let $\mathcal{G}$ be a connected $3$-node REI network of valence $2$ with $N^E = \{1,2\}$, $N^I = \{ 3\}$, and arrow-types $A^E$ and $A^I$. Assume that 
$\sim_I = \left\{ \{ 1\}, \{ 2\}, \{3\}\right\}$ 
 and that nodes $1$ and $3$ receive one arrow of each type.   
Then the subnetwork of $\mathcal{G}$ containing only arrow-type $A^I$ is the network in Figure~{\rm \ref{f:13_arrow_I}}.
\end{lemma}

\begin{proof}
Since nodes $1,3$ of $\mathcal{G}$ 
receive both an arrow of type $A^I$ 
 and $\mathcal{G}$ is REI, node $3$ is the only one of type $N^I$, and nodes $1,3$ receive one arrow of type $A^I$ from node $3$.
\end{proof}

\begin{figure}
\begin{center}
{\tiny 
\begin{tabular}{|c|c|c|c|}
\hline
 & & & \\
 \begin{tikzpicture}
 [scale=.15,auto=left, node distance=1.5cm, 
 ]
 \node[fill=white,style={circle,draw}] (n1) at (4,6) {\small{2}};
 \node[fill=black!30,style={circle,draw}] (n3) at (4,-6) {\small{3}};
  \node[fill=white,style={circle,draw}] (n2) at (-4,0) {\small{1}};
    \path (n2) [->,thick] edge[] node {} (n3);
    \path (n2) [->,thick] edge[loop  left=90] node {} (n2);
    \path (n1) [->,thick] edge[loop left=90,below] node {{\tiny $b$}} (n1);
 \end{tikzpicture} &
 \begin{tikzpicture}
 [scale=.15,auto=left, node distance=1.5cm, 
 ]
 \node[fill=white,style={circle,draw}] (n1) at (4,6) {\small{2}};
 \node[fill=black!30,style={circle,draw}] (n3) at (4,-6) {\small{3}};
  \node[fill=white,style={circle,draw}] (n2) at (-4,0) {\small{1}};
    \path (n1) [->,thick] edge[] node {} (n3);
    \path (n2) [->,thick] edge[loop  left=90] node {} (n2);
    \path (n1) [->,thick] edge[loop left=90,below] node {{\tiny $b$}} (n1);
 \end{tikzpicture} 
 &
\begin{tikzpicture}
 [scale=.15,auto=left, node distance=1.5cm, 
 ]
 \node[fill=white,style={circle,draw}] (n1) at (4,6) {\small{2}};
 \node[fill=black!30,style={circle,draw}] (n3) at (4,-6) {\small{3}};
  \node[fill=white,style={circle,draw}] (n2) at (-4,0) {\small{1}};
   \path (n2) [->,thick] edge[] node {{\tiny $b$}} (n1);
    \path (n2) [->,thick] edge[loop left=90] node {} (n2);
      \path (n2) [->,thick] edge[] node {} (n3);
 \end{tikzpicture} &
 \begin{tikzpicture}
 [scale=.15,auto=left, node distance=1.5cm, 
 ]
 \node[fill=white,style={circle,draw}] (n1) at (4,6) {\small{2}};
 \node[fill=black!30,style={circle,draw}] (n3) at (4,-6) {\small{3}};
  \node[fill=white,style={circle,draw}] (n2) at (-4,0) {\small{1}};
   \path (n2) [->,thick] edge[] node {{\tiny $b$}} (n1);
    \path (n2) [->,thick] edge[loop left=90] node {} (n2);
      \path (n1) [->,thick] edge[] node {} (n3);
 \end{tikzpicture} \\
 \hline
 \hline 
 \begin{tikzpicture}
 [scale=.15,auto=left, node distance=1.5cm, 
 ]
 \node[fill=white,style={circle,draw}] (n1) at (4,6) {\small{2}};
 \node[fill=black!30,style={circle,draw}] (n3) at (4,-6) {\small{3}};
  \node[fill=white,style={circle,draw}] (n2) at (-4,0) {\small{1}};
    \path (n2) [->,thick] edge[bend right=30,below] node {{\tiny $b$}} (n1);
    \path (n1) [->,thick] edge[] node {} (n2);
 \path (n2) [->,thick] edge[] node {} (n3);
 \end{tikzpicture}  &
 \begin{tikzpicture}
 [scale=.15,auto=left, node distance=1.5cm, 
 ]
 \node[fill=white,style={circle,draw}] (n1) at (4,6) {\small{2}};
 \node[fill=black!30,style={circle,draw}] (n3) at (4,-6) {\small{3}};
  \node[fill=white,style={circle,draw}] (n2) at (-4,0) {\small{1}};
    \path (n2) [->,thick] edge[bend right=30,below] node {{\tiny $b$}} (n1);
    \path (n1) [->,thick] edge[] node {} (n2);
 \path (n1) [->,thick] edge[] node {} (n3);
 \end{tikzpicture} &
 \begin{tikzpicture}
 [scale=.15,auto=left, node distance=1.5cm, 
 ]
 \node[fill=white,style={circle,draw}] (n1) at (4,6) {\small{2}};
 \node[fill=black!30,style={circle,draw}] (n3) at (4,-6) {\small{3}};
  \node[fill=white,style={circle,draw}] (n2) at (-4,0) {\small{1}};
  \path (n1) [->,thick] edge[] node {} (n2);
    \path (n2) [->,thick] edge[] node {} (n3);
    \path (n1) [->,thick] edge[loop left=90,below] node {{\tiny $b$}} (n1);
 \end{tikzpicture} &
 \begin{tikzpicture}
 [scale=.15,auto=left, node distance=1.5cm, 
 ]
 \node[fill=white,style={circle,draw}] (n1) at (4,6) {\small{2}};
 \node[fill=black!30,style={circle,draw}] (n3) at (4,-6) {\small{3}};
  \node[fill=white,style={circle,draw}] (n2) at (-4,0) {\small{1}};
    \path (n1) [->,thick] edge[] node {} (n2);
    \path (n1) [->,thick] edge[] node {} (n3);
    \path (n1) [->,thick] edge[loop left=90,below] node {{\tiny $b$}} (n1);
 \end{tikzpicture} \\
  & & &  \\
 \hline
 \end{tabular}}
 \end{center}
 \caption{The $3$-node networks in which
nodes $1, 2$ are excitatory, node $3$ is inhibitory, and
nodes $1,3$ receive an  excitatory arrow which can be from node $1$ or node $2$.  Here  $b \in \{0,1,2\}$ represents the total number of 
 excitatory inputs that node $2$ receives (from nodes $1,2$).}
 \label{f:13_arrow_E}
 \end{figure}
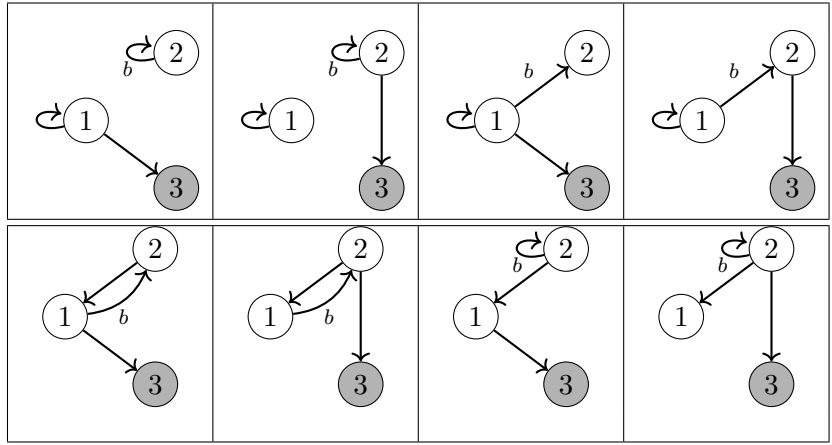

Recall from Definition \ref{Def:input_equiv} that $\sim_I$ denotes input equivalence.
We now prove:

\begin{lemma} \label{lemma:13_strategy_sub_activation}
Let $\mathcal{G}$ be a connected $3$-node REI network of valence 
$2$, with $N^E = \{1,2\}$, $N^I = \{ 3\}$, and arrow-types $A^E$ and $A^I$. 
Assume that 
$\sim_I = \left\{ \{ 1\}, \{ 2\}, \{3\}\right\}$  
and that nodes $1$ and $3$ receive one arrow of each type.   
Then the subnetwork of $\mathcal{G}$ containing only arrow-type $A^E$ is one of the networks in Figure~{\rm \ref{f:13_arrow_E}}. 
\end{lemma}

\begin{proof}
Since $\mathcal{G}$ is REI and nodes  $1$ and $2$ are those of type $N^E$,
 every node $1,3$ receives one arrow of type $A^E$ from node $1$ or $2$.
\end{proof}

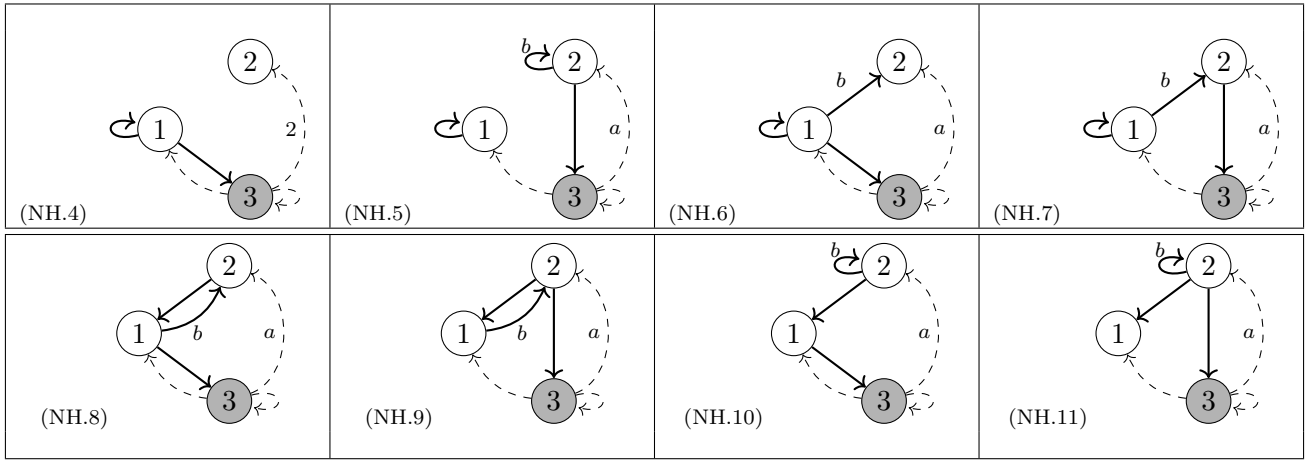
\begin{figure}
\begin{center}
{\tiny 
\begin{tabular}{|c|c|c|c|}
\hline 
 & & & \\
 (NH.4) \begin{tikzpicture}
 [scale=.15,auto=left, node distance=1.5cm, 
 ]
 \node[fill=white,style={circle,draw}] (n1) at (4,6) {\small{2}};
 \node[fill=black!30,style={circle,draw}] (n3) at (4,-6) {\small{3}};
  \node[fill=white,style={circle,draw}] (n2) at (-4,0) {\small{1}};
    \path (n2) [->,thick] edge[] node {} (n3);
    \path (n2) [->,thick] edge[loop  left=90] node {} (n2);
     \path (n3) [->,dashed] edge[bend left=30] node {} (n2);
    \path (n3) [->,dashed] edge[loop right=90] node {} (n3);
    \path (n3) [->,dashed] edge[bend right=70] node {{\tiny $2$}} (n1);
 \end{tikzpicture} &
(NH.5)  \begin{tikzpicture}
 [scale=.15,auto=left, node distance=1.5cm, 
 ]
 \node[fill=white,style={circle,draw}] (n1) at (4,6) {\small{2}};
 \node[fill=black!30,style={circle,draw}] (n3) at (4,-6) {\small{3}};
  \node[fill=white,style={circle,draw}] (n2) at (-4,0) {\small{1}};
    \path (n1) [->,thick] edge[] node {} (n3);
    \path (n2) [->,thick] edge[loop  left=90] node {} (n2);
    \path (n1) [->,thick] edge[loop left=90,above] node {{\tiny $b$}} (n1);
     \path (n3) [->,dashed] edge[bend left=30] node {} (n2);
    \path (n3) [->,dashed] edge[loop right=90] node {} (n3);
    \path (n3) [->,dashed] edge[bend right=70] node {{\tiny $a$}} (n1);
 \end{tikzpicture} 
 &
(NH.6) \begin{tikzpicture}
 [scale=.15,auto=left, node distance=1.5cm, 
 ]
 \node[fill=white,style={circle,draw}] (n1) at (4,6) {\small{2}};
 \node[fill=black!30,style={circle,draw}] (n3) at (4,-6) {\small{3}};
  \node[fill=white,style={circle,draw}] (n2) at (-4,0) {\small{1}};
   \path (n2) [->,thick] edge[] node {{\tiny $b$}} (n1);
    \path (n2) [->,thick] edge[loop left=90] node {} (n2);
      \path (n2) [->,thick] edge[] node {} (n3);
       \path (n3) [->,dashed] edge[bend left=30] node {} (n2);
    \path (n3) [->,dashed] edge[loop right=90] node {} (n3);
    \path (n3) [->,dashed] edge[bend right=70] node {{\tiny $a$}} (n1);
 \end{tikzpicture} &
 (NH.7) \begin{tikzpicture}
 [scale=.15,auto=left, node distance=1.5cm, 
 ]
 \node[fill=white,style={circle,draw}] (n1) at (4,6) {\small{2}};
 \node[fill=black!30,style={circle,draw}] (n3) at (4,-6) {\small{3}};
  \node[fill=white,style={circle,draw}] (n2) at (-4,0) {\small{1}};
   \path (n2) [->,thick] edge[] node {{\tiny $b$}} (n1);
    \path (n2) [->,thick] edge[loop left=90] node {} (n2);
      \path (n1) [->,thick] edge[] node {} (n3);
       \path (n3) [->,dashed] edge[bend left=30] node {} (n2);
    \path (n3) [->,dashed] edge[loop right=90] node {} (n3);
    \path (n3) [->,dashed] edge[bend right=70] node {{\tiny $a$}} (n1);
 \end{tikzpicture} \\
 \hline
 \hline 
(NH.8) \begin{tikzpicture}
 [scale=.15,auto=left, node distance=1.5cm, 
 ]
 \node[fill=white,style={circle,draw}] (n1) at (4,6) {\small{2}};
 \node[fill=black!30,style={circle,draw}] (n3) at (4,-6) {\small{3}};
  \node[fill=white,style={circle,draw}] (n2) at (-4,0) {\small{1}};
    \path (n2) [->,thick] edge[bend right=30,below] node {{\tiny $b$}} (n1);
    \path (n1) [->,thick] edge[] node {} (n2);
 \path (n2) [->,thick] edge[] node {} (n3);
  \path (n3) [->,dashed] edge[bend left=30] node {} (n2);
    \path (n3) [->,dashed] edge[loop right=90] node {} (n3);
    \path (n3) [->,dashed] edge[bend right=70] node {{\tiny $a$}} (n1);
 \end{tikzpicture}  &
 (NH.9) \begin{tikzpicture}
 [scale=.15,auto=left, node distance=1.5cm, 
 ]
 \node[fill=white,style={circle,draw}] (n1) at (4,6) {\small{2}};
 \node[fill=black!30,style={circle,draw}] (n3) at (4,-6) {\small{3}};
  \node[fill=white,style={circle,draw}] (n2) at (-4,0) {\small{1}};
    \path (n2) [->,thick] edge[bend right=30,below] node {{\tiny $b$}} (n1);
    \path (n1) [->,thick] edge[] node {} (n2);
 \path (n1) [->,thick] edge[] node {} (n3);
  \path (n3) [->,dashed] edge[bend left=30] node {} (n2);
    \path (n3) [->,dashed] edge[loop right=90] node {} (n3);
    \path (n3) [->,dashed] edge[bend right=70] node {{\tiny $a$}} (n1);
 \end{tikzpicture} &
(NH.10)  \begin{tikzpicture}
 [scale=.15,auto=left, node distance=1.5cm, 
 ]
 \node[fill=white,style={circle,draw}] (n1) at (4,6) {\small{2}};
 \node[fill=black!30,style={circle,draw}] (n3) at (4,-6) {\small{3}};
  \node[fill=white,style={circle,draw}] (n2) at (-4,0) {\small{1}};
  \path (n1) [->,thick] edge[] node {} (n2);
    \path (n2) [->,thick] edge[] node {} (n3);
    \path (n1) [->,thick] edge[loop left=90,above] node {{\tiny $b$}} (n1);
     \path (n3) [->,dashed] edge[bend left=30] node {} (n2);
    \path (n3) [->,dashed] edge[loop right=90] node {} (n3);
    \path (n3) [->,dashed] edge[bend right=70] node {{\tiny $a$}} (n1);
 \end{tikzpicture} &
 (NH.11) \begin{tikzpicture}
 [scale=.15,auto=left, node distance=1.5cm, 
 ]
 \node[fill=white,style={circle,draw}] (n1) at (4,6) {\small{2}};
 \node[fill=black!30,style={circle,draw}] (n3) at (4,-6) {\small{3}};
  \node[fill=white,style={circle,draw}] (n2) at (-4,0) {\small{1}};
    \path (n1) [->,thick] edge[] node {} (n2);
    \path (n1) [->,thick] edge[] node {} (n3);
    \path (n1) [->,thick] edge[loop left=90,above] node {{\tiny $b$}} (n1);
     \path (n3) [->,dashed] edge[bend left=30] node {} (n2);
    \path (n3) [->,dashed] edge[loop right=90] node {} (n3);
    \path (n3) [->,dashed] edge[bend right=70] node {{\tiny $a$}} (n1);
 \end{tikzpicture} \\
  & & & \\ 
 \hline
 \end{tabular}}
 \end{center}
 \caption{The connected $3$-node REI networks with valence $2$, 
where nodes $1, 2$ are excitatory and node $3$ is inhibitory, there are
two arrow-types $A^E$ and $A^I$, 
all nodes are not input equivalent, and nodes $1,3$ 
 receive one input of each arrow-type. For networks $(NH.5) - (NH.11)$, $a,b$ are nonnegative integers such that $a+b =2$ and $a \not=b$. That is, $a=0,\, b=2$ or $a=2,\, b=0$.}
 \label{fig:13_3NCNREIV2}
 \end{figure}

\begin{prop} \label{prop:hom_REI3}
Any connected $3$-node REI network of valence $2$, 
where two nodes are excitatory, one node is inhibitory, and all three nodes are not input equivalent but where one excitatory node and the inhibitory node receive one arrow of each type, 
is one of the networks listed in Figure~{\rm \ref{fig:13_3NCNREIV2}}. 
\end{prop}

\begin{proof} We can assume that nodes $1$ and $2$  have type $N^E$ and node $3$ has type $N^I$. 
Let $\mathcal{G}$ be a connected $3$-node REI network with valence $2$, 
and input equivalence relation 
$\sim_I = \left\{ \{1\}, \{3\}, \{2\}\right\}$,
where  nodes $1$ and $3$ receive one arrow of each type. Then 
the subnetwork containing only arrow-type $A^I$ is  
the network in Figure~\ref{f:13_arrow_I}, see 
Lemma~\ref{lemma:13_strategy_sub_inhibition} and the subnetwork of $\mathcal{G}$ containing only 
arrow-type $A^E$ is one of the networks listed in 
Figure~\ref{f:13_arrow_E}, see 
Lemma~\ref{lemma:13_strategy_sub_activation}. 
We obtain the networks in Figure~\ref{fig:13_3NCNREIV2}.
\end{proof}

We consider 3-node REI networks of valence 2, where
we now assume that, given any two nodes, there is no arrow-type preserving bijection between their input sets: 
\begin{equation}
\label{E:C}
\begin{array}{l}
\mbox{One node receives one arrow of each type}\ A^E\  \mbox{and}\ A^I,\\
\mbox{another node receives two arrows of type}\ A^E\\
\mbox{and the other node receives two arrows of type}\ A^I.
\end{array}
\end{equation}

Thus, each node lies in a different input equivalence class, 
that is, $\sim_I = \left\{ \{ 1\}, \{ 2\}, \{ 3\} \right\}$.

\begin{figure}
\begin{center}
{\tiny 
\begin{tabular}{|c|c|c|}
\hline 
 & &  \\
 $(a)$ 
\begin{tikzpicture}
 [scale=.15,auto=left, node distance=1.5cm, 
 ]
 \node[fill=white,style={circle,draw}] (n2) at (4,6) {\small{2}};
 \node[fill=black!30,style={circle,draw}] (n3) at (4,-6) {\small{3}};
  \node[fill=white,style={circle,draw}] (n1) at (-4,0) {\small{1}};
    \path (n3) [->,dashed] edge[] node {} (n1);
     \path (n3) [->,dashed] edge[bend left=10] node {} (n2);
     \path (n3) [->,dashed] edge[bend right=10] node {} (n2);
 \end{tikzpicture} &
$(b)$  
\begin{tikzpicture}
 [scale=.15,auto=left, node distance=1.5cm, 
 ]
 \node[fill=white,style={circle,draw}] (n2) at (4,6) {\small{2}};
 \node[fill=black!30,style={circle,draw}] (n3) at (4,-6) {\small{3}};
  \node[fill=white,style={circle,draw}] (n1) at (-4,0) {\small{1}};
    \path (n3) [->,dashed] edge[] node {} (n2);
     \path (n3) [->,dashed] edge[bend left=10] node {} (n1);
     \path (n3) [->,dashed] edge[bend right=10] node {} (n1);
 \end{tikzpicture}
 &
$(c)$
\begin{tikzpicture}
 [scale=.15,auto=left, node distance=1.5cm, 
 ]
 \node[fill=white,style={circle,draw}] (n2) at (4,6) {\small{2}};
 \node[fill=black!30,style={circle,draw}] (n3) at (4,-6) {\small{3}};
  \node[fill=white,style={circle,draw}] (n1) at (-4,0) {\small{1}};
    \path (n3) [->,dashed] edge[] node {} (n1);
\path
     (n3) [->,dashed] edge[loop left=90] node {} (n3)
     (n3) [->,dashed] edge[loop right=90] node {} (n3);
 \end{tikzpicture} 
  \\
 \hline
 \hline 
$(d)$
\begin{tikzpicture}
 [scale=.15,auto=left, node distance=1.5cm, 
 ]
 \node[fill=white,style={circle,draw}] (n2) at (4,6) {\small{2}};
 \node[fill=black!30,style={circle,draw}] (n3) at (4,-6) {\small{3}};
  \node[fill=white,style={circle,draw}] (n1) at (-4,0) {\small{1}};
     \path (n3) [->,dashed] edge[bend left=10] node {} (n1);
     \path (n3) [->,dashed] edge[bend right=10] node {} (n1); 
\path
     (n3) [->,dashed] edge[loop right=90] node {} (n3);
 \end{tikzpicture}  &
 $(e)$ 
\begin{tikzpicture}
 [scale=.15,auto=left, node distance=1.5cm, 
 ]
 \node[fill=white,style={circle,draw}] (n2) at (4,6) {\small{2}};
 \node[fill=black!30,style={circle,draw}] (n3) at (4,-6) {\small{3}};
  \node[fill=white,style={circle,draw}] (n1) at (-4,0) {\small{1}};
    \path (n3) [->,dashed] edge[] node {} (n2);
\path
     (n3) [->,dashed] edge[loop left=90] node {} (n3)
     (n3) [->,dashed] edge[loop right=90] node {} (n3);
 \end{tikzpicture}  &
$(f)$  
\begin{tikzpicture}
 [scale=.15,auto=left, node distance=1.5cm, 
 ]
 \node[fill=white,style={circle,draw}] (n2) at (4,6) {\small{2}};
 \node[fill=black!30,style={circle,draw}] (n3) at (4,-6) {\small{3}};
  \node[fill=white,style={circle,draw}] (n1) at (-4,0) {\small{1}};
     \path (n3) [->,dashed] edge[bend left=10] node {} (n2);
     \path (n3) [->,dashed] edge[bend right=10] node {} (n2); 
\path
     (n3) [->,dashed] edge[loop right=90] node {} (n3);
 \end{tikzpicture} \\
  & &   \\
 \hline
 \end{tabular}}
 \end{center}
 \caption{The $3$-node networks 
 in which nodes $1, 2$ are excitatory and node $3$ is inhibitory, 
 only two nodes receive arrows, and
one node receives one arrow and another node receives two arrows. The arrows are 
 inhibitory and are outputs from node $3$.}
 \label{f:123_arrow_I}
 \end{figure}
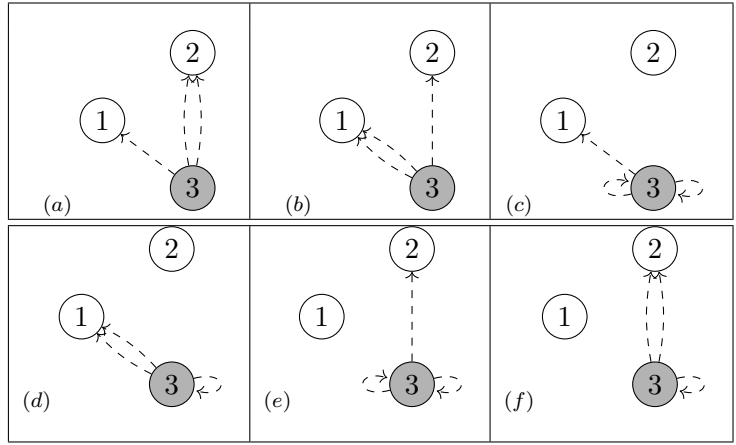

\begin{lemma}\label{lemma:123_strategy_sub_inhibition}
Let $\mathcal{G}$ be a connected $3$-node REI network of valence $2$, with node set $N^E \cup N^I$ where $N^E = \{1,2\}$,  $N^I = \{ 3\}$, arrow-types $A^E$, $A^I$ 
and satisfying 
\eqref{E:C}. 
 Then the subnetwork of $\mathcal{G}$ containing only arrow  type $A^I$ is one of the networks in {\rm Figure}~{\rm \ref{f:123_arrow_I}}.
\end{lemma}
\begin{proof}
The network $\mathcal{G}$ is REI and node $3$ is the only one of type $N^I$. By \eqref{E:C}, only two nodes receive arrows of type $A^I$ from node $3$. Moreover, one receives one arrow and the other two arrows.
\end{proof}

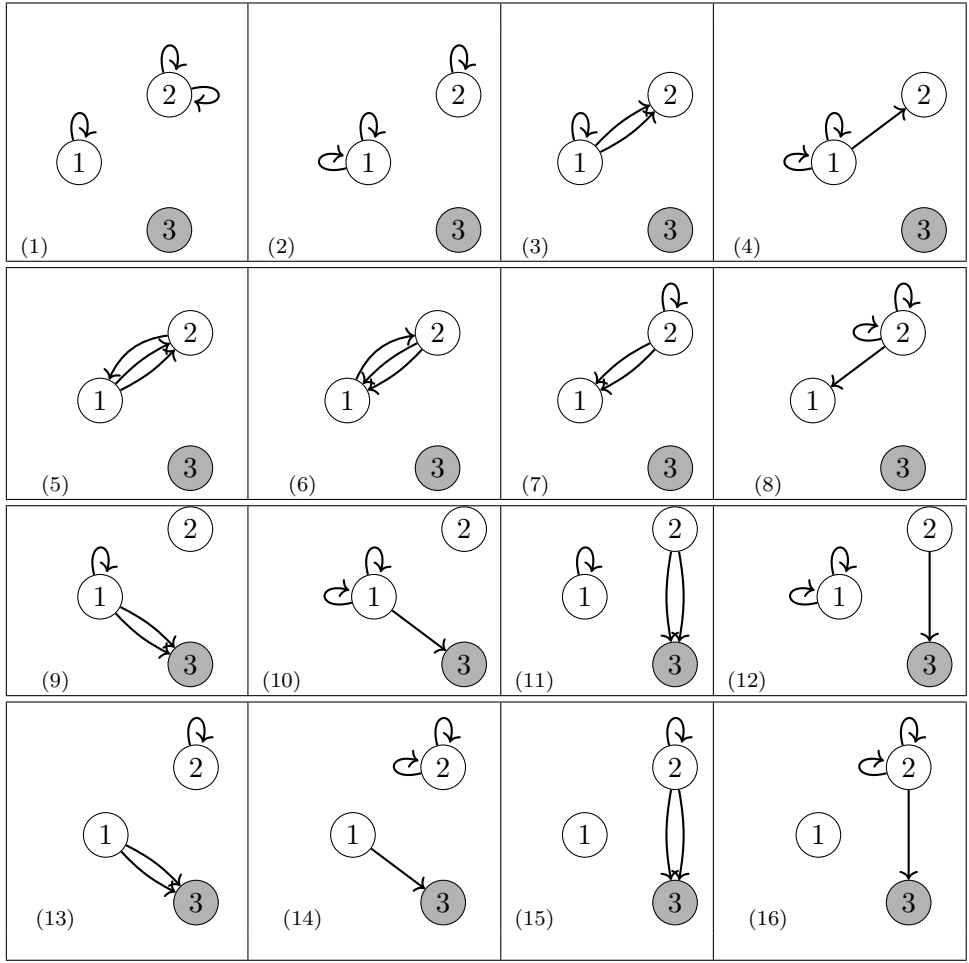
\begin{figure}
\begin{center}
{\tiny 
\begin{tabular}{|c|c|c|c|}
\hline 
 & & &  \\
 $(1)$ 
\begin{tikzpicture}
 [scale=.15,auto=left, node distance=1.5cm, 
 ]
 \node[fill=white,style={circle,draw}] (n2) at (4,6) {\small{2}};
 \node[fill=black!30,style={circle,draw}] (n3) at (4,-6) {\small{3}};
  \node[fill=white,style={circle,draw}] (n1) at (-4,0) {\small{1}};
   \path    (n1) [->,thick] edge[loop above=90] node {} (n1);
   \path    (n2) [->,thick] edge[loop above=90] node {} (n2);
   \path    (n2) [->,thick] edge[loop right=90] node {} (n2);
 \end{tikzpicture} &
$(2)$  
\begin{tikzpicture}
 [scale=.15,auto=left, node distance=1.5cm, 
 ]
 \node[fill=white,style={circle,draw}] (n2) at (4,6) {\small{2}};
 \node[fill=black!30,style={circle,draw}] (n3) at (4,-6) {\small{3}};
  \node[fill=white,style={circle,draw}] (n1) at (-4,0) {\small{1}};
   \path    (n1) [->,thick] edge[loop above=90] node {} (n1);
   \path    (n2) [->,thick] edge[loop above=90] node {} (n2);
   \path    (n1) [->,thick] edge[loop left=90] node {} (n1);
 \end{tikzpicture} 
 &
$(3)$
\begin{tikzpicture}
 [scale=.15,auto=left, node distance=1.5cm, 
 ]
 \node[fill=white,style={circle,draw}] (n2) at (4,6) {\small{2}};
 \node[fill=black!30,style={circle,draw}] (n3) at (4,-6) {\small{3}};
  \node[fill=white,style={circle,draw}] (n1) at (-4,0) {\small{1}};
   \path    (n1) [->,thick] edge[loop above=90] node {} (n1);
   \path    (n1) [->,thick] edge[bend left=10] node {} (n2);
   \path    (n1) [->,thick] edge[bend right=10] node {} (n2);
 \end{tikzpicture}
 &
$(4)$
\begin{tikzpicture}
 [scale=.15,auto=left, node distance=1.5cm, 
 ]
 \node[fill=white,style={circle,draw}] (n2) at (4,6) {\small{2}};
 \node[fill=black!30,style={circle,draw}] (n3) at (4,-6) {\small{3}};
  \node[fill=white,style={circle,draw}] (n1) at (-4,0) {\small{1}};
   \path    (n1) [->,thick] edge[loop above=90] node {} (n1);
   \path    (n1) [->,thick] edge[] node {} (n2);
   \path    (n1) [->,thick] edge[loop left=90] node {} (n1);
 \end{tikzpicture} 
  \\
 \hline
 \hline 
$(5)$
\begin{tikzpicture}
 [scale=.15,auto=left, node distance=1.5cm, 
 ]
 \node[fill=white,style={circle,draw}] (n2) at (4,6) {\small{2}};
 \node[fill=black!30,style={circle,draw}] (n3) at (4,-6) {\small{3}};
  \node[fill=white,style={circle,draw}] (n1) at (-4,0) {\small{1}};
   \path    (n2) [->,thick] edge[bend right=30] node {} (n1);
   \path    (n1) [->,thick] edge[bend left=10] node {} (n2);
   \path    (n1) [->,thick] edge[bend right=10] node {} (n2);
 \end{tikzpicture}
  &
 $(6)$ 
\begin{tikzpicture}
 [scale=.15,auto=left, node distance=1.5cm, 
 ]
 \node[fill=white,style={circle,draw}] (n2) at (4,6) {\small{2}};
 \node[fill=black!30,style={circle,draw}] (n3) at (4,-6) {\small{3}};
  \node[fill=white,style={circle,draw}] (n1) at (-4,0) {\small{1}};
   \path    (n1) [->,thick] edge[bend left=30] node {} (n2);
   \path    (n2) [->,thick] edge[bend left=10] node {} (n1);
   \path    (n2) [->,thick] edge[bend right=10] node {} (n1);
 \end{tikzpicture}
  &
$(7)$  
\begin{tikzpicture}
 [scale=.15,auto=left, node distance=1.5cm, 
 ]
 \node[fill=white,style={circle,draw}] (n2) at (4,6) {\small{2}};
 \node[fill=black!30,style={circle,draw}] (n3) at (4,-6) {\small{3}};
  \node[fill=white,style={circle,draw}] (n1) at (-4,0) {\small{1}};
   \path    (n2) [->,thick] edge[loop above=90] node {} (n2);
   \path    (n2) [->,thick] edge[bend left=10] node {} (n1);
   \path    (n2) [->,thick] edge[bend right=10] node {} (n1);
 \end{tikzpicture}
&
$(8)$ 
\begin{tikzpicture}
 [scale=.15,auto=left, node distance=1.5cm, 
 ]
 \node[fill=white,style={circle,draw}] (n2) at (4,6) {\small{2}};
 \node[fill=black!30,style={circle,draw}] (n3) at (4,-6) {\small{3}};
  \node[fill=white,style={circle,draw}] (n1) at (-4,0) {\small{1}};
   \path    (n2) [->,thick] edge[loop above=90] node {} (n2);
   \path    (n2) [->,thick] edge[] node {} (n1);
   \path    (n2) [->,thick] edge[loop left=90] node {} (n2);
 \end{tikzpicture} 
  \\
 \hline
 \hline 
$(9)$
\begin{tikzpicture}
 [scale=.15,auto=left, node distance=1.5cm, 
 ]
 \node[fill=white,style={circle,draw}] (n2) at (4,6) {\small{2}};
 \node[fill=black!30,style={circle,draw}] (n3) at (4,-6) {\small{3}};
  \node[fill=white,style={circle,draw}] (n1) at (-4,0) {\small{1}};
   \path    (n1) [->,thick] edge[loop above=90] node {} (n1);
   \path    (n1) [->,thick] edge[bend left=10] node {} (n3);
   \path    (n1) [->,thick] edge[bend right=10] node {} (n3);
 \end{tikzpicture}
&
$(10)$
\begin{tikzpicture}
 [scale=.15,auto=left, node distance=1.5cm, 
 ]
 \node[fill=white,style={circle,draw}] (n2) at (4,6) {\small{2}};
 \node[fill=black!30,style={circle,draw}] (n3) at (4,-6) {\small{3}};
  \node[fill=white,style={circle,draw}] (n1) at (-4,0) {\small{1}};
   \path    (n1) [->,thick] edge[loop above=90] node {} (n1);
   \path    (n1) [->,thick] edge[loop left=90] node {} (n1); 
   \path    (n1) [->,thick] edge[] node {} (n3);
 \end{tikzpicture}
&
$(11)$
\begin{tikzpicture}
 [scale=.15,auto=left, node distance=1.5cm, 
 ]
 \node[fill=white,style={circle,draw}] (n2) at (4,6) {\small{2}};
 \node[fill=black!30,style={circle,draw}] (n3) at (4,-6) {\small{3}};
  \node[fill=white,style={circle,draw}] (n1) at (-4,0) {\small{1}};
   \path    (n1) [->,thick] edge[loop above=90] node {} (n1);
   \path    (n2) [->,thick] edge[bend left=10] node {} (n3);
   \path    (n2) [->,thick] edge[bend right=10] node {} (n3);
 \end{tikzpicture}
&
$(12)$
\begin{tikzpicture}
 [scale=.15,auto=left, node distance=1.5cm, 
 ]
 \node[fill=white,style={circle,draw}] (n2) at (4,6) {\small{2}};
 \node[fill=black!30,style={circle,draw}] (n3) at (4,-6) {\small{3}};
  \node[fill=white,style={circle,draw}] (n1) at (-4,0) {\small{1}};
   \path    (n1) [->,thick] edge[loop above=90] node {} (n1);
   \path    (n1) [->,thick] edge[loop left=90] node {} (n1); 
   \path    (n2) [->,thick] edge[] node {} (n3);
 \end{tikzpicture}
 \\
 \hline
 \hline 
$(13)$
\begin{tikzpicture}
 [scale=.15,auto=left, node distance=1.5cm, 
 ]
 \node[fill=white,style={circle,draw}] (n2) at (4,6) {\small{2}};
 \node[fill=black!30,style={circle,draw}] (n3) at (4,-6) {\small{3}};
  \node[fill=white,style={circle,draw}] (n1) at (-4,0) {\small{1}};
   \path    (n2) [->,thick] edge[loop above=90] node {} (n2);
   \path    (n1) [->,thick] edge[bend left=10] node {} (n3);
   \path    (n1) [->,thick] edge[bend right=10] node {} (n3);
 \end{tikzpicture}
&
$(14)$
\begin{tikzpicture}
 [scale=.15,auto=left, node distance=1.5cm, 
 ]
 \node[fill=white,style={circle,draw}] (n2) at (4,6) {\small{2}};
 \node[fill=black!30,style={circle,draw}] (n3) at (4,-6) {\small{3}};
  \node[fill=white,style={circle,draw}] (n1) at (-4,0) {\small{1}};
   \path    (n2) [->,thick] edge[loop above=90] node {} (n2);
   \path    (n2) [->,thick] edge[loop left=90] node {} (n2); 
   \path    (n1) [->,thick] edge[] node {} (n3);
 \end{tikzpicture}
&
$(15)$
\begin{tikzpicture}
 [scale=.15,auto=left, node distance=1.5cm, 
 ]
 \node[fill=white,style={circle,draw}] (n2) at (4,6) {\small{2}};
 \node[fill=black!30,style={circle,draw}] (n3) at (4,-6) {\small{3}};
  \node[fill=white,style={circle,draw}] (n1) at (-4,0) {\small{1}};
   \path    (n2) [->,thick] edge[loop above=90] node {} (n2);
   \path    (n2) [->,thick] edge[bend left=10] node {} (n3);
   \path    (n2) [->,thick] edge[bend right=10] node {} (n3);
 \end{tikzpicture}
&
$(16)$
\begin{tikzpicture}
 [scale=.15,auto=left, node distance=1.5cm, 
 ]
 \node[fill=white,style={circle,draw}] (n2) at (4,6) {\small{2}};
 \node[fill=black!30,style={circle,draw}] (n3) at (4,-6) {\small{3}};
  \node[fill=white,style={circle,draw}] (n1) at (-4,0) {\small{1}};
   \path    (n2) [->,thick] edge[loop above=90] node {} (n2);
   \path    (n2) [->,thick] edge[loop left=90] node {} (n2); 
   \path    (n2) [->,thick] edge[] node {} (n3);
 \end{tikzpicture}
 \\
  & & &  \\
 \hline
 \end{tabular}}
 \end{center}
 \caption{The $3$-node networks where nodes $1, 2$ are excitatory and node $3$ is inhibitory, only two nodes receive arrows, and
one node receives one arrow and another node receives two arrows. The arrows are of type $A^E$ and are outputs from nodes $1$ and/or $2$.}
 \label{f:123_arrow_E}
 \end{figure}

\begin{lemma}\label{lemma:123_strategy_sub_activation}
Let $\mathcal{G}$ be a connected $3$-node REI network of valence $2$, with 
$N^E = \{1,2\}$,  $N^I = \{ 3\}$, two distinct arrow-types $A^E$, $A^I$,
and satisfying \eqref{E:C}. 
Then the subnetwork of $\mathcal{G}$ containing only arrow-type $A^E$ is one of the networks in Figure~{\rm \ref{f:123_arrow_E}}.
\end{lemma}
\begin{proof}
The network $\mathcal{G}$ is REI and nodes $1$ and $2$ are those of type $N^E$. By \eqref{E:C} only two nodes receive arrows of type $A^E$ from nodes $1$ and $2$. Moreover, one receives one arrow and the other receives two arrows.
\end{proof}

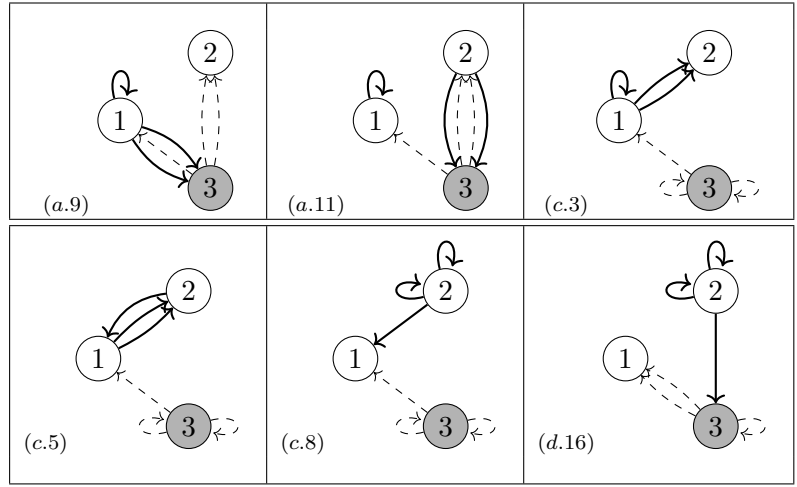
\begin{figure}
\begin{center}
{\tiny 
\begin{tabular}{|c|c|c|}
\hline 
 & &  \\
 $(a.9)$ 
\begin{tikzpicture}
 [scale=.15,auto=left, node distance=1.5cm, 
 ]
 \node[fill=white,style={circle,draw}] (n2) at (4,6) {\small{2}};
 \node[fill=black!30,style={circle,draw}] (n3) at (4,-6) {\small{3}};
  \node[fill=white,style={circle,draw}] (n1) at (-4,0) {\small{1}};
    \path (n3) [->,dashed] edge[] node {} (n1);
     \path (n3) [->,dashed] edge[bend left=10] node {} (n2);
     \path (n3) [->,dashed] edge[bend right=10] node {} (n2);
   \path    (n1) [->,thick] edge[loop above=90] node {} (n1);
   \path    (n1) [->,thick] edge[bend left=20] node {} (n3);
   \path    (n1) [->,thick] edge[bend right=20] node {} (n3);
 \end{tikzpicture}
 &
$(a.11)$  
\begin{tikzpicture}
 [scale=.15,auto=left, node distance=1.5cm, 
 ]
 \node[fill=white,style={circle,draw}] (n2) at (4,6) {\small{2}};
 \node[fill=black!30,style={circle,draw}] (n3) at (4,-6) {\small{3}};
  \node[fill=white,style={circle,draw}] (n1) at (-4,0) {\small{1}};
    \path (n3) [->,dashed] edge[] node {} (n1);
     \path (n3) [->,dashed] edge[bend left=10] node {} (n2);
     \path (n3) [->,dashed] edge[bend right=10] node {} (n2);
   \path    (n1) [->,thick] edge[loop above=90] node {} (n1);
   \path    (n2) [->,thick] edge[bend left=25] node {} (n3);
   \path    (n2) [->,thick] edge[bend right=25] node {} (n3); 
 \end{tikzpicture} 
 &
$(c.3)$
\begin{tikzpicture}
 [scale=.15,auto=left, node distance=1.5cm, 
 ]
 \node[fill=white,style={circle,draw}] (n2) at (4,6) {\small{2}};
 \node[fill=black!30,style={circle,draw}] (n3) at (4,-6) {\small{3}};
  \node[fill=white,style={circle,draw}] (n1) at (-4,0) {\small{1}};
    \path (n3) [->,dashed] edge[] node {} (n1);
\path
     (n3) [->,dashed] edge[loop left=90] node {} (n3)
     (n3) [->,dashed] edge[loop right=90] node {} (n3);
   \path    (n1) [->,thick] edge[loop above=90] node {} (n1);
   \path    (n1) [->,thick] edge[bend left=10] node {} (n2);
   \path    (n1) [->,thick] edge[bend right=10] node {} (n2);
 \end{tikzpicture} 
  \\
 \hline
 \hline 
$(c.5)$
\begin{tikzpicture}
 [scale=.15,auto=left, node distance=1.5cm, 
 ]
 \node[fill=white,style={circle,draw}] (n2) at (4,6) {\small{2}};
 \node[fill=black!30,style={circle,draw}] (n3) at (4,-6) {\small{3}};
  \node[fill=white,style={circle,draw}] (n1) at (-4,0) {\small{1}};
    \path (n3) [->,dashed] edge[] node {} (n1);
\path
     (n3) [->,dashed] edge[loop left=90] node {} (n3)
     (n3) [->,dashed] edge[loop right=90] node {} (n3);
   \path    (n2) [->,thick] edge[bend right=30] node {} (n1);
   \path    (n1) [->,thick] edge[bend left=10] node {} (n2);
   \path    (n1) [->,thick] edge[bend right=10] node {} (n2);
 \end{tikzpicture} 
&
 $(c.8)$
\begin{tikzpicture}
 [scale=.15,auto=left, node distance=1.5cm, 
 ]
 \node[fill=white,style={circle,draw}] (n2) at (4,6) {\small{2}};
 \node[fill=black!30,style={circle,draw}] (n3) at (4,-6) {\small{3}};
  \node[fill=white,style={circle,draw}] (n1) at (-4,0) {\small{1}};
    \path (n3) [->,dashed] edge[] node {} (n1);
\path
     (n3) [->,dashed] edge[loop left=90] node {} (n3)
     (n3) [->,dashed] edge[loop right=90] node {} (n3);
   \path    (n2) [->,thick] edge[loop above=90] node {} (n2);
   \path    (n2) [->,thick] edge[] node {} (n1);
   \path    (n2) [->,thick] edge[loop left=90] node {} (n2);
 \end{tikzpicture} 
&
$(d.16)$
\begin{tikzpicture}
 [scale=.15,auto=left, node distance=1.5cm, 
 ]
 \node[fill=white,style={circle,draw}] (n2) at (4,6) {\small{2}};
 \node[fill=black!30,style={circle,draw}] (n3) at (4,-6) {\small{3}};
  \node[fill=white,style={circle,draw}] (n1) at (-4,0) {\small{1}};
     \path (n3) [->,dashed] edge[bend left=10] node {} (n1);
     \path (n3) [->,dashed] edge[bend right=10] node {} (n1); 
   \path    (n2) [->,thick] edge[loop above=90] node {} (n2);
   \path    (n2) [->,thick] edge[loop left=90] node {} (n2); 
   \path    (n2) [->,thick] edge[] node {} (n3);
\path
     (n3) [->,dashed] edge[loop right=90] node {} (n3);
 \end{tikzpicture} \\
  & &   \\
 \hline
 \end{tabular}}
 \end{center}
 \caption{The connected $3$-node REI networks with valence $2$, in which
 nodes $1, 2$ are excitatory and node $3$ is inhibitory, having  
two arrow-types $A^E$ and $A^I$, input equivalence relation $\sim_I\, = \left\{ \{ 1\}, \{ 2\}, \{ 3\} \right\}$,
and satisfying \eqref{E:C}.
}
 \label{fig:123_3NCNREIV2}
 \end{figure}

\begin{prop} \label{prop:non_hom_REI123}
The set of connected $3$-node REI networks of valence $2$ with input equivalence relation $\sim_I = \left\{ \{ 1\}, \{ 2\}, \{ 3\} \right\}$ 
and satisfying \eqref{E:C} 
is listed in  {\rm Figure}~{\rm \ref{fig:123_3NCNREIV2}}.
\end{prop}
\begin{proof}
 Up to duality, REI networks have nodes $1$ and $2$  of type $N^E$ and node $3$ of type $N^I$. If $\mathcal{G}$ is a connected $3$-node REI network with valence $2$, input equivalence relation $\sim_I = \left\{ \{ 1\}, \{ 2\}, \{ 3\} \right\}$ 
 and satisfying \eqref{E:C},
 then the subnetwork containing only arrow-type $A^I$ is one of the networks listed in Lemma~\ref{lemma:123_strategy_sub_inhibition}, and the subnetwork of $\mathcal{G}$ containing only arrow  type $A^E$ is one of the networks listed in Lemma~\ref{lemma:123_strategy_sub_activation}. We obtain the networks in Figure~\ref{fig:123_3NCNREIV2}.
\end{proof}

\section{Conclusions}

Motivated by the growing interest in network motifs and their functionality in biological networks, and following the work in \cite{ADS24}, we give a characterization of the connected 3-node restricted excitatory-inhibitory networks.
Our classifications are up to renumbering of the nodes and duality -- switch nodes and arrows types from `excitatory' to 'inhibitory', and vice versa.
Although there is an infinity of connected 3-node restricted excitatory-inhibitory networks, when we restrict to networks with valence less or equal to $2$ -- each node receives at most $2$ inputs -- we get a finite number. 
Taking our characterization further, we also list those networks with valence exactly equal to 2, under different conditions on the input arrows of the 3 nodes, ranging from all nodes receiving an arrow of each type to all having non-isomorphic input sets.
Both, for all connected 3-node restricted excitatory-inhibitory networks and those with valence less or equal to $2$, we give their characterization under ODE-equivalence. Moreover, we give a minimal representative for each ODE-class.

The next step for future work in our systematic study is to explore the dynamics, in particular the bifurcations, of these 3-node restricted excitatory-inhibitory motifs. This will be in line to what is done in \cite{MRS24} for six particular motifs that occur as functional building blocks in gene regulatory networks, where the state of each
gene is modeled in terms of two variables: mRNA and protein concentration. The study in \cite{MRS24} explores the patterns of synchrony (fibration symmetries) of the motifs and considers both all possible network admissible models as well as special specializations to simple models based on Hill functions and linear degradation.

\vspace{5mm}

\noindent {\bf Acknowledgments} \\
MA and AD were partially supported by CMUP, member of LASI, which is financed by national funds through FCT -- Funda\c c\~ao para a Ci\^encia e a Tecnologia, I.P., under the projects with reference UIDB/00144/2020 and UIDP/00144/2020.

\end{document}